\documentclass[11pt]{article}

\usepackage[T1]{fontenc}
\usepackage[english]{babel}
\usepackage[letterpaper]{geometry}
\usepackage{helvet}

\usepackage{amsmath}
\usepackage{amsbsy}
\usepackage{amssymb}
\usepackage{nccmath}
\usepackage{amsthm}

\usepackage{graphicx}
\usepackage{float}
\usepackage{wrapfig}
\usepackage{sidecap}
\usepackage{adjustbox}
\usepackage{subcaption}
\usepackage{pdfpages}

\usepackage{array}
\usepackage{multirow}
\usepackage{tabularx}
\usepackage{booktabs}
\usepackage{longtable}
\usepackage[table]{xcolor}
\usepackage{colortbl}

\usepackage{caption}
\usepackage{tocloft}

\usepackage{setspace}
\usepackage{enumitem}

\usepackage{textgreek}
\usepackage{breakcites}
\usepackage{cite}

\PassOptionsToPackage{normalem}{ulem}
\usepackage{ulem}

\usepackage{algorithm}
\usepackage{algpseudocode}

\usepackage{tikz}
\usepackage{tikz-cd}
\usepackage{pgfplots}
\usepackage{pgfgantt}
\usepackage{pgf-pie}

\pgfplotsset{compat=1.18}

\usetikzlibrary{
    3d,
    shapes.geometric,
    arrows,
    calc,
    positioning,
    decorations.pathmorphing,
    fit,
    patterns,
    decorations.pathreplacing,
    through,
    mindmap,
    shadows,
    backgrounds
}

\usepackage{framed}
\usepackage{tcolorbox}

\usepackage{units}
\usepackage{lipsum}

\usepackage{hyperref}
\usepackage[nameinlink,noabbrev]{cleveref}

\makeatletter

\renewcommand{\section}{%
    \@startsection{section}{1}{\z@}%
    {-3.0ex \@plus -1ex \@minus -.2ex}%
    {1.5ex \@plus .2ex}%
    {\normalfont\large\bfseries}%
}

\renewcommand{\subsection}{%
    \@startsection{subsection}{2}{\z@}%
    {-2.5ex \@plus -1ex \@minus -.2ex}%
    {1.0ex \@plus .2ex}%
    {\normalfont\normalsize\bfseries}%
}

\renewcommand{\subsubsection}{%
    \@startsection{subsubsection}{3}{\z@}%
    {-2.0ex \@plus -1ex \@minus -.2ex}%
    {0.8ex \@plus .2ex}%
    {\normalfont\normalsize\bfseries\itshape}%
}

\renewcommand{\@biblabel}[1]{[#1]}

\makeatother

\setlist{
    nosep,
    topsep=2pt,
    partopsep=0pt,
    parsep=0pt,
    itemsep=1pt
}

\floatstyle{plaintop}
\restylefloat{algorithm}

\makeatletter
\@ifundefined{creflastconjunction}
    {}
    {}
\makeatother

\newenvironment{keywords}
{%
    \par\noindent
    \textbf{Keywords: }\ignorespaces
}
{%
    \par
}

\newtheorem{theorem}{Theorem}

\newtheorem{lemma}{Lemma}

\newtheorem{assumption}{Assumption}

\newtheorem{corollary}{Corollary}

\theoremstyle{definition}
\newtheorem{definition}{Definition}

\theoremstyle{remark}
\newtheorem{remark}{Remark}

\theoremstyle{definition}
\newtheorem{example}{Example}

\begin{document}

\title{LEvITAS: An Operator-Theoretic Variational model for Generative Modeling and Structured Synthesis}
\author{
  Lokendra S. Thakur\textsuperscript{1,2*}, 
  James Dale Berry\textsuperscript{2}
}

\date{
\raggedright
\textsuperscript{1} Innovation Engine Lab (iEL), MDtRI, Hyderabad, INDIA.\\ 
\textsuperscript{2} Department of Neurology, Massachusetts General Hospital, Harvard Medical School, Massachusetts, USA.\\
*Corresponding Author: lthakur@mdtri.org, Copyright \textcopyright{} with MDtRI Foundation,
\href{https://www.mdtri.org}{Visit MDtRI}
}

\maketitle

\begin{abstract}
We introduce LEvITAS (\textbf{L}atent \textbf{S}pace 
\textbf{E}le\textbf{v}ation and \textbf{I}terative \textbf{T}opological 
\textbf{A}pproximation \textbf{S}ystem), a generalized variational lifting 
model for structure-preserving probabilistic evolution in generative 
modeling. The formulation represents generative systems through compositions 
of measurable lifting--projection operators, structure-preserving refinement 
maps, and variational objectives defined on probability measures and their 
trajectories. An intermediate synthetic structure space is introduced between 
latent and observable representations, allowing geometric, temporal, 
semantic, and topological constraints to be imposed explicitly during 
probabilistic evolution. The resulting model connects operator-theoretic 
lifting, variational dynamics, Wasserstein transport, measure-valued 
evolution, and continuum probability flows within a common formulation. We 
establish stability and convergence properties of lifted trajectories and 
variational refinements, together with probability preservation, weak 
stability, transport-geometric structure, and continuum formulations 
including continuity-equation, Fokker--Planck, and Wasserstein gradient-flow 
dynamics. Several widely used generative architectures, including generative 
adversarial networks, variational autoencoders, diffusion models, normalizing 
flows, and energy-based models, arise as constrained reductions of the 
general LEvITAS formulation. We further evaluate the model across 
heterogeneous biomedical data modalities, including laparoscopic imaging, 
ultrasound phenotyping, multi-omic data, genomic--morphological 
representations, and multimodal clinical data. Across these settings, the 
empirical analyses provide evidence that our structure-preserving formulation 
retains dominant latent organization, cross-modal relationships, and 
distributional characteristics between observed and generated 
representations. These results support LEvITAS as a unified mathematical 
model for structure-preserving generative modeling and probabilistic 
evolution, connecting latent representations, measurable operators, 
variational refinement, optimal transport, and continuum dynamics within a 
common data-scientific formulation.
\end{abstract}


\section{Introduction}
Modern generative modeling encompasses several complementary paradigms, 
including adversarial learning~ 
\cite{goodfellow2014generative,arjovsky2017wasserstein},
variational inference~\cite{hoffman2013stochastic,kingma2014autoencoding,
blei2017variational}, stochastic diffusion 
models~\cite{sohldickstein2015deep,ho2020denoising,
song2021scorebased}, invertible flow-based models~
\cite{dinh2017density,kingma2018glow,kobyzev2021normalizing}, and continuous 
dynamical systems~\cite{chen2018neural,grathwohl2019ffjord,
kidger2020neural}. Despite architectural differences, these approaches share 
a common objective: the construction of transformations from latent 
probability spaces into structured observable phenomena. 
Most existing formulations model generation as a direct transformation from 
latent variables into observable sample space. However, many physical, 
biological, and engineered systems naturally involve intermediate structural 
representations prior to observation. Examples include latent physical states 
prior to measurement, semantic scene representations prior to rendering, and 
hidden dynamical trajectories prior to sensing.

Motivated by these considerations, we introduce LEvITAS, an operator-
theoretic variational model in which generation proceeds through:
(i) latent lifting, (ii) structured synthesis, (iii) observable projection, 
(iv) structural consistency enforcement, and (v) dynamical trajectory preservation.

The principal contribution of this work is the formulation of a generalized 
variational operator system, termed LEvITAS, from which several major 
generative model classes emerge as special cases under suitable structural 
and variational conditions. A central component and also one of the novelties 
of this formulation is an intermediate lifting space that establishes a 
common mathematical representation between the observational data space and 
the generated outputs.  
Observational data are mapped into this intermediate space, where 
admissibility constraints, regularity properties, and variational energy 
conditions can be imposed and analyzed through controlled refinement 
mechanisms. For example, when the lifted representations consist of 
probability densities, constraints such as non-negativity, unit probability 
mass, and bounded variational energy restrict the admissible states and 
consequently the outputs attainable through the observation operator. The 
refined lifted representations are subsequently mapped back into the 
observational data space to produce synthetic datasets. Furthermore, the 
framework provides a basis for quantitatively assessing the distributional 
fidelity of generated datasets relative to observed data through Wasserstein-
based comparisons, while recognizing that distributional distance alone does 
not fully characterize structural fidelity. The development of dedicated 
structural quality-assessment metrics, particularly for evaluating clinically 
meaningful properties of synthetic datasets, is reserved for the clinical 
extension of LEvITAS. Through this construction, LEvITAS provides a unified 
operator-theoretic foundation connecting generative model architectures with 
measure-theoretic dynamics, Wasserstein geometry, variational gradient flows, 
and their continuum PDE formulations.

\section{Main Results}
\label{sec:main}

The principal analytical contributions of the LEvITAS model are
developed through six complementary components, each addressing a
distinct but interconnected aspect of its mathematical structure.
Together, these components establish a unified analytical foundation
for LEvITAS, connecting operator-theoretic representations,
trajectory-level dynamics, variational refinement, probability
evolution, optimal transport, generative-model reductions, and
continuum partial differential equation (PDE) formulations.

The first component develops the \emph{dynamical theory} of LEvITAS,
establishing the existence and continuity of lifted trajectories,
their uniform stability under latent perturbations, weak convergence
of the associated path measures, and the existence of dynamically
consistent variational minimizers. These results provide the
analytical foundation for the evolution of lifted representations.

The second component establishes the \emph{variational refinement
theory}, addressing the local existence and uniqueness of refinement
trajectories, monotone energy dissipation, characterization of
stationary equilibrium states, exponential stability under strong
convexity, preservation of topological invariants, and monotone
convergence of discrete refinement iterations. This component
characterizes the continuous and discrete behavior of the LEvITAS
refinement mechanism.

The third component develops the \emph{operator-theoretic reduction
theory}, demonstrating how major generative-model architectures,
including generative adversarial networks (GANs), variational
autoencoders (VAEs), diffusion systems, normalizing flows, and
energy-based models, arise as particular realizations of the
general LEvITAS formulation. A hierarchical containment result
further establishes the structural relationships among these
architectures within the unified framework.

The fourth component establishes the \emph{measure-theoretic and
probabilistic foundations} of LEvITAS. It addresses probability
preservation under measurable pushforward mappings, weak stability
of induced probability measures, Wasserstein metric structure,
compactness through Prokhorov tightness, existence of synthetic
probability measures for stochastic liftings, measure disintegration
on Polish spaces, and mass preservation along measure-valued
refinement flows. These results provide the probabilistic
consistency required for the stochastic and measure-valued
formulations.

The fifth component develops the \emph{optimal transport and
Wasserstein-geometric theory}. It establishes the existence of
Monge-type and Kantorovich transport minimizers, the geodesic
structure of Wasserstein probability spaces, continuity of
displacement interpolation flows, the dynamical characterization
of transport through the Benamou--Brenier formulation, and energy
dissipation along Wasserstein gradient flows. This component
provides the geometric and variational framework connecting
probability distributions through optimal transport.

The sixth component establishes the \emph{PDE and continuum
formulation} of LEvITAS, extending the preceding operator-theoretic,
dynamical, measure-theoretic, and variational structures to
continuous evolution systems. In particular, the continuity
equation, continuum gradient flow, and Fokker--Planck dynamics
provide deterministic and stochastic continuum realizations.
The associated entropy dissipation law, weak-solution existence
theory, and equilibrium characterization establish the analytical
foundations of these continuum systems and connect stationary
continuum states with the underlying variational structure.

Collectively, these six components form a coherent analytical
hierarchy in which lifted representations, refinement dynamics,
probability measures, optimal transport, generative architectures,
and continuum evolution equations are interconnected through
common operator-theoretic, variational, geometric, and
probabilistic principles. The resulting model provides a
unified perspective on the existence, stability, convergence,
structural preservation, and equilibrium properties of LEvITAS
across both finite-dimensional and continuum formulations.

\subsection{Main Dynamical Results}
\label{subsec:main_dynamical_results}

The principal analytical results of this section establish the existence, stability, convergence, and variational consistency of the lifted LEvITAS dynamics. Specifically, we establish the existence of continuous lifted trajectories (Theorem~\ref{thm:existence_continuous_lifted_paths}), uniform stability with respect to perturbations in the latent representation (Theorem~\ref{thm:uniform_stability_lifted_paths}), and weak convergence of the associated lifted path measures (Theorem~\ref{thm:weak_stability_lifted_path_measures}). We further prove the existence of dynamically consistent variational minimizers (Theorem~\ref{thm:existence_dynamically_consistent_liftings}).  
Together, these results provide a unified analytical foundation for the lifted dynamical framework. The existence and stability results characterize the well-posedness and robustness of individual lifted trajectories, while weak convergence extends these properties to the level of path distributions. The variational existence result completes the analysis by establishing the compatibility of the lifted dynamics with the underlying variational structure of LEvITAS.

\subsection{Main Refinement Results}
\label{subsec:main_refinement_results}

The principal analytical results of this section establish the well-posedness, variational behavior, stability, structural preservation, and convergence properties of the LEvITAS refinement mechanism. Specifically, we prove the local existence and uniqueness of refinement trajectories (Theorem~\ref{thm:local_existence_refinement_flows}) and establish a monotone energy dissipation law along the resulting continuous refinement flows (Theorem~\ref{thm:energy_dissipation_law}). We further characterize the stationary equilibrium states of the refinement dynamics (Theorem~\ref{thm:characterization_equilibrium_states}) and demonstrate exponential stability under strong convexity assumptions (Theorem~\ref{thm:exponential_stability}). In addition, we establish the preservation of topological invariants under continuous refinement (Theorem~\ref{thm:topological_preservation_continuous_refinement}) and prove monotone convergence of the associated discrete refinement iterations (Theorem~\ref{thm:monotone_descent_discrete_refinement}). 
Together, these results provide a unified analytical foundation for the LEvITAS refinement framework across both continuous and discrete formulations. The existence and uniqueness result establishes the local well-posedness of the refinement evolution, while energy dissipation and equilibrium characterization describe its variational structure and stationary behavior. Exponential stability further quantifies convergence toward equilibrium under strong convexity, whereas topological preservation ensures that the prescribed structural invariants remain intact throughout continuous refinement. Finally, the discrete convergence result extends the descent analysis to iterative refinement schemes, connecting the continuous evolution theory with its computational realization.

\subsection{Main Reduction Results}
\label{subsec:main_reduction_results}

The principal analytical results of this section establish the general operator-theoretic reduction principle and its application to several major generative-model architectures. Specifically, we formulate a general reduction principle for identifying particular generative systems within the unified LEvITAS model (Theorem~\ref{thm:general_reduction_theorem}). Building on this principle, we establish the reduction of LEvITAS to generative adversarial networks (GANs) (Theorem~\ref{thm:gan_reduction}), variational autoencoders (VAEs) (Theorem~\ref{thm:vae_reduction}), diffusion systems (Theorem~\ref{thm:diffusion_reduction}), normalizing flows (Theorem~\ref{thm:normalizing_flow_reduction}), and energy-based models (Theorem~\ref{thm:energy_based_reduction}). We further establish a hierarchical containment result that organizes these generative architectures within the broader LEvITAS formulation (Theorem~\ref{thm:hierarchical_containment}). 
Together, these results provide a unified operator-theoretic perspective on the relationship between LEvITAS and established generative-model architectures. The general reduction principle supplies the mathematical foundation for recovering individual generative systems, while the architecture-specific results demonstrate how GANs, VAEs, diffusion models, normalizing flows, and energy-based models arise as particular realizations of the unified formulation. The hierarchical containment result further clarifies the structural relationships among these architectures and their placement within the broader LEvITAS model. Consequently, the reduction theory establishes a common analytical foundation for understanding diverse generative architectures through a unified operator-theoretic formulation.

\subsection{Main Measure-Theoretic Results}
\label{subsec:main_measure_theoretic_results}

The principal analytical results of this section establish the probabilistic consistency, measure-theoretic stability, and structural foundations of the LEvITAS operator framework. Specifically, we establish probability preservation under measurable pushforward mappings (Theorem~\ref{thm:probability_preservation}) and weak stability of the induced pushforward probability measures (Theorem~\ref{thm:stability_pushforward_measures}). We further characterize the metric structure of Wasserstein probability geometry (Theorem~\ref{thm:wasserstein_metric_structure}) and establish compactness properties through Prokhorov tightness (Theorem~\ref{thm:prokhorov_compactness}). For the stochastic setting, we prove the existence of induced synthetic probability measures associated with stochastic liftings (Theorem~\ref{thm:existence_induced_synthetic_measures}) and establish the existence of measure disintegration structures on Polish spaces (Theorem~\ref{thm:existence_disintegration}). Finally, we demonstrate mass preservation along measure-valued LEvITAS refinement flows (Theorem~\ref{thm:mass_preservation}). 
Together, these results provide a unified measure-theoretic foundation for the probabilistic formulation of LEvITAS. Probability preservation and weak stability establish the consistency and continuity of induced probability distributions under measurable transformations, while Wasserstein geometry and Prokhorov compactness provide the metric and topological structures needed for analyzing convergence and stability in probability spaces. The existence of induced synthetic measures and disintegration structures extends this foundation to stochastic liftings and conditional probabilistic representations. Finally, mass preservation ensures the probabilistic consistency of measure-valued refinement dynamics, connecting the measure-theoretic formulation with the broader variational and dynamical structure of LEvITAS.

\subsection{Main Optimal Transport Results}
\label{subsec:main_optimal_transport_results}

The principal analytical results of this section establish the existence, geometric structure, and variational dynamics of optimal transport within the LEvITAS model. Specifically, we establish the existence of Monge-type LEvITAS transport minimizers (Theorem~\ref{thm:existence_monge_levitas_transport}) and Kantorovich transport minimizers (Theorem~\ref{thm:existence_kantorovich_minimizers}). We further characterize the geodesic structure of Wasserstein probability spaces (Theorem~\ref{thm:geodesic_structure_wasserstein_space}) and establish the continuity of the associated displacement interpolation flows (Theorem~\ref{thm:continuity_displacement_flows}). In addition, we derive a dynamical characterization of optimal transport through the Benamou--Brenier formulation (Theorem~\ref{thm:benamou_brenier_formulation}) and establish energy dissipation along Wasserstein gradient flows (Theorem~\ref{thm:energy_dissipation_wasserstein_flows}). 
Together, these results provide a unified geometric and variational foundation for optimal transport within LEvITAS. The Monge and Kantorovich existence results establish the availability of optimal transport solutions under their respective formulations, while the geodesic structure of Wasserstein space provides the geometric framework for interpolation between probability measures. The continuity of displacement interpolation connects this geometry with continuous probability evolution, whereas the Benamou--Brenier formulation characterizes optimal transport through a dynamical variational principle. Finally, the Wasserstein gradient-flow result establishes the corresponding energy-dissipation structure, linking optimal transport geometry with the variational refinement and probabilistic evolution mechanisms of LEvITAS.

\subsection{Main PDE and Continuum Formulation Results}
\label{subsec:pde_connections_main_results}

The principal analytical contribution of this section is to establish the continuum realization of the LEvITAS operator-theoretic and variational framework. In particular, the continuity equation
\eqref{eq:levitas_continuity}
and continuum gradient flow
\eqref{eq:continuum_levitas}
extend the refinement dynamics developed in
Section~\ref{sec:iterative_topological_refinement}
to a continuous evolution setting. The associated Wasserstein-gradient
structure connects this formulation with the optimal transport geometry
established in
Section~\ref{sec:optimal_transport_wasserstein_geometry}.
Furthermore, the Fokker--Planck dynamics
\eqref{eq:fokker_planck_levitas}
provide a continuum stochastic realization of the probability evolution
and stochastic lifting structures developed in
Section~\ref{sec:measure_theoretic_foundations}. 
The analytical properties of these continuum systems are established
through several complementary results. The entropy dissipation law
(Theorem~\ref{thm:entropy_dissipation})
extends the variational energy-dissipation principles of the refinement
framework to the continuum setting. The weak solution formulation
(Definition~\ref{def:weak_solution})
and its corresponding existence result
(Theorem~\ref{thm:existence_weak_solutions})
provide the analytical foundation for continuum LEvITAS evolution.
Moreover, the equilibrium characterization
(Theorem~\ref{thm:equilibrium_characterization})
connects stationary continuum states with the variational minimization
structure of the generalized LEvITAS functional introduced in
Section~\ref{sec:unified_reduction_theory}. 
Together, these formulations and results establish the analytical
connection between discrete refinement mechanisms, measure-valued
probability evolution, Wasserstein transport geometry, stochastic
dynamics, and infinite-dimensional variational optimization.
The continuity equation describes the transport of probability mass,
the continuum gradient flow characterizes variational refinement, and
the Fokker--Planck equation provides the corresponding stochastic
evolution model. The weak-solution theory establishes the existence
framework for these continuum dynamics, while entropy dissipation and
equilibrium characterization describe their variational and
asymptotic properties.

\subsection{Unified analytical structure of LEvITAS}

\begin{wrapfigure}[27]{r}{0.5\textwidth}
\centering
\vspace{-0.5cm}
\begin{tikzpicture}[
    scale=0.58,
    transform shape
]

\node[
    draw,
    rectangle,
    minimum width=8cm,
    minimum height=1cm,
    align=center
] at (0,0) {
    \textbf{LEvITAS Unified Operator}
};

\node[
    draw,
    rectangle,
    minimum width=5cm,
    minimum height=2.5cm,
    align=center
] at (-3,-2.75) {
    \textbf{Dynamical Structure}\\
    Lifted trajectories\\
    Refinement trajectories\\
    Existence\\
    Stability\\
    Convergence
};

\node[
    draw,
    rectangle,
    minimum width=5cm,
    minimum height=2.5cm,
    align=center
] at (3,-2.75) {
    \textbf{Measure-Theoretic Structure}\\
    Pushforward measures\\
    Probability preservation\\
    Weak stability\\
    Disintegration\\
    Mass preservation
};

\node[
    draw,
    rectangle,
    minimum width=8cm,
    minimum height=2.5cm,
    align=center
] at (0,-5.75) {
    \textbf{Variational and Wasserstein--Transport Core}\\
    Energy functionals\\
    Optimal transport\\
    Wasserstein geometry\\
    Gradient flows
};

\node[
    draw,
    rectangle,
    minimum width=5cm,
    minimum height=2.5cm,
    align=center
] at (-3,-8.75) {
    \textbf{Generative Models}\\
    GAN\\
    VAE\\
    Diffusion\\
    Normalizing flows\\
    Energy-based models
};

\node[
    draw,
    rectangle,
    minimum width=5cm,
    minimum height=2.5cm,
    align=center
] at (3,-8.75) {
    \textbf{Continuum Formulation}\\
    Continuity equation\\
    Gradient flow\\
    Fokker--Planck\\
    Weak solutions\\
    Equilibrium
};

\node[
    draw,
    rectangle,
    minimum width=7cm,
    minimum height=1cm,
    align=center
] at (0,-11.5) {
    \textbf{Equilibrium, Convergence, and Stability}
};

\draw[->] (-1.2,-0.6) -- (-2.2,-1.3);
\draw[->] (1.2,-0.6) -- (2.2,-1.3);

\draw[->] (-3,-3.8) -- (-1.8,-4.5);
\draw[->] (3,-3.8) -- (1.8,-4.5);

\draw[->] (-2,-6.8) -- (-2.5,-7.2);
\draw[->] (2,-6.8) -- (2.5,-7.2);

\draw[->] (-3,-10.2) -- (-1.5,-10.8);
\draw[->] (3,-10.2) -- (1.5,-10.8);

\end{tikzpicture}

\caption{
Unified analytical structure of LEvITAS. The LEvITAS Unified Operator
provides the generalize operator-theoretic foundation for the dynamical and
measure-theoretic structures. These structures are coupled through the
variational and Wasserstein--transport core, which provides the
geometric and variational model for probability evolution and
refinement. The resulting model admits reductions to major
generative-model architectures and a continuum realization through
transport equations, gradient flows, and Fokker--Planck dynamics.
The different formulations are ultimately connected through common
principles of equilibrium, convergence, and stability.
}
\label{fig:levitas_unified_model}
\end{wrapfigure}
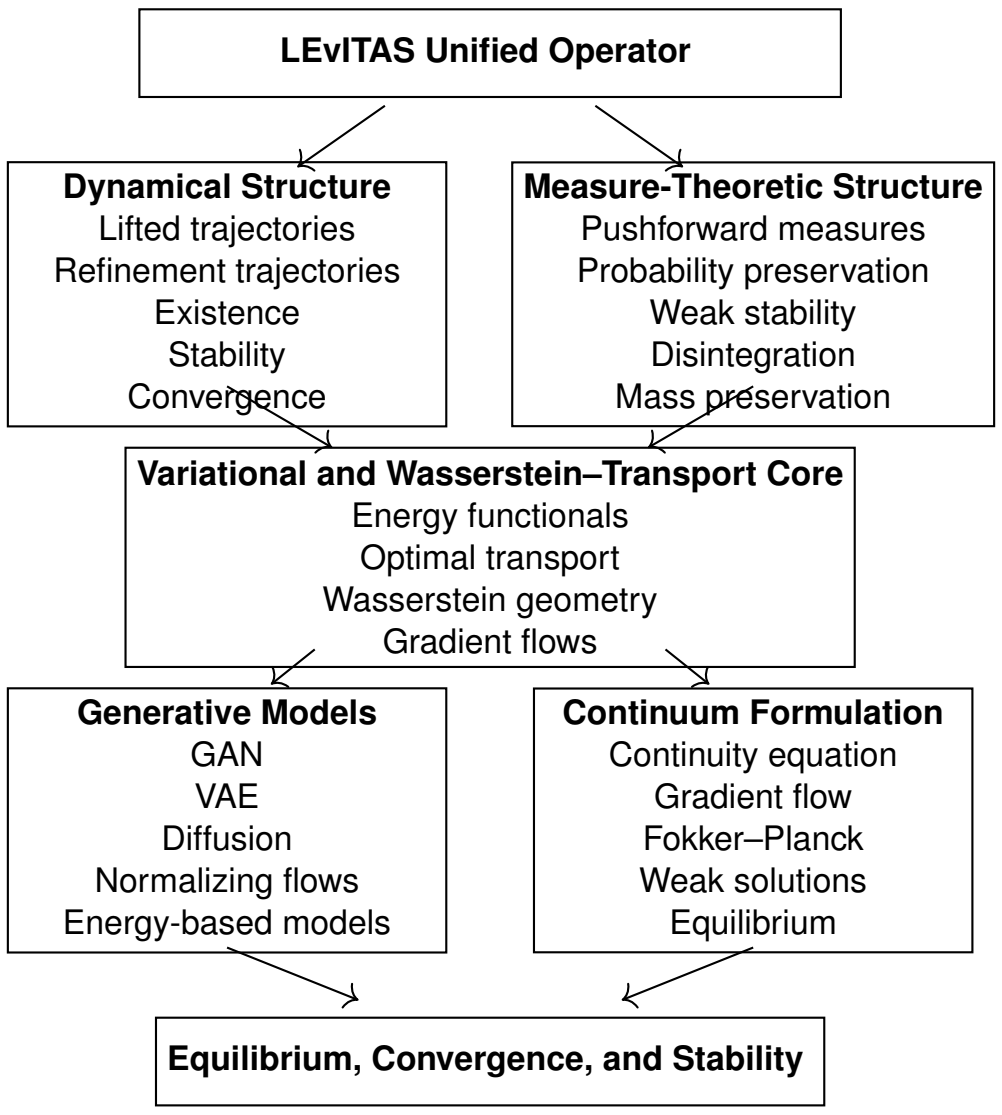

The relationships among these six analytical components are summarized
in Figure~\ref{fig:levitas_unified_model}, which illustrates the
structural hierarchy of LEvITAS and the connections among its
operator-theoretic, dynamical, measure-theoretic, variational,
generative, and continuum formulations.

At the top of this hierarchy, the \textbf{LEvITAS Unified Operator}
provides the generalized operator-theoretic foundation from which
two complementary analytical structures emerge:
the \textbf{Dynamical Structure} and the
\textbf{Measure-Theoretic Structure}.
The dynamical structure encompasses lifted trajectories,
refinement trajectories, existence, stability, and convergence,
describing the evolution of individual states under the LEvITAS
mechanism. In contrast, the measure-theoretic structure encompasses
pushforward probability measures, probability preservation, weak
stability, disintegration, and mass preservation, describing the
corresponding evolution of probability distributions.

These two structures are connected through the central
\textbf{Variational and Wasserstein--Transport Core}, which provides
the geometric and variational framework for coupling state-level
dynamics with probability-level evolution. Energy functionals
determine the variational objectives, optimal transport characterizes
the movement between probability measures, Wasserstein geometry
provides the associated metric structure, and gradient flows describe
the resulting continuous variational evolution.

From this central core, two complementary realizations follow.
The \textbf{Generative Models} branch comprises GANs, VAEs,
diffusion systems, normalizing flows, and energy-based models,
representing particular reductions of the unified LEvITAS
operator formulation. The \textbf{Continuum Formulation} branch
comprises the continuity equation, continuum gradient flows,
Fokker--Planck dynamics, weak solutions, and equilibrium
characterizations, representing the PDE-based realization of the
underlying dynamical, probabilistic, and variational structures.

Both branches are connected through the common analytical principles
of \textbf{Equilibrium, Convergence, and Stability}, highlighting
the relationships between the different formulations at their
respective mathematical levels. Thus, the figure provides a unified
structural interpretation of the preceding results, connecting
individual state evolution, probability-measure dynamics,
variational transport, generative-model reductions, and continuum
PDE systems within a common operator-theoretic framework.

Consequently, the continuum formulation brings together
\emph{variational refinement},
\emph{Wasserstein geometry},
\emph{stochastic probability evolution},
\emph{measure-valued dynamics},
\emph{continuum transport systems},
and
\emph{infinite-dimensional optimization}
within a unified continuum realization of LEvITAS.

The complete collection of results therefore establishes a coherent
analytical hierarchy in which lifted trajectories, refinement flows,
probability measures, optimal transport, generative architectures,
and continuum PDE systems are interconnected manifestations of the
underlying LEvITAS operator framework.

\section{Experimental Results}

Building upon the mathematical formulation of the LEvITAS operator model,
which establishes the connections between lifted trajectories, refinement
flows, probability measures, optimal transport, generative architectures, and
continuum PDE systems, this section investigates its practical applicability
to synthetic biomedical data generation across heterogeneous modalities. The
experimental evaluation examines whether the developed operator formulation
translates into consistent structural preservation, distributional alignment,
and cross-domain generalization when applied to real-world biomedical datasets.

Contemporary biomedical generative modelling faces significant challenges
arising from data heterogeneity and high dimensionality across multi-omics and
multimodal datasets~\cite{baiao2025multiomics, qoku2023multimodal}, limited 
sample availability and population
underrepresentation~\cite{vanbreugel2024synthetic}, complex statistical
dependencies~\cite{nathoo2019imaging, qoku2023multimodal}, and variations in
underlying biological and clinical characteristics that can introduce
distribution shifts across populations and healthcare
settings~\cite{zhang2022shifting, vanbreugel2024synthetic}. Biomedical
datasets encompass diverse modalities, including imaging, ultrasound, genomic
and proteomic profiles, morphological descriptors, and clinical and cognitive
measurements, which exhibit distinct dimensionality, statistical
distributions, and structural properties~\cite{baiao2025multiomics,
qoku2023multimodal, nathoo2019imaging, vanbreugel2024synthetic}. 
Consequently, existing generative approaches, particularly those relying on 
modality-specific architectures or processing assumptions, may struggle to 
preserve meaningful feature interactions and coordinated cross-modal 
relationships during synthetic data 
generation~\cite{stahlschmidt2022multimodal, ibrahim2025generative, 
marchesi2026coherent}.

To investigate the applicability of LEvITAS across heterogeneous biomedical settings, the model was evaluated using five independent biomedical applications: endometriosis laparoscopic imaging, antral follicle count (AFC) ultrasound phenotyping, amyotrophic lateral sclerosis (ALS) multi-omic profiling, antimicrobial resistance (AMR) genomic--morphological modelling, and Alzheimer's disease (AD) multimodal clinical--cognitive representation learning. These applications provide diverse experimental settings for examining the generality of the LEvITAS operator framework across distinct biological contexts, data representations, and modality configurations.

Particular emphasis is placed on cross-domain generalization and within-domain structural preservation, including statistical fidelity, preservation of latent feature hierarchies, and maintenance of relationships between heterogeneous modalities within individual disease datasets. Importantly, the five biomedical applications are evaluated independently rather than combined into a jointly modelled disease cohort. The common LEvITAS formulation provides the underlying generative model, while the associated data representations and structural relationships are examined within the context of each application. The basic matrix-factorization implementation provides the experimental realization of the LEvITAS generative framework, with its pivot-free factorization approach drawing methodological justification from the previously developed ci-fGBD~\cite{Thakur2025ciFGBD}, which establishes an fGBD-derived basis for pivot-free matrix factorization. Additionally, the Merkmal Treiber Plot (MTP) visualization approach~\cite{Thakur2025ciFGBD} is employed as a complementary analytical tool to investigate feature-level organization and structural consistency beyond conventional distributional similarity measures. Collectively, these experiments examine the potential of LEvITAS to provide a unified, structure-aware generative model applicable to heterogeneous biomedical data synthesis, multimodal representation learning, and robustness evaluation in biomedical artificial intelligence systems.

\subsection{Cross-Domain Generalization and Within-Domain Multimodal Structural Preservation}

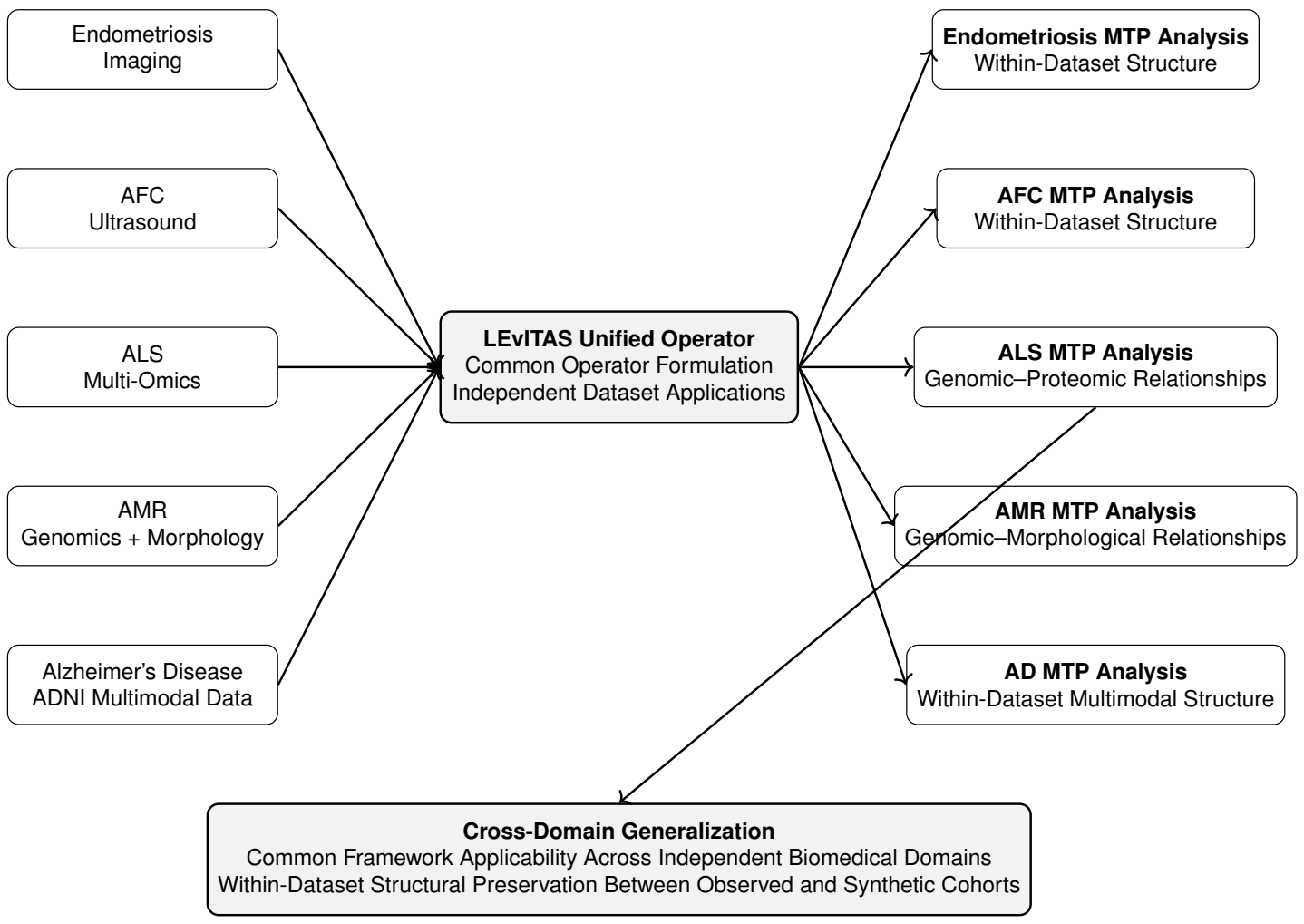
\begin{wrapfigure}{r}{0.52\textwidth}
\centering

\resizebox{0.95\linewidth}{!}{%
\begin{tikzpicture}[
    box/.style={
        draw, rounded corners, align=center,
        minimum width=3.4cm,
        minimum height=1.0cm,
        font=\scriptsize
    },
    framework/.style={
        draw, thick, rounded corners,
        fill=gray!10, align=center,
        minimum width=4.5cm,
        minimum height=1.4cm,
        font=\scriptsize
    },
    analysis/.style={
        draw, rounded corners, align=center,
        minimum width=4.0cm,
        minimum height=1.0cm,
        font=\scriptsize
    },
    arrow/.style={->, thick}
]

\node[box] (endo) at (0,2) {Endometriosis\\ Imaging};
\node[box] (afc) at (0,0) {AFC\\ Ultrasound};
\node[box] (als) at (0,-2) {ALS\\ Multi-Omics};
\node[box] (amr) at (0,-4) {AMR\\ Genomics + Morphology};
\node[box] (ad) at (0,-6) {Alzheimer's Disease\\ ADNI Multimodal Data};

\node[framework] (levitas) at (6,-2)
{\textbf{LEvITAS Unified Operator}\\
\scriptsize Common Operator Formulation\\
\scriptsize Independent Dataset Applications};

\node[analysis] (mtp1) at (12,2)
{\textbf{Endometriosis MTP Analysis}\\
\scriptsize Within-Dataset Structure};

\node[analysis] (mtp2) at (12,0)
{\textbf{AFC MTP Analysis}\\
\scriptsize Within-Dataset Structure};

\node[analysis] (mtp3) at (12,-2)
{\textbf{ALS MTP Analysis}\\
\scriptsize Genomic--Proteomic Relationships};

\node[analysis] (mtp4) at (12,-4)
{\textbf{AMR MTP Analysis}\\
\scriptsize Genomic--Morphological Relationships};

\node[analysis] (mtp5) at (12,-6)
{\textbf{AD MTP Analysis}\\
\scriptsize Within-Dataset Multimodal Structure};

\draw[arrow] (endo.east) -- (levitas.west);
\draw[arrow] (afc.east) -- (levitas.west);
\draw[arrow] (als.east) -- (levitas.west);
\draw[arrow] (amr.east) -- (levitas.west);
\draw[arrow] (ad.east) -- (levitas.west);

\draw[arrow] (levitas.east) -- (mtp1.west);
\draw[arrow] (levitas.east) -- (mtp2.west);
\draw[arrow] (levitas.east) -- (mtp3.west);
\draw[arrow] (levitas.east) -- (mtp4.west);
\draw[arrow] (levitas.east) -- (mtp5.west);

\node[framework, minimum width=10cm] (summary)
at (6,-8.2)
{\textbf{Cross-Domain Generalization}\\
\scriptsize Common Framework Applicability Across Independent Biomedical Domains\\
\scriptsize Within-Dataset Structural Preservation Between Observed and Synthetic Cohorts};

\draw[arrow] (mtp3.south) -- (summary.north);

\end{tikzpicture}%
}

\caption{
\textbf{Cross-domain generalization and within-domain structural preservation using LEvITAS and MTP.}
Five independent biomedical datasets are evaluated using the common LEvITAS operator framework. Each application retains its own dataset-specific modality configuration and associated feature relationships. MTP analysis is conducted within each respective application to examine preservation of latent feature hierarchy and relationships between heterogeneous modalities, where applicable, between observed and synthetic cohorts. Cross-domain generalization is assessed through the applicability of the common framework and the consistency of structural preservation across these independent biomedical applications, rather than through joint modelling of the five disease datasets.
}

\label{fig:mtp_cross_disease_safe}

\end{wrapfigure}

To evaluate the generality of LEvITAS across heterogeneous biomedical 
domains, we applied the unified operator framework independently to five 
distinct biomedical applications as shown in the Figure~\ref{fig:mtp_cross_disease_safe}, using the previously developed ci-fGBD 
model and Merkmal Treiber Plot (MTP) visualization 
approach~\cite{Thakur2025ciFGBD}. The evaluation included:
(i) endometriosis laparoscopic imaging, 
(ii) antral follicle count (AFC) ultrasound phenotyping, 
(iii) amyotrophic lateral sclerosis (ALS) multi-omic profiling, 
(iv) antimicrobial resistance (AMR) genomic--morphological modelling, and 
(v) Alzheimer's disease (AD) multimodal clinical--cognitive representation 
learning.

The purpose of this evaluation was to examine whether the same underlying 
LEvITAS operator formulation could be applied across distinct biomedical 
contexts while maintaining the structural organization of each individual 
dataset. Accordingly, cross-domain generalization refers to the applicability 
of the unified framework across independently evaluated disease domains, 
whereas multimodal structural preservation refers to the maintenance of 
relationships between heterogeneous modalities within each individual 
disease domain's multimodal datasets. Thus, the evaluation does not assume a jointly learned latent representation across the five diseases.

Despite differences in modality configuration, dimensionality, and biological context, LEvITAS maintained structural organization between observed and synthetic cohorts within the evaluated applications. The analysis focused on preservation of latent feature hierarchies, relationships between modalities where applicable, and distributional stability rather than direct biological mechanism inference.

In endometriosis imaging, LEvITAS preserved organization among imaging-derived features associated with lesion morphology, inflammatory characteristics, vascular patterns, and spatial heterogeneity. MTP analysis demonstrated consistent alignment between generated and observed cohorts, indicating preservation of imaging-derived feature relationships.

In AFC ultrasound cohorts, generated representations maintained relationships between follicular morphology, intensity distributions, and structural descriptors associated with ovarian reserve characterization. The synthetic cohorts retained the dominant organization observed in the real population while reducing unstable variance in sparsely represented regions.

In ALS multi-omic datasets, LEvITAS preserved hierarchical organization across molecular representations and maintained consistency between genomic and proteomic feature relationships. The resulting latent representations showed improved structural alignment compared with independent feature-wise generation.

In AMR datasets, LEvITAS preserved relationships between genomic signatures derived from k-mer representations and morphology-associated cellular phenotypes obtained from Cell Painting profiles. The generated representations maintained coordinated organization between genotype-associated and phenotype-associated features.

In AD multimodal cohorts derived from ADNI representations, LEvITAS preserved organization among cognitive, functional, imaging, and genetic features. MTP analysis demonstrated consistent alignment of multimodal feature hierarchies between synthetic and observed cohorts.

Across the five independently evaluated applications, LEvITAS preserved dominant latent organization while reducing stochastic variation associated with individual cohorts. These findings support the applicability of the unified operator framework across heterogeneous biomedical domains while demonstrating the preservation of dataset-specific structural organization within individual disease contexts.

Conventional generative evaluation metrics quantify distributional similarity 
but do not directly evaluate whether generated samples preserve relationships 
between heterogeneous features. To address this limitation, we applied the 
MTP as an interpretability analysis model for examining latent feature 
organization within individual observed and synthetic cohorts.

Across the evaluated datasets, MTP analysis demonstrated preservation of 
dominant feature hierarchies and relationships between heterogeneous 
features. Imaging datasets maintained organization among morphology-related 
features, while molecular datasets preserved relationships among genomic, 
proteomic, and resistance-associated representations. In multimodal 
applications, the analysis further examined the preservation of relationships 
between the corresponding modality-specific representations.

Importantly, MTP analysis indicated that LEvITAS reduces excessive stochastic 
dispersion while maintaining major latent organization patterns within each 
evaluated dataset. This behaviour suggests that the model regularizes 
synthetic representations without introducing artificial separation of 
biological groups.

The MTP analyses were performed within the respective biomedical 
applications. Their purpose was to assess structural preservation between 
observed and generated representations, rather than to establish a single 
jointly learned latent space shared by the five disease domains. Together, 
these observations demonstrate that LEvITAS preserves interpretable 
structural relationships beyond conventional distributional similarity 
measures.

\subsection{Summary of Multimodal Structural Preservation}

Collectively, the ci-fGBD and MTP analyses demonstrate that LEvITAS maintains 
clinically relevant latent organization within independently evaluated 
biomedical datasets. The framework preserves relationships among imaging-
derived morphology, molecular profiles, genomic signatures, and clinical 
phenotypes while generating synthetic cohorts with improved structural 
consistency.

The findings establish two complementary aspects of the framework's 
applicability: generalization of the underlying LEvITAS operator formulation 
across distinct biomedical domains, and preservation of dataset-specific 
multimodal or feature-level structural relationships within each domain. 
These findings support the use of LEvITAS for structure-aware synthetic 
cohort generation, multimodal representation learning, and robustness 
evaluation in biomedical artificial intelligence systems.

\subsection{Statistical Fidelity, Structural Relevance and Interpretability Analysis}

LEvITAS is designed not only to generate statistically plausible synthetic cohorts, but also to preserve clinically relevant structural organization across heterogeneous biomedical modalities. This capability is particularly important in biomedical settings characterized by limited cohort size, heterogeneous measurements, privacy constraints, and complex multimodal relationships, where synthetic data generation requires preservation of meaningful feature interactions beyond marginal distribution matching.

To evaluate whether LEvITAS preserves latent organization rather than reproducing superficial statistical similarity alone, we applied our previously developed interpretability model, the \textbf{MTP}~\cite{Thakur2025ciFGBD}, to both observed and synthetic cohorts.

MTP is not introduced as a new methodological contribution in this study; instead, it is used as an established visualization and interpretability model to examine latent feature hierarchy, directional organization, and relationships between heterogeneous features within generated representations. MTP enables visualization of feature-level contributions and hierarchical organization across biomedical data representations, including imaging-derived features, morphology descriptors, inflammatory measurements, genomic and proteomic profiles, cognitive assessments, neuroimaging biomarkers, and antimicrobial resistance-associated signatures. By comparing MTP representations between observed and LEvITAS-generated cohorts, we assess whether synthetic augmentation maintains the structural relationships present within the original datasets.

Across all evaluated disease domains, LEvITAS preserved dominant feature organization while reducing excessive cohort-specific variation. MTP analysis indicated that synthetic representations maintained coordinated relationships between heterogeneous features rather than independently perturbing individual variables. This property is important for biomedical modelling because clinically relevant phenotypes frequently arise from interactions among multiple biological and observational measurements.

To examine the generality of the LEvITAS model, we evaluated its application across five independently modelled biomedical datasets:
(i) endometriosis laparoscopic imaging, 
(ii) antral follicle count (AFC) ultrasound-based reproductive phenotyping, 
(iii) amyotrophic lateral sclerosis (ALS) multi-omic genomic--proteomic profiling, 
(iv) Alzheimer’s disease (AD) multimodal cognitive, imaging, CSF, and genetic representation learning, and 
(v) antimicrobial resistance (AMR) datasets integrating genomic signatures with morphological response profiles.

Representative MTP analyses are shown in Figures~\ref{fig:mtp_endometriosis_cifgbd}--\ref{fig:cluster_amr} (Appendix~B). Within each application, LEvITAS maintained alignment between synthetic and observed latent feature organization despite substantial differences in modality, dimensionality, and biological context.

Importantly, the comparison across these applications concerns the generality of the underlying operator framework and the preservation of structural organization within each respective dataset. It does not imply that the five disease datasets were combined into a single training cohort or that their modality-specific representations were jointly optimized in a common disease-independent latent space.

Collectively, these analyses demonstrate that LEvITAS preserves disease-associated latent structure and relationships among heterogeneous biomedical features across diverse independently evaluated datasets. These findings support the potential utility of structure-preserving synthetic generation for biomedical data augmentation, representation learning, robustness evaluation, and future translational modelling applications.

\section{Conclusion}

This work introduced LEvITAS as a generalized variational lifting model
on probability spaces and developed a unified mathematical model for
structure-preserving probabilistic evolution. The formulation connects
measurable lifting--projection operators, dynamically consistent
refinement systems, path-space transport evolution, variational
probability dynamics, Wasserstein transport geometry, continuum PDE
flows, and measure-valued gradient evolution within a common operator-theoretic model. 
The mathematical development establishes several complementary layers
of the LEvITAS theory. At the dynamical level, the model provides
existence of continuous lifted trajectories, uniform stability under
latent perturbations, weak convergence of lifted path measures, and
existence of dynamically consistent variational minimizers. At the
refinement level, the theory establishes local existence and uniqueness,
monotone energy dissipation, equilibrium characterization, exponential
stability under strong convexity, preservation of topological
invariants, and monotone convergence of discrete refinement iterations.
These results provide the analytical foundation for interpreting LEvITAS
as a stable variational evolution rather than solely as a generative
mapping. 
At the measure-theoretic level, LEvITAS preserves probability under
measurable pushforward mappings and provides weak stability, Wasserstein
metric structure, Prokhorov-type compactness, induced synthetic
probability measures, disintegration structures, and mass preservation
for measure-valued refinement flows. The optimal-transport formulation
further establishes Monge-type and Kantorovich minimizers, Wasserstein
geodesic structure, continuity of displacement interpolation, the
Benamou--Brenier dynamical formulation, and energy dissipation along
Wasserstein gradient flows. These results identify probability geometry
and optimal transport as the variational mechanism through which
structural evolution is formulated. 
The resulting theory establishes a mathematical connection between
operator theory, probability geometry, optimal transport, variational dynamics, generative modeling, and continuum PDE evolution. 
A principal consequence of this unified formulation is that a broad
family of generative systems can be interpreted as constrained
reductions of the same generalized variational lifting model. In this
sense, GAN systems, variational autoencoders, diffusion systems,
normalizing flows, and energy-based models are not treated merely as
independent generative architectures, but as particular realizations
within a common operator-theoretic and variational hierarchy. 
The continuum formulation extends this structure from finite-dimensional
or discrete representations to probability evolution on spaces of
measures. In particular, the LEvITAS construction provides the
conceptual progression 
$
P_Z
\longrightarrow
(\pi\circ\mathcal L_\theta)_\#P_Z
\longrightarrow
\mu_t
\longrightarrow
\rho_t,
$ 
thereby connecting latent probability structure, measurable
lifting--projection mappings, transport evolution, measure-valued
dynamics, and continuum refinement within a single variational system.
The continuity equation, continuum gradient flow, Fokker--Planck
dynamics, weak-solution theory, and equilibrium characterization provide
the corresponding continuum realization of this construction. 
Importantly, the mathematical formulation is complemented by an
empirical evaluation across heterogeneous biomedical data modalities.
The experimental analysis evaluated LEvITAS using endometriosis
laparoscopic imaging, antral follicle count ultrasound phenotyping,
amyotrophic lateral sclerosis multi-omic profiling, antimicrobial
resistance genomic--morphological modelling, and Alzheimer's disease
multimodal clinical--cognitive representation learning. Despite the
substantial differences in modality, dimensionality, and biological
context, the analyses consistently examined preservation of latent
feature hierarchy, cross-modal relationships, and distributional
stability between observed and synthetic cohorts. 
The cross-domain evaluation provides empirical evidence that the
structure-preserving principle underlying LEvITAS is not restricted to
a single modality. In endometriosis imaging, the generated
representations preserved organization associated with lesion
morphology, inflammatory characteristics, vascular patterns, and
spatial heterogeneity. In AFC ultrasound data, the synthetic
representations retained relationships among follicular morphology,
intensity distributions, and structural descriptors. In ALS multi-omic
data, LEvITAS maintained hierarchical organization across molecular
representations and genomic--proteomic relationships. In AMR data, the
model preserved coordinated organization between genomic signatures and
morphology-associated cellular phenotypes. In Alzheimer's disease
representations, LEvITAS maintained organization across demographi,
clinical, and genetic features. 
The Merkmal Treiber Plot analysis further provides an interpretable
empirical view of this structural preservation. Rather than evaluating
synthetic generation only through marginal distributional similarity,
the analysis examines latent feature hierarchy, directional organization,
and cross-modal relationships between observed and generated cohorts.
Across the evaluated domains, the results indicate preservation of
dominant latent organization together with a reduction in excessive
stochastic dispersion. This supports the interpretation of LEvITAS as a
structure-aware generative model whose objective extends beyond
independent feature-wise distribution matching. 
The agreement between the mathematical formulation and the empirical
evaluation is therefore an important aspect of the LEvITAS model.
The theory establishes a variational and transport-based mechanism for
structure-preserving probabilistic evolution, while the multimodal
experiments demonstrate that the resulting synthetic representations
can retain dominant latent organization across substantially different
biomedical domains. Collectively, these results support the interpretation
of LEvITAS as a generalized structure-preserving generative model
rather than as a modality-specific synthetic data generator. 
At the same time, the empirical results define important boundaries of
the present formulation. The generated representations remain
constrained by the information content and resolution of the original
measurements. For imaging-based applications, LEvITAS can preserve
imaging-visible patterns as well as statistical structures that may
reflect underlying biological processes that are not directly
observable in the acquired modality. Likewise, preservation of
genomic, proteomic, phenotypic, or multimodal organization does not by
itself imply mechanistic modelling of the underlying biological
pathways. These measurements constitute incomplete observations of
complex biological systems. Rare pathological states and highly
heterogeneous subgroups may also remain underrepresented when cohort
sizes are limited. Accordingly, LEvITAS should be interpreted as a
structure-preserving generative model rather than as a complete
biological simulation model. 
The combined theoretical and empirical findings therefore position
LEvITAS at the intersection of mathematical probability theory,
variational analysis, optimal transport, generative modelling, and
continuum dynamics. The model provides a generalize language in which
latent representations, probability measures, transport paths,
variational refinements, generative architectures, and continuum
evolution can be understood as different levels of the same underlying
operator structure. 
Several directions remain open. These include regularity theory for
nonlinear refinement flows, convergence analysis of variational
refinement semigroups, high-dimensional Wasserstein transport geometry,
stochastic stability of lifted trajectories, numerical approximation of
continuum LEvITAS systems, and more systematic characterization of the
relationship between mathematical structural preservation and
domain-specific biological relevance. Further work is also required to
evaluate rare and heterogeneous subpopulations, strengthen multimodal
validation, and establish quantitative links between latent structural
preservation and downstream biomedical modeling performance.

In conclusion, LEvITAS provides a unified mathematical and empirical
model for generalized probabilistic structural evolution. Its
operator-theoretic formulation supplies the analytical foundation,
variational and Wasserstein geometry provide the mechanism for
structure-preserving evolution, generative-model reductions establish
its architectural generality, and continuum PDE formulations extend
the model to infinite-dimensional probability dynamics. The
cross-domain biomedical evaluation further demonstrates that this
structure-preserving formulation can retain meaningful latent
organization across heterogeneous imaging, ultrasound, multi-omic,
genomic--morphological, and multimodal clinical representations. 
Thus, LEvITAS extends beyond conventional generative modeling toward a
general mathematical theory of structure-preserving probabilistic
evolution on spaces of measures, with a unified connection between
measurable operators, variational refinement, transport geometry,
generative representations, and continuum probability dynamics.

\newpage

\section{Appendix A: Mathematical Formulations, Theorems, Proofs}
\subsection{Overview}
\label{sec:overview_supplementary_results}
The LEvITAS model development begins with the dynamical theory of lifted
trajectory systems, where we establish existence of continuous lifted paths
(Theorem~\ref{thm:existence_continuous_lifted_paths}), uniform stability
under latent perturbations
(Theorem~\ref{thm:uniform_stability_lifted_paths}), weak convergence of
lifted path measures (Theorem~\ref{thm:weak_stability_lifted_path_measures}),
and existence of dynamically consistent variational minimizers
(Theorem~\ref{thm:existence_dynamically_consistent_liftings}). 
Building on this foundation, we develop the theory of iterative
topological refinement and variational evolution. In particular,
we prove local existence and uniqueness of refinement trajectories
(Theorem~\ref{thm:local_existence_refinement_flows}),
monotone energy dissipation along refinement flows
(Theorem~\ref{thm:energy_dissipation_law}),
characterization of stationary equilibrium states
(Theorem~\ref{thm:characterization_equilibrium_states}),
exponential stability under strong convexity
(Theorem~\ref{thm:exponential_stability}),
preservation of topological invariants during refinement
(Theorem~\ref{thm:topological_preservation_continuous_refinement}),
and monotone convergence of discrete refinement iterations
(Theorem~\ref{thm:monotone_descent_discrete_refinement}). 
We next establish the probabilistic and measure-theoretic
foundations of LEvITAS through probability preservation under
measurable pushforward mappings
(Theorem~\ref{thm:probability_preservation}),
weak stability of pushforward probability measures
(Theorem~\ref{thm:stability_pushforward_measures}),
metric structure of Wasserstein probability geometry
(Theorem~\ref{thm:wasserstein_metric_structure}),
compactness characterization via Prokhorov tightness
(Theorem~\ref{thm:prokhorov_compactness}),
existence of induced synthetic probability measures for stochastic
liftings
(Theorem~\ref{thm:existence_induced_synthetic_measures}),
existence of disintegration structures on Polish spaces
(Theorem~\ref{thm:existence_disintegration}),
and mass preservation for measure-valued LEvITAS refinement flows
(Theorem~\ref{thm:mass_preservation}). 
These probabilistic foundations naturally lead to the optimal
transport and Wasserstein-geometric formulation of the model.
Within this setting, we establish existence of Monge-type LEvITAS
transport minimizers
(Theorem~\ref{thm:existence_monge_levitas_transport}),
existence of Kantorovich transport minimizers
(Theorem~\ref{thm:existence_kantorovich_minimizers}),
geodesic structure of Wasserstein probability spaces
(Theorem~\ref{thm:geodesic_structure_wasserstein_space}),
continuity of displacement interpolation flows
(Theorem~\ref{thm:continuity_displacement_flows}),
dynamical transport characterizations via the
Benamou--Brenier formulation
(Theorem~\ref{thm:benamou_brenier_formulation}),
and energy dissipation laws for Wasserstein gradient flows
(Theorem~\ref{thm:energy_dissipation_wasserstein_flows}). 
We further derive the continuum PDE formulation of LEvITAS through
the continuity equation
\eqref{eq:levitas_continuity},
the continuum gradient-flow system
\eqref{eq:continuum_levitas},
and the Fokker--Planck evolution equation~\cite{Jordan1998VariationalFokkerPlanck}
\eqref{eq:fokker_planck_levitas}. Within this continuum setting,
we establish entropy dissipation laws
(Theorem~\ref{thm:entropy_dissipation}),
existence theory for weak continuum solutions
(Theorem~\ref{thm:existence_weak_solutions}),
and equilibrium characterizations for stationary continuum states
(Theorem~\ref{thm:equilibrium_characterization}). These results
provide the infinite-dimensional continuum counterpart of the
variational, transport, and stochastic formulations developed
throughout the supplementary model. 
Finally, we establish a unified operator-theoretic reduction
model for modern generative architectures. In particular,
we prove a general reduction principle
(Theorem~\ref{thm:general_reduction_theorem}),
together with rigorous reductions of LEvITAS to GAN systems
(Theorem~\ref{thm:gan_reduction}),
variational autoencoders
(Theorem~\ref{thm:vae_reduction}),
diffusion systems
(Theorem~\ref{thm:diffusion_reduction}),
normalizing flows
(Theorem~\ref{thm:normalizing_flow_reduction}),
and energy-based models
(Theorem~\ref{thm:energy_based_reduction}).
These reductions culminate in a hierarchical containment theorem
for major generative architectures
(Theorem~\ref{thm:hierarchical_containment}). 
Collectively, these results establish LEvITAS as a unified
operator-theoretic model connecting variational refinement,
probability evolution, Wasserstein geometry, stochastic dynamics,
continuum PDE systems, optimal transport theory~\cite{Villani2003TopicsOT, Santambrogio2015OptimalTransport, Villani2009OTOldNew}, and modern
generative modeling within a common analytical structure.

\subsection{Generalized Lifting Operators}
\label{sec:generalized_lifting}

The LEvITAS model models generative systems through
operator-induced transformations between latent probability
spaces, structured synthetic spaces, and observable realization
spaces. The principal mathematical distinction of the model
is the explicit factorization
\begin{equation}
\mathcal Z
\xrightarrow{\mathcal L_\theta}
\mathcal X_{syn}
\xrightarrow{\pi}
\mathcal X_{real},
\label{eq:levitas_factorization}
\end{equation}
which separates latent probabilistic structure from synthetic
structural generation and observable realization. 
This section develops the analytical foundations of generalized
lifting systems, including existence, stability, variational
optimization, and approximation properties.

\subsection{Preliminaries}
\label{subsec:preliminaries}

Let $(\mathcal Z,\Sigma_Z,P_Z)$ be a complete probability space.
Let $\mathcal X_{syn}$ and $\mathcal X_{real}$ be separable metric spaces equipped with Borel sigma-algebras $\mathcal B(\mathcal X_{syn}),
\mathcal B(\mathcal X_{real})$. Denote by
$\mathcal P(\mathcal X)$ the space of Borel probability measures on a measurable space
$\mathcal X$. Throughout this section, weak convergence of probability
measures is denoted by $\Rightarrow$.

\subsection{Generalized Lifting Systems}
\label{subsec:generalized_lifting_systems}

\begin{definition}[Generalized Lifting Operator]
\label{def:generalized_lifting_operator}


A unified generalized stochastic lifting operator is a parameterized family of Markov probability kernels that maps latent representations into probability distributions over synthetic representations.

Let \((\mathcal Z_r,\mathcal B_{\mathcal Z_r})\) and \((\mathcal X_r,\mathcal B_{\mathcal X_r})\) denote measurable latent and synthetic spaces, respectively, where \(r\) identifies the representation type, including individual observations, probability distributions, or longitudinal trajectories.

For each representation type \(r\), the stochastic lifting operator is defined as 
$
K_\theta^{(r)}:
\mathcal Z_r\times\mathcal B(\mathcal X_r)
\rightarrow [0,1],
$ 
parameterized by \(\theta\in\Theta\), satisfying: 
(1) For every \(z\in\mathcal Z_r\), \(K_\theta^{(r)}(\cdot\mid z)\) is a probability measure on \(\mathcal X_r\).
(2) For every \(A\in\mathcal B(\mathcal X_r)\), \(K_\theta^{(r)}(A\mid\cdot)\) is measurable. 
Given a latent probability measure \(P_Z^{(r)}\), the induced synthetic probability measure is defined by 
$
P_{\mathrm{syn}}^{(r)}(A)
=
\int_{\mathcal Z_r}
K_\theta^{(r)}(A\mid z)
P_Z^{(r)}(dz),
$ 
for all \(A\in\mathcal B(\mathcal X_r)\).

In particular, the operator supports the following representations:
(i) Point-level lifting: For an individual latent observation \(z\), the synthetic observation is generated according to
$
X_{\mathrm{syn}}\mid Z=z
\sim K_\theta^{(\mathrm{point})}(\cdot\mid z).
$
(ii) Distribution-level lifting: For a latent population distribution \(P_Z\), the induced synthetic population distribution is
$
P_{\mathrm{syn}}
=
\mathfrak L_\theta(P_Z)
=
P_Z K_\theta^{(\mathrm{point})}.
$
(iii) Trajectory-level lifting: For a longitudinal latent trajectory \(z_{1:T}\), the synthetic trajectory is generated according to the joint conditional distribution
$
X_{1:T}^{\mathrm{syn}}\mid Z_{1:T}=z_{1:T}
\sim
K_\theta^{(\mathrm{trajectory})}
(\cdot\mid z_{1:T}),
$ 
where the joint kernel specifies the dependence structure among synthetic observations across visits.

The stochastic lifting mechanism may incorporate persistent latent patient characteristics, visit-specific measurement uncertainty, and temporally correlated stochastic effects.

The deterministic generalized lifting operator is recovered as a special case when the kernel is concentrated at a single output:
$
K_\theta^{(r)}(dx\mid z)
=
\delta_{\mathcal L_\theta^{(r)}(z)}(dx),
$ 
where \(\delta_a\) denotes the Dirac probability measure concentrated at \(a\).
\end{definition}
Thus, the unified generalized stochastic lifting operator accommodates deterministic and stochastic transformations, individual observations, population distributions, and longitudinal trajectories within a common probabilistic framework.

\begin{definition}[Projection Operator]
\label{def:projection_operator}

A projection operator is a measurable map
\begin{equation}
\pi:
\mathcal X_{syn}
\rightarrow
\mathcal X_{real}.
\label{eq:projection_operator}
\end{equation}

\end{definition}

The observable probability measure induced by the LEvITAS
system is
\begin{equation}
P_{obs}
=
(\pi\circ\mathcal L_\theta)_\#P_Z.
\label{eq:observable_pushforward}
\end{equation}

\subsection{Existence of Lifted Measures}
\label{subsec:existence_lifted_measures}

\begin{theorem}[Existence of Synthetic Measures]
\label{thm:existence_synthetic_measure}
Let $
\mathcal L_\theta:
\mathcal Z
\rightarrow
\mathcal X_{syn}
$ be measurable. Then there exists a unique probability measure
$P_{syn}\in\mathcal P(\mathcal X_{syn})$ satisfying $P_{syn}=(\mathcal L_\theta)_\#P_Z$.
\end{theorem}

\begin{proof}

Define $P_{syn}(A)=P_Z(\mathcal L_\theta^{-1}(A))$ for each $A\in\mathcal B(\mathcal X_{syn}).$ Since $\mathcal L_\theta$ is measurable, $\mathcal L_\theta^{-1}(A)\in\Sigma_Z.$ Hence $P_{syn}$ is well-defined. Countable additivity follows directly from the countable
additivity of $P_Z.$ Moreover, $P_{syn}(\mathcal X_{syn})=P_Z(\mathcal Z)=1.$ Thus
$P_{syn}$ is a probability measure. Uniqueness follows from uniqueness of pushforward measures.
\end{proof}

\subsection{Operator Factorization}
\label{subsec:operator_factorization}

The central structural principle of LEvITAS is the
decomposition in \eqref{eq:levitas_factorization},
which separates:
1. latent probabilistic structure, 2. synthetic structural generation, 3. observable realization.


\begin{theorem}[Measurable Factorization]
\label{thm:measurable_factorization}

Suppose
$\mathcal L_\theta:
\mathcal Z
\rightarrow
\mathcal X_{syn}
$
and
$
\pi:
\mathcal X_{syn}
\rightarrow
\mathcal X_{real}
$ are measurable. Then
$\pi\circ\mathcal L_\theta
$ is measurable.
\end{theorem}

\begin{proof}
Let $A\in\mathcal B(\mathcal X_{real}).$ Since $\pi$ is measurable, $\pi^{-1}(A)\in
\mathcal B(\mathcal X_{syn}).$ Since
$\mathcal L_\theta$ is measurable, $\mathcal L_\theta^{-1}(\pi^{-1}(A))\in \Sigma_Z.$
Therefore, $(\pi\circ\mathcal L_\theta)^{-1}(A)
=
\mathcal L_\theta^{-1}(\pi^{-1}(A))$ is measurable. Hence
$\pi\circ\mathcal L_\theta$ is measurable.
\end{proof}

\subsection{Variational Formulation}
\label{subsec:variational_formulation}

Let $P_{data}\in\mathcal P(\mathcal X_{real})$ denote the target observable distribution. Define the baseline LEvITAS functional
$\mathcal J:\Theta\rightarrow\mathbb R\cup\{+\infty\}$ by
\vspace{-0.5cm}
\begin{align}
\mathcal J(\theta)
&=
D\!\left(
(\pi\circ\mathcal L_\theta)_\# P_Z,
P_{data}
\right)
\nonumber\\
&\quad
+
\lambda_1\mathcal D_{feat}(\theta)
+
\lambda_2\mathcal D_{struct}(\theta)
+
\lambda_3\mathcal D_{dyn}(\theta),
\label{eq:levitas_functional_base}
\end{align}
where:
$D$ is a statistical discrepancy,
$\mathcal D_{feat}$ is a feature discrepancy,
$\mathcal D_{struct}$ is a structural consistency penalty,
and $\mathcal D_{dyn}$ is a dynamical consistency penalty.

The baseline formulation corresponds to the unimodal or modality-independent setting and may be regarded as a special case of the multimodal LEvITAS model.

For multimodal learning, the functional is extended with an additional cross-modal consistency regularizer:
\begin{align}
\mathcal J_{full}(\theta)
&=
D\!\left(
(\pi\circ\mathcal L_\theta)_\# P_Z,
P_{data}
\right)
\nonumber\\
&\quad
+
\lambda_1\mathcal D_{feat}(\theta)
+
\lambda_2\mathcal D_{struct}(\theta)
+
\lambda_3\mathcal D_{dyn}(\theta)
\nonumber\\
&\quad
+
\lambda_4\mathcal D_{crossmodal}(\theta),
\label{eq:levitas_functional_full}
\end{align}
where
$\mathcal D_{crossmodal}$ enforces consistency, alignment, and shared representation structure across heterogeneous modalities.

The baseline functional in \eqref{eq:levitas_functional_base}
is recovered as the special case 
$\lambda_4 = 0.$
\subsection{Existence of Optimal Liftings}
\label{subsec:existence_optimal_liftings}

\begin{assumption}[Reflexive Parameter Space]
\label{assump:reflexive_parameter_space}
The parameter space $\Theta$ is a reflexive Banach space.
\end{assumption}


\begin{assumption}[Coercivity]
\label{assump:coercivity}

The functional
$\mathcal J$ is coercive:
$\|\theta_n\|_\Theta
\rightarrow
\infty
\quad
\Longrightarrow
\quad
\mathcal J(\theta_n)
\rightarrow
\infty.
$
\end{assumption}


\begin{assumption}[Weak Lower Semicontinuity]
\label{assump:wls}

The functional
$\mathcal J$ is weakly lower semicontinuous on $\Theta.$

\end{assumption}


\begin{lemma}[Boundedness of Minimizing Sequences]
\label{lem:bounded_minimizing_sequences}

Let
$
\{\theta_n\}\subseteq\Theta
$
satisfy
$
\mathcal J(\theta_n)
\rightarrow
\inf_{\theta\in\Theta}\mathcal J(\theta).
$ 
Then
$
\{\theta_n\}
$
is bounded in
$
\Theta.
$

\end{lemma}

\begin{proof}

Assume otherwise. Then
$
\|\theta_n\|_\Theta
\rightarrow
\infty.
$ 
By Assumption~\ref{assump:coercivity},
$
\mathcal J(\theta_n)
\rightarrow
\infty,
$
contradicting the minimizing property.

\end{proof}


\begin{lemma}[Weak Compactness]
\label{lem:weak_compactness}

Let
$
\Theta
$
be reflexive and let
$
\{\theta_n\}
$
be bounded in
$
\Theta.
$ 
Then there exists a subsequence
$
\{\theta_{n_k}\}
$
and
$
\theta^\ast\in\Theta
$
such that
$
\theta_{n_k}
\rightharpoonup
\theta^\ast
$
weakly in
$
\Theta.
$

\end{lemma}

\begin{proof}

Since
$
\Theta
$
is reflexive, bounded subsets are weakly relatively compact.

\end{proof}


\begin{theorem}[Existence of Optimal Lifting Operators]
\label{thm:existence_optimal_lifting}

Suppose:
Assumptions~\ref{assump:reflexive_parameter_space},\ref{assump:coercivity}, and \ref{assump:wls} holds, Then there exists
$
\theta^\ast\in\Theta
$
such that
$
\mathcal J(\theta^\ast)
=
\inf_{\theta\in\Theta}
\mathcal J(\theta).
$

\end{theorem}

\begin{proof}

Define
$
m
=
\inf_{\theta\in\Theta}
\mathcal J(\theta).
$ 
Let
$
\{\theta_n\}
$ 
be a minimizing sequence satisfying
$
\mathcal J(\theta_n)\rightarrow m.
$ 
By Lemma~\ref{lem:bounded_minimizing_sequences},
$
\{\theta_n\}
$
is bounded. 
By Lemma~\ref{lem:weak_compactness},
there exists a subsequence satisfying
$
\theta_n\rightharpoonup\theta^\ast.
$ 
Weak lower semicontinuity implies
$
\mathcal J(\theta^\ast)
\leq
\liminf_{n\to\infty}\mathcal J(\theta_n)
=
m.
$ 
Since
$
m
\leq
\mathcal J(\theta^\ast),
$
we conclude
$
\mathcal J(\theta^\ast)=m.
$

\end{proof}


To illustrate how the assumptions of Theorem~\ref{thm:existence_optimal_lifting} arise in practice, we consider a finite-dimensional parameterization of the lifting operator. Let $\Theta=\mathbb R^d$ and define the empirical loss functional 
$
\mathcal J(\theta)
=
\frac{1}{N}\sum_{i=1}^N
\ell\big(f_\theta(x_i),y_i\big)
+
\frac{\lambda}{2}\|\theta\|_2^2,
$ 
where $\{(x_i,y_i)\}_{i=1}^N$ denotes a finite dataset, $f_\theta$ represents the parameterized lifting operator, $\ell$ is a continuous and nonnegative loss function, and $\lambda>0$ is a regularization parameter. 
We verify that the assumptions of Theorem~\ref{thm:existence_optimal_lifting} hold for this objective. 
First, $\Theta=\mathbb R^d$ is a finite-dimensional Banach space and is therefore reflexive. Second, since the empirical loss is nonnegative, we have
$
\mathcal J(\theta)
\geq
\frac{\lambda}{2}\|\theta\|_2^2.
$ 
It follows that
$
\|\theta_n\|_2\rightarrow\infty
\quad\Longrightarrow\quad
\mathcal J(\theta_n)\rightarrow\infty,
$ 
which establishes the coercivity of $\mathcal J$. 
Finally, suppose that the parameterized lifting operator $\theta\mapsto f_\theta(x_i)$ is continuous for every training sample $x_i$. Since $\ell$ is continuous and the regularization term is continuous, $\mathcal J$ is continuous on $\mathbb R^d$. Moreover, weak convergence and strong convergence coincide in finite-dimensional spaces. Therefore, $\mathcal J$ is weakly lower semicontinuous. 
Hence, all assumptions of Theorem~\ref{thm:existence_optimal_lifting} are satisfied, and there exists at least one optimal parameter $\theta^\ast\in\Theta$ such that
$
\mathcal J(\theta^\ast)
=
\min_{\theta\in\Theta}\mathcal J(\theta).
$ 
This example demonstrates that the abstract existence result applies directly to regularized finite-dimensional parameterizations of lifting operators.

\subsection{Stability of Lifted Measures}
\label{subsec:stability_lifted_measures}
\begin{assumption}[Stability of Lifting Operators and Readout Continuity]
\label{assump:stability_conditions}

Suppose that $\mathcal X_{syn}$ is a Polish space and that the lifting operators satisfy 
$
\mathcal L_{\theta_n}(z)
\rightarrow
\mathcal L_\theta(z)
$ 
for $P_Z$-almost every $z\in\mathcal Z$. 
Furthermore, let 
$
\pi:\mathcal X_{syn}\rightarrow\mathcal Y
$ 
be a measurable readout map, where $\mathcal Y$ is a Polish space. Assume that $\pi$ is continuous at $\mathcal L_\theta(z)$ for $P_Z$-almost every $z\in\mathcal Z$.

\end{assumption}

\subsection{Approximation Properties}
\label{subsec:approximation_properties}

\begin{definition}[Universal Structured Approximation]
\label{def:universal_structured_approximation}

A lifting family
$
\{\mathcal L_\theta\}_{\theta\in\Theta}
$
is universally structure-approximating if for every
$
P^\ast\in\mathcal P(\mathcal X_{real})
$
and every
$
\varepsilon>0,
$
there exists
$
\theta\in\Theta
$
such that
$
D\left(
(\pi\circ\mathcal L_\theta)_\#P_Z,
P^\ast
\right)
<
\varepsilon.
$

\end{definition}

\subsection{Universal Structured Approximation}
\label{subsec:universal_structured_approximation}

\begin{definition}[Universal Structured Approximation]
\label{def:universal_structured_approximation}

A lifting family 
$
\{\mathcal L_\theta\}_{\theta\in\Theta}
$ 
is universally structure-approximating if for every
$
P^\ast\in\mathcal P(\mathcal X_{real})
$
and every $\varepsilon>0$, there exists
$\theta\in\Theta$ such that
$
D\left(
(\pi\circ\mathcal L_\theta)_\#P_Z,
P^\ast
\right)
<
\varepsilon.
$

\end{definition}

\subsection{Constructive Verification of Universal Approximation}
\label{subsec:constructive_universal_approximation}

The universal approximation property in
Definition~\ref{def:universal_structured_approximation}
is an abstract capacity condition. To demonstrate that this property can hold for a concrete lifting architecture, we introduce a simple finite-dimensional toy model and establish a constructive approximation result.

\subsubsection{Scalar toy model.}

Let 
$
\mathcal Z=[0,1],\qquad
Z\sim\operatorname{Unif}[0,1],
$ 
and consider 
$
\mathcal X_{syn}
=
\mathcal X_{real}
=
[0,1].
$ 
We choose the observable readout to be the identity:
$
\pi(x)=x.
$ 
Let $\mathfrak L_{pc}$ denote the family of all measurable, finite piecewise-constant lifting operators 
$
\mathcal L:[0,1]\rightarrow[0,1].
$ 
Each lifting operator is defined by a finite partition of the latent interval and assigns a constant observable value to each partition element.

For example, consider the target probability measure 
$
P^\ast
=
0.3\delta_{0.2}
+
0.7\delta_{0.8}.
$
Define 
$
\mathcal L(u)
=
\begin{cases}
0.2,&0\leq u<0.3,\\
0.8,&0.3\leq u\leq1.
\end{cases}
$ 
Since $Z$ is uniformly distributed on $[0,1]$, we obtain 
$
(\pi\circ\mathcal L)_\#P_Z
=
0.3\delta_{0.2}
+
0.7\delta_{0.8}
=
P^\ast.
$ 
Thus, this target distribution is represented exactly by a piecewise-constant lifting operator.

We next establish that the same family can approximate arbitrary probability distributions on $[0,1]$.


\begin{theorem}[Universal Approximation by Piecewise-Constant Liftings]
\label{thm:piecewise_constant_universal_approximation}

Let 
$
\mathcal Z=\mathcal X_{syn}=\mathcal X_{real}=[0,1],
$ 
let $P_Z$ be the uniform probability measure on $[0,1]$, and let $\pi$ be the identity readout.

For every target probability measure 
$
P^\ast\in\mathcal P([0,1])
$ 
and every integer $n\geq1$, there exists a piecewise-constant lifting operator 
$
\mathcal L_n\in\mathfrak L_{pc}
$ 
such that
$
W_1\left(
(\pi\circ\mathcal L_n)_\#P_Z,
P^\ast
\right)
\leq\frac{1}{2n}.
$ 
Consequently,
$
\inf_{\mathcal L\in\mathfrak L_{pc}}
W_1\left(
(\pi\circ\mathcal L)_\#P_Z,
P^\ast
\right)
=0
$ 
for every $P^\ast\in\mathcal P([0,1])$. 
Hence, the piecewise-constant lifting family is universally structure-approximating with respect to the Wasserstein-1 distance.

\end{theorem}

\begin{proof}

Fix an arbitrary target distribution
$P^\ast\in\mathcal P([0,1])$.

Partition the observable space into $n$ intervals: 
$
I_j=
\left[
\frac{j-1}{n},\frac{j}{n}
\right),
\qquad j=1,\ldots,n-1,
$ 
with the final interval including the endpoint $1$.

Define 
$
p_j=P^\ast(I_j),
\qquad
c_j=\frac{2j-1}{2n},
$ 
where $c_j$ is the midpoint of the $j$-th interval. 
Next, define the cumulative probabilities
$
a_0=0,\qquad
a_j=\sum_{k=1}^{j}p_k,
\qquad j=1,\ldots,n.
$ 
Since $\sum_{j=1}^np_j=1$, we have $a_n=1$.

Construct the lifting operator
$
\mathcal L_n(u)=c_j
\quad\text{if}\quad
a_{j-1}\leq u<a_j,
$
with the endpoint $u=1$ assigned to the last interval. 
This is a measurable, finite piecewise-constant lifting operator. 
Because $Z$ is uniformly distributed on $[0,1]$, the probability of selecting the value $c_j$ is precisely $p_j$. Therefore, 
$
(\mathcal L_n)_\#P_Z
=
\sum_{j=1}^{n}p_j\delta_{c_j}.
$ 
Since $\pi$ is the identity,
$
(\pi\circ\mathcal L_n)_\#P_Z
=
\sum_{j=1}^{n}p_j\delta_{c_j}.
$ 
To estimate the approximation error, consider a random variable $X\sim P^\ast$ and define its quantized counterpart
$
T_n(X)=c_j
\quad\text{whenever }X\in I_j.
$ 
By construction, $T_n(X)$ has distribution
$(\pi\circ\mathcal L_n)_\#P_Z$. 
Furthermore, every point in $I_j$ lies within distance $1/(2n)$ of its midpoint. Thus,
$
|X-T_n(X)|
\leq\frac{1}{2n}
\qquad\text{almost surely}.
$ 
The joint distribution of $(X,T_n(X))$ is a coupling of the target and approximating measures. By the definition of the Wasserstein-1 distance,

$
\begin{aligned}
W_1\left(
P^\ast,
(\pi\circ\mathcal L_n)_\#P_Z
\right)
&\leq
\mathbb E\left[|X-T_n(X)|\right]\\
&\leq\frac{1}{2n}.
\end{aligned}
$ 
Taking $n\rightarrow\infty$ establishes 
$
W_1\left(
(\pi\circ\mathcal L_n)_\#P_Z,
P^\ast
\right)\rightarrow0.
$ 
Since the target distribution was arbitrary, the result follows.
\end{proof}


\subsubsection{Extension to multimodal observable distributions.}

The scalar construction illustrates the approximation mechanism. In LEvITAS, however, the observable space may consist of several heterogeneous modalities. 
For example, consider the joint observable space 
$
\mathcal X_{real}
=
\mathcal X_{img}\times\mathcal X_{cog},
$ 
where $\mathcal X_{img}$ represents imaging observations and $\mathcal X_{cog}$ represents cognitive measurements. 
We equip these spaces with bounded compatible metrics $d_{img}$ and $d_{cog}$, respectively, and define the product metric 
$
d\big((x_1,y_1),(x_2,y_2)\big)
=
\alpha d_{img}(x_1,x_2)
+
\beta d_{cog}(y_1,y_2),
$ 
where $\alpha,\beta>0$ and $\alpha+\beta=1$.

\begin{theorem}[Universal Approximation of Joint Multimodal Distributions]
\label{thm:multimodal_universal_approximation}

Let $\mathcal X_{img}$ and $\mathcal X_{cog}$ be Polish spaces equipped with bounded compatible metrics. Let 
$
\mathcal X_{real}
=
\mathcal X_{img}\times\mathcal X_{cog}
$
be equipped with the corresponding bounded product metric. 
Suppose that 
$
\mathcal Z=[0,1],
\qquad
Z\sim\operatorname{Unif}[0,1],
$ 
and that the synthetic space is 
$
\mathcal X_{syn}=\mathcal X_{real},
$ 
with identity readout 
$
\pi=\operatorname{id}_{\mathcal X_{real}}.
$ 
Consider the family of all finite piecewise-constant measurable lifting operators 
$
\mathcal L:[0,1]\rightarrow\mathcal X_{real}.
$ 
Then, for every joint probability distribution 
$
P^\ast\in
\mathcal P(\mathcal X_{img}\times\mathcal X_{cog})
$ 
and every $\varepsilon>0$, there exists a lifting operator $\mathcal L_\varepsilon$ in this family such that 
$
W_1\left(
(\pi\circ\mathcal L_\varepsilon)_\#P_Z,
P^\ast
\right)
<\varepsilon.
$ 
Therefore, the piecewise-constant lifting family can approximate arbitrary joint imaging-cognition distributions in the Wasserstein-1 metric.

\end{theorem}

\begin{proof}

Let $P^\ast$ be an arbitrary probability measure on the joint observable space. 
Since the observable space is Polish, $P^\ast$ is tight. Hence, for any $\delta>0$, there exists a compact set $K$ such that
$
P^\ast(K)\geq1-\delta.
$ 
By compactness, for any $\eta>0$, the set $K$ admits a finite measurable partition
$
K=A_1\cup\cdots\cup A_m
$ 
such that each partition element has diameter at most $\eta$. 
Let
$
A_0=\mathcal X_{real}\setminus K.
$ 
Choose a representative point $x_j\in A_j$ for each nonempty $A_j$, and choose an arbitrary $x_0\in\mathcal X_{real}$. 
Define
$
p_j=P^\ast(A_j),
\qquad j=0,\ldots,m.
$ 
Partition the latent interval $[0,1]$ into consecutive intervals of lengths $p_0,\ldots,p_m$. Define a piecewise-constant lifting operator that outputs $x_j$ on the interval corresponding to $p_j$. 
The induced observable distribution is
$
Q_\varepsilon
=
\sum_{j=0}^{m}p_j\delta_{x_j}.
$ 
Couple a random variable $X\sim P^\ast$ with its representative point $T(X)=x_j$ whenever $X\in A_j$. 
For $X\in K$, we have
$
d(X,T(X))\leq\eta.
$
On the complement of $K$, boundedness of the metric gives
$
d(X,T(X))\leq M,
$
where $M$ is the diameter bound of the observable space.
Consequently,

$
\begin{aligned}
W_1(P^\ast,Q_\varepsilon)
&\leq
\mathbb E[d(X,T(X))]\\
&\leq
\eta P^\ast(K)+MP^\ast(K^c)\\
&\leq\eta+M\delta.
\end{aligned}
$

Choosing $\eta$ and $\delta$ sufficiently small yields
$
W_1(P^\ast,Q_\varepsilon)<\varepsilon.
$ 
Since $Q_\varepsilon$ is generated by a finite piecewise-constant lifting operator and the readout is the identity, the desired result follows.
\end{proof}


\begin{corollary}[Approximation Capacity]
\label{cor:approximation_capacity}

If the lifting family is universally structure-approximating with respect to the discrepancy $D$, then, for every target observable probability measure
$P^\ast\in\mathcal P(\mathcal X_{real})$
and every $\varepsilon>0$, there exists
$\theta\in\Theta$ such that
$
D\left(
(\pi\circ\mathcal L_\theta)_\#P_Z,
P^\ast
\right)<\varepsilon.
$

\end{corollary}

\begin{proof}
The result follows immediately from
Definition~\ref{def:universal_structured_approximation},
which establishes the stated approximation property for every
target observable probability measure and every positive
tolerance $\varepsilon$.
\end{proof}

\subsubsection{Implications for LEvITAS.}

The preceding results establish that universal approximation is achievable for an explicit class of lifting operators and readouts. In particular, the multimodal result concerns joint probability distributions rather than independent marginal distributions. 
Thus, the approximation claim applies not only to imaging and cognitive distributions separately, but also to their joint dependence structure, including statistical associations between modalities.

For a general LEvITAS architecture, however, the approximation property depends on the expressive capacity of its lifting family and readout maps. The toy model establishes a constructive sufficient example.

In particular, an implemented multimodal model must be sufficiently expressive to represent or approximate the joint observable configurations required by the target distribution through its synthetic lifting and readout factorization.

\subsection{Structure Operators}
\label{sec:structure_operators}

The LEvITAS model extends conventional generative systems by introducing
explicit structure operators acting on synthetic and observable spaces.
Unlike purely distributional formulations that optimize only sample realism,
the structure operator formalism permits direct enforcement of geometric,
topological, semantic, relational, and dynamical consistency.

This section develops the analytical foundations of
structure-preserving generative systems within the LEvITAS model.

\subsection{Structure Spaces}
\label{subsec:structure_spaces}

Let $\mathcal X$ denote a measurable sample space and let $\mathcal M$ be a 
metric structure space equipped with metric $d_{\mathcal M}$. The structure space may encode:
geometric invariants, topological structure, semantic consistency, relational
connectivity, physical constraints, temporal dependencies, dynamical
invariants.
Throughout this section, assume $\mathcal M$ is complete and separable.

\subsection{Definition of Structure Operators}
\label{subsec:definition_structure_operators}

\begin{definition}[Structure Operator]
\label{def:structure_operator}
A structure operator is a measurable map
$
\mathcal S:
\mathcal X
\rightarrow
\mathcal M.
$
\end{definition}

The operator
$
\mathcal S
$
extracts structural information from observable or synthetic samples.

Examples include:
graph Laplacians, homology descriptors, geometric embeddings, semantic
feature maps, physical state operators, temporal dependency representations.

\subsection{Structural Consistency}
\label{subsec:structural_consistency}

Let
$
\pi:
\mathcal X_{syn}
\rightarrow
\mathcal X_{real}
$
be the LEvITAS projection operator.

The central structural requirement of LEvITAS is consistency between
projected synthetic structure and observable structure.


\begin{definition}[Structural Consistency]
\label{def:structural_consistency}
A lifting-projection system is structurally consistent if
$
\mathcal S(\pi(x_{syn}))
\approx
\mathcal S(x_{real}).
$
\end{definition}
To quantify consistency, define the structural discrepancy functional
$
\mathcal D_{struct}:
\mathcal M
\times
\mathcal M
\rightarrow
\mathbb R_+.
$
A canonical choice is
\begin{equation}
\mathcal D_{struct}(x,y)
=
d_{\mathcal M}(x,y)^p,
\quad p\geq1.
\label{eq:structural_discrepancy}
\end{equation}

\subsection{Structure-Preserving LEvITAS Functionals}
\label{subsec:structure_preserving_functionals}
Let
$\mathcal S:\mathcal X_{real}\to\mathcal X_{struct}$
denote a measurable structural representation map, where
$\mathcal X_{struct}$ is the corresponding structural feature
space. To avoid requiring an arbitrary pointwise correspondence
between synthetic and real observations, we define the structural
component of the LEvITAS functional through a discrepancy between
the distributions of their structural representations:
\begin{equation}
\mathcal J_{struct}(\theta)
=
D_{struct}
\left(
(\mathcal S\circ\pi\circ\mathcal L_\theta)_\#P_Z,
\mathcal S_\#P_{data}
\right).
\label{eq:structural_functional}
\end{equation}
Here, $(\mathcal S\circ\pi\circ\mathcal L_\theta)_\#P_Z$
denotes the distribution of structural representations generated
by the lifting and readout operators, whereas
$\mathcal S_\#P_{data}$ denotes the corresponding distribution
obtained from real observations. Thus, the structural term
compares the structural properties of the generated and observed
data at the distributional level, without requiring a
sample-to-sample pairing.

The full variational formulation becomes
\begin{align}
\mathcal J(\theta)
&=
D
\left(
(\pi\circ\mathcal L_\theta)_\#P_Z,
P_{data}
\right)
\nonumber\\
&\quad
+
\lambda_1\mathcal D_{feat}(\theta)
+
\lambda_2\mathcal J_{struct}(\theta)
\nonumber\\
&\quad
+
\lambda_3\mathcal D_{dyn}(\theta).
\label{eq:structural_levitas}
\end{align}

\subsection{Continuity of Structure Operators}
\label{subsec:continuity_structure_operators}

\begin{assumption}[Continuity of Structure Operators]
\label{assump:continuity_structure_operator}

The structure operator
$
\mathcal S:
\mathcal X
\rightarrow
\mathcal M
$
is continuous.
\end{assumption}


\begin{lemma}[Preservation of Convergent Sequences]
\label{lem:preservation_convergent_sequences}

Suppose
$
x_n \rightarrow x
$
in
$
\mathcal X.
$ 
If
$
\mathcal S
$
is continuous, then
$
\mathcal S(x_n)
\rightarrow
\mathcal S(x)
$
in
$
\mathcal M.
$
\end{lemma}

\begin{proof}
The result follows directly from continuity of
$
\mathcal S.
$
\end{proof}

\subsection{Structural Consistency Convergence}
\label{subsec:structural_consistency_convergence}

\begin{theorem}[Structural Consistency Theorem]
\label{thm:structural_consistency}

Let $(\mathcal M,d_{\mathcal M})$ be a metric space, and let
$\mathcal S$ be a continuous structural representation map.
Suppose that the structural discrepancy
$\mathcal D_{\mathrm{struct}}:\mathcal M\times\mathcal M
\to[0,\infty)$ satisfies 
$
\mathcal D_{\mathrm{struct}}(u,v)=0
\iff u=v,
$ 
and that convergence with respect to
$\mathcal D_{\mathrm{struct}}$ implies convergence in
$d_{\mathcal M}$; that is, 
$
\mathcal D_{\mathrm{struct}}(u_n,u)\to0
\quad\Longrightarrow\quad
d_{\mathcal M}(u_n,u)\to0.
$ 
If
$
\mathcal D_{\mathrm{struct}}
\left(
\mathcal S(\pi(x_n)),
\mathcal S(x)
\right)
\to0,
$ 
then 
$
\mathcal S(\pi(x_n))
\to
\mathcal S(x)
\qquad\text{in }\mathcal M.
$

\end{theorem}

\begin{proof}
By assumption, 
$
\mathcal D_{\mathrm{struct}}
\left(
\mathcal S(\pi(x_n)),
\mathcal S(x)
\right)
\to0.
$ 
Since convergence with respect to
$\mathcal D_{\mathrm{struct}}$ implies convergence in
the metric $d_{\mathcal M}$, it follows that
$
d_{\mathcal M}
\left(
\mathcal S(\pi(x_n)),
\mathcal S(x)
\right)
\to0.
$
Therefore,
$
\mathcal S(\pi(x_n))
\to
\mathcal S(x)
\qquad\text{in }\mathcal M.
$
\end{proof}

\subsection{Path Lifting and Dynamical Systems}
\label{sec:path_lifting_dynamical_systems}

The LEvITAS model naturally extends from static generative systems to
trajectory-valued and dynamical generative systems. In this setting, latent
trajectories are lifted into structured synthetic trajectories prior to
observable realization. 
This section develops the analytical foundations of path lifting, dynamical
consistency, trajectory stability, and variational dynamical optimization.

\subsection{Trajectory Spaces}
\label{subsec:trajectory_spaces}

Let
$
T=[0,\tau]
$
be a compact temporal interval. 
Define the latent trajectory space
$
\Gamma_Z
=
C(T;\mathcal Z)
$
and the synthetic trajectory space
$
\Gamma_{syn}
=
C(T;\mathcal X_{syn}),
$
where
$
C(T;\mathcal X)
$
denotes the space of continuous paths equipped with the uniform norm
$
\|\gamma\|_\infty
=
\sup_{t\in T}
 d(\gamma(t),0).
$ 
Similarly define the observable trajectory space
$
\Gamma_{real}
=
C(T;\mathcal X_{real}).
$ 
Throughout this section, assume: 
$
\mathcal Z,
\mathcal X_{syn},
\mathcal X_{real}
$
are complete separable metric spaces, all trajectories are continuous unless otherwise stated.

\subsection{Path Lifting Operators}
\label{subsec:path_lifting_operators}

The central dynamical object in LEvITAS is the lifted trajectory.


\begin{definition}[Path Lifting Operator]
\label{def:path_lifting_operator}

Let
$
\mathcal L_\theta:
\mathcal Z
\rightarrow
\mathcal X_{syn}
$
be a generalized lifting operator. 
The induced path lifting operator is defined by
$
\widetilde{\mathcal L}_\theta:
\Gamma_Z
\rightarrow
\Gamma_{syn}
$
through
$
(\widetilde{\mathcal L}_\theta\gamma)(t)
=
\mathcal L_\theta(\gamma(t)).
$
\end{definition}

For a latent trajectory
$
\gamma_Z(t),
$
the lifted synthetic trajectory is therefore
$
\widetilde\gamma(t)
=
\mathcal L_\theta(\gamma_Z(t)).
$ 
The corresponding observable trajectory becomes
$
\gamma_{obs}(t)
=
\pi(\widetilde\gamma(t)).
$

\subsection{Well-Posedness of Lifted Trajectories}
\label{subsec:wellposedness_lifted_trajectories}

We first establish existence and continuity of lifted trajectories.


\begin{assumption}[Continuity of the Lifting Operator]
\label{ass:continuity_lifting_operator}

Suppose
$
\mathcal L_\theta:
\mathcal Z
\rightarrow
\mathcal X_{syn}
$
is continuous.
\end{assumption}


\begin{theorem}[Existence of Continuous Lifted Paths]
\label{thm:existence_continuous_lifted_paths}

Suppose:
$
\gamma_Z\in C(T;\mathcal Z),
$ 
$
\mathcal L_\theta
$
is continuous. 
Then
$
\widetilde\gamma
=
\mathcal L_\theta\circ\gamma_Z
\in
C(T;\mathcal X_{syn}).
$
\end{theorem}

\begin{proof}
Since
$
\gamma_Z
$
is continuous and
$
\mathcal L_\theta
$
is continuous, the composition
$
\widetilde\gamma
=
\mathcal L_\theta\circ\gamma_Z
$
is continuous. 
Hence
$
\widetilde\gamma
\in
C(T;\mathcal X_{syn}).
$
\end{proof}

\subsection{Dynamical Consistency}
\label{subsec:dynamical_consistency}

\subsection{Patient-Specific and Population-Level
Dynamical Consistency}
\label{subsec:refined_dynamical_consistency}

Dynamical consistency requires that generated trajectories
reproduce the longitudinal behavior of the underlying
real-world population. However, the appropriate notion
of consistency depends on whether synthetic trajectories
can be paired with individual real patients.

We distinguish two settings:
\begin{definition}[i. Patient-Specific Dynamical Consistency]
\label{def:patient_specific_dynamic_consistency}

Suppose a model receives the baseline information of a
real patient and generates a predicted longitudinal
trajectory. 
Let $t_0<t_1<\cdots<t_m$ denote the observation times, and let 
$y_{ij}\in\mathcal X_{\mathrm{real}}$ denote the observed measurement for patient $i$ at time $t_j$.

Let 
$
\widehat y_{ij}(\theta)
=
\pi\left(
\mathcal L_\theta
\left(
\gamma_Z^{(i)}(t_j)
\right)
\right)
$ 
denote the corresponding model-generated prediction,
where the latent trajectory is conditioned on the
patient's baseline information. 
The patient-specific dynamical discrepancy is defined by 
$
\mathcal D_{\mathrm{dyn}}^{\mathrm{pair}}(\theta)
=
\frac{1}{N}
\sum_{i=1}^{N}
\frac{
\displaystyle\sum_{j=1}^{m}
w_{ij}
d_{\mathcal X_{\mathrm{real}}}
\left(
\widehat y_{ij}(\theta),y_{ij}
\right)^2
}{
\displaystyle\sum_{j=1}^{m}w_{ij}
},
$ 
where $w_{ij}\geq 0$ are observation weights.

The model is patient-specifically dynamically consistent
when this discrepancy is sufficiently small.
\end{definition}

The squared trajectory discrepancy penalizes deviations
across the entire observation period rather than assessing
individual visits independently.

In particular, for scalar cognitive scores, consider
the observed trajectory $(20,18,16)$ and the generated trajectory $(20,19,18).$ 
Their mean squared trajectory error is $\mathcal D_{\mathrm{dyn}}^{\mathrm{pair}}
=
\frac{1}{3}
\left[
(20-20)^2+
(19-18)^2+
(18-16)^2
\right]
=
\frac{5}{3}.
$

Thus, the discrepancy explicitly penalizes the
underestimation of longitudinal decline.


\begin{definition}[ii. Population-Level Dynamical Consistency]
\label{def:population_dynamic_consistency}

Suppose synthetic patients are generated without a
one-to-one correspondence to real patients. 
Let 
$P_{\mathrm{real}}^\Gamma
\in\mathcal P(\Gamma_{\mathrm{real}})
$ 
denote the probability measure on real observable
trajectories, and let
$
P_{\theta}^\Gamma
=
(\pi\circ\mathcal L_\theta)_\#P_\Gamma
$ 
denote the corresponding synthetic observable path
measure.

For observation times $T_m=\{t_1,\ldots,t_m\},$ define the trajectory evaluation map 
$E_{T_m}(\gamma)
=
\big(
\gamma(t_1),\ldots,\gamma(t_m)
\big).
$

The induced finite-dimensional trajectory
distributions are
$
Q_{\mathrm{real}}
=
(E_{T_m})_\#P_{\mathrm{real}}^\Gamma
$
and
$
Q_\theta
=
(E_{T_m})_\#P_\theta^\Gamma.
$

The population-level dynamical discrepancy is defined by
$
\mathcal D_{\mathrm{dyn}}^{\mathrm{dist}}(\theta)
=
W_p^p
\left(
Q_\theta,Q_{\mathrm{real}}
\right),
$
where $W_p$ denotes the $p$-Wasserstein distance. 
The model is population-level dynamically consistent
when
$
\mathcal D_{\mathrm{dyn}}^{\mathrm{dist}}(\theta)
$
is sufficiently small.
\end{definition}

This formulation compares the joint distributions
of longitudinal measurements. Consequently, it can
capture differences in the magnitude and variability
of decline across visits, as well as dependencies
between measurements at different time points.

Importantly, the comparison is not restricted to
matching the marginal distributions at individual
visits.

For example, similarity between the six-month and
twelve-month score distributions separately does not
necessarily imply similarity in the joint distribution
of six-month and twelve-month trajectories.

\subsection{Stability of Path Liftings}
\label{subsec:stability_path_liftings}

We next establish stability of lifted trajectories under perturbations of latent trajectories.


\begin{definition}[Lipschitz Lifting Operator]
\label{def:lipschitz_lifting_operator}

A lifting operator
$
\mathcal L_\theta
$
is Lipschitz if there exists
$
C>0
$

such that
$
d_{\mathcal X_{syn}}
(\mathcal L_\theta(z_1),\mathcal L_\theta(z_2))
\leq
C d_{\mathcal Z}(z_1,z_2)
$
for all
$
z_1,z_2\in\mathcal Z.
$
\end{definition}


\begin{theorem}[Uniform Stability of Lifted Paths]
\label{thm:uniform_stability_lifted_paths}

Suppose: 
$
\mathcal L_\theta
$
is Lipschitz with constant
$
C,
$ 
$
\gamma_n
\rightarrow
\gamma
$
uniformly in
$
C(T;\mathcal Z).
$ 
Then
$
\widetilde{\mathcal L}_\theta(\gamma_n)
\rightarrow
\widetilde{\mathcal L}_\theta(\gamma)
$
uniformly in
$
C(T;\mathcal X_{syn}).
$
\end{theorem}

\begin{proof}
For each
$
t\in T,
$
Lipschitz continuity yields
\begin{align*}
&d_{\mathcal X_{syn}}
\left(
\mathcal L_\theta(\gamma_n(t)),
\mathcal L_\theta(\gamma(t))
\right)
\\
&\leq
C d_{\mathcal Z}(\gamma_n(t),\gamma(t)).
\end{align*} 
Taking the supremum over
$
t\in T
$
gives
$
\|\widetilde{\mathcal L}_\theta(\gamma_n)
-
\widetilde{\mathcal L}_\theta(\gamma)\|_\infty
\leq
C
\|\gamma_n-\gamma\|_\infty.
$

Since
$
\gamma_n
\rightarrow
\gamma
$
uniformly, the result follows.
\end{proof}

\subsection{Weak Convergence of Path Measures}
\label{subsec:weak_convergence_path_measures}

The LEvITAS model naturally induces probability measures on path spaces. 
Let
$
P_\Gamma
\in
\mathcal P(\Gamma_Z)
$
be a latent trajectory measure. 
The induced lifted path measure is
$
(\widetilde{\mathcal L}_\theta)_\#P_\Gamma.
$


\begin{theorem}[Weak Stability of Lifted Path Measures]
\label{thm:weak_stability_lifted_path_measures}

Suppose:
$
\Gamma_{syn}
$
is Polish, 
$
\widetilde{\mathcal L}_{\theta_n}
\rightarrow
\widetilde{\mathcal L}_\theta
$
uniformly on compact subsets,
 the induced path measures are tight. 
Then
$
(\widetilde{\mathcal L}_{\theta_n})_\#P_\Gamma
\Rightarrow
(\widetilde{\mathcal L}_\theta)_\#P_\Gamma.
$
\end{theorem}

\begin{proof}
Let
$
f:\Gamma_{syn}\rightarrow\mathbb R
$
be bounded and continuous. 
Uniform convergence on compact subsets implies
$
f(\widetilde{\mathcal L}_{\theta_n}(\gamma))
\rightarrow
f(\widetilde{\mathcal L}_\theta(\gamma))
$
pointwise. 
Since
$
f
$
is bounded, dominated convergence yields
\begin{align*}
&\int_{\Gamma_Z}
 f(\widetilde{\mathcal L}_{\theta_n}(\gamma))
\,dP_\Gamma(\gamma)
\\
&\rightarrow
\int_{\Gamma_Z}
 f(\widetilde{\mathcal L}_\theta(\gamma))
\,dP_\Gamma(\gamma).
\end{align*}

Hence
$
(\widetilde{\mathcal L}_{\theta_n})_\#P_\Gamma
\Rightarrow
(\widetilde{\mathcal L}_\theta)_\#P_\Gamma.
$
\end{proof}

\subsection{Variational Dynamical Optimization}
\label{subsec:variational_dynamical_optimization}

The dynamical component of LEvITAS optimization is
designed to ensure that generated trajectories reproduce
the longitudinal behavior of real-world observations. 
Depending on the generation setting, dynamical consistency
can be formulated either at the individual patient level
or at the population level. 
For patient-specific prediction, the dynamical objective
penalizes deviations between predicted and observed
measurements at the corresponding follow-up visits. 
For population-level generation, the objective compares
the joint distributions of real and synthetic longitudinal
trajectories.

Accordingly, define the unified dynamical discrepancy by
$
\mathcal D_{\mathrm{dyn}}(\theta)
=
\begin{cases}
\mathcal D_{\mathrm{dyn}}^{\mathrm{pair}}(\theta),
& \text{patient-specific generation},\\[6pt]
\mathcal D_{\mathrm{dyn}}^{\mathrm{dist}}(\theta),
& \text{population-level generation}.
\end{cases}
$

where the two discrepancy functionals are defined in
Section~\ref{subsec:refined_dynamical_consistency}.

The dynamical optimization objective is then
$
\mathcal J_{\mathrm{dyn}}(\theta)
=
\mathcal D_{\mathrm{dyn}}(\theta).
$ 
When an application requires simultaneous control of
individual prediction accuracy and population-level
trajectory fidelity, a combined dynamical discrepancy
can alternatively be introduced as
$
\mathcal D_{\mathrm{dyn}}^{\mathrm{combined}}(\theta)
=
\alpha
\mathcal D_{\mathrm{dyn}}^{\mathrm{pair}}(\theta)
+
(1-\alpha)
\mathcal D_{\mathrm{dyn}}^{\mathrm{dist}}(\theta),
$
where
$
0\leq\alpha\leq1
$
controls the relative contribution of patient-specific
and population-level consistency. 
The combined formulation is applicable when both
paired longitudinal observations and population-level
trajectory distributions are available.


The complete LEvITAS optimization functional is
therefore defined by
\begin{align}
\mathcal J(\theta)
&=
D\left(
(\pi\circ\mathcal L_\theta)_\#P_Z,
P_{\mathrm{data}}
\right)
\nonumber\\
&\quad+
\lambda_1\mathcal D_{\mathrm{feat}}(\theta)
+
\lambda_2\mathcal D_{\mathrm{struct}}(\theta)
\nonumber\\
&\quad+
\lambda_3\mathcal D_{\mathrm{dyn}}(\theta).
\label{eq:dynamical_levitas}
\end{align}

Here,
\[
\lambda_1,\lambda_2,\lambda_3\geq0
\]
are regularization weights controlling the relative
contributions of feature-level, structural, and
dynamical consistency.

The dynamical component ensures that optimization
accounts not only for the plausibility of generated
measurements but also for their longitudinal evolution.

In particular, the patient-specific formulation
penalizes errors in predicted follow-up measurements,
whereas the population-level formulation encourages
the synthetic population to reproduce the joint
statistical behavior of real longitudinal trajectories.

Thus, the LEvITAS functional incorporates both
cross-sectional fidelity and longitudinal dynamical
consistency within a unified variational optimization
framework.

\subsection{Existence of Dynamically Consistent Liftings}
\label{subsec:existence_dynamically_consistent_liftings}

We now establish the existence of minimizers for the
refined dynamical LEvITAS functional. 
The existence result applies to both patient-specific
and population-level dynamical consistency, provided
that the corresponding dynamical discrepancy is
incorporated into the complete optimization functional
and the resulting functional satisfies the appropriate
variational assumptions.

In particular, the existence of an optimal lifting
requires control of minimizing sequences and lower
semicontinuity of the complete LEvITAS objective with
respect to the underlying optimization parameters.


\begin{assumption}[Existence of Dynamic Minimizers]
\label{ass:existence_dynamic_minimizers}

Suppose that the following conditions hold:

\begin{enumerate}

\item The parameter space $\Theta$ is a nonempty,
weakly closed subset of a reflexive Banach space
$\mathcal B$.

\item The LEvITAS functional
$
\mathcal J:\Theta\rightarrow
(-\infty,+\infty]
$
is proper and bounded below.

\item The functional $\mathcal J$ is coercive, in the
sense that
$
\|\theta\|_{\mathcal B}\rightarrow\infty
\quad\Longrightarrow\quad
\mathcal J(\theta)\rightarrow+\infty.
$

\item The functional $\mathcal J$ is sequentially
weakly lower semicontinuous, meaning that
$
\theta_n\rightharpoonup\theta
\quad\Longrightarrow\quad
\mathcal J(\theta)
\leq
\liminf_{n\rightarrow\infty}
\mathcal J(\theta_n).
$

\end{enumerate}

Here, properness means that $\mathcal J$ is finite
at least at one admissible parameter value.

\end{assumption}


\begin{theorem}[Existence of Dynamically Consistent Liftings]
\label{thm:existence_dynamically_consistent_liftings}

Under Assumption
\ref{ass:existence_dynamic_minimizers}, there exists
$\theta^\ast\in\Theta$ such that
$
\boxed{
\mathcal J(\theta^\ast)
=
\inf_{\theta\in\Theta}\mathcal J(\theta).
}
$ 
The result holds for either the patient-specific or
population-level dynamical formulation, as well as
their combined formulation, provided that the
corresponding complete LEvITAS functional satisfies
the stated assumptions.

\end{theorem}

\begin{proof}

The proof follows from the direct method of the
calculus of variations. 
Since $\mathcal J$ is proper and bounded below,
its infimum is finite. Define
$
J^\ast
=
\inf_{\theta\in\Theta}\mathcal J(\theta).
$ 
Let
$
\{\theta_n\}_{n=1}^{\infty}
\subset\Theta
$
be a minimizing sequence satisfying
$
\mathcal J(\theta_n)
\longrightarrow J^\ast.
$ 
By coercivity, the minimizing sequence
$\{\theta_n\}$ is bounded in $\mathcal B$. 
Since $\mathcal B$ is reflexive, every bounded
sequence admits a weakly convergent subsequence.
Consequently, there exist a subsequence, which we
continue to denote by $\{\theta_n\}$, and some
$\theta^\ast\in\mathcal B$ such that
$
\theta_n\rightharpoonup\theta^\ast.
$ 
Because $\Theta$ is weakly closed and
$\theta_n\in\Theta$, the weak limit satisfies
$
\theta^\ast\in\Theta.
$ 
By sequential weak lower semicontinuity of
$\mathcal J$,
$
\mathcal J(\theta^\ast)
\leq
\liminf_{n\rightarrow\infty}
\mathcal J(\theta_n).
$ 
Using the minimizing property of the sequence,
we obtain
$
\mathcal J(\theta^\ast)
\leq
J^\ast.
$ 
On the other hand, since $\theta^\ast\in\Theta$,
the definition of the infimum gives
$
\mathcal J(\theta^\ast)
\geq
\inf_{\theta\in\Theta}\mathcal J(\theta)
=
J^\ast.
$

Therefore,
$
\boxed{
\mathcal J(\theta^\ast)
=
\inf_{\theta\in\Theta}\mathcal J(\theta).
}
$ 
Hence, there exists a minimizer of the refined
dynamical LEvITAS functional.

\end{proof}

\subsection{Main Dynamical Results}
\label{subsec:main_dynamical_results}

The principal analytical results of this section are summarized as follows:
Existence of continuous lifted trajectories
(Theorem~\ref{thm:existence_continuous_lifted_paths}), 
Uniform stability under latent perturbations
(Theorem~\ref{thm:uniform_stability_lifted_paths}), 
Weak convergence of lifted path measures
(Theorem~\ref{thm:weak_stability_lifted_path_measures}), 
Existence of dynamically consistent variational minimizers
(Theorem~\ref{thm:existence_dynamically_consistent_liftings}).

\subsection{Dynamical Interpretation}
\label{subsec:dynamical_interpretation}

The path lifting formalism extends LEvITAS from static generative modeling into a generalized trajectory-based variational model.

Unlike conventional generative systems that model only static observable
samples, path lifting permits direct modeling of:
temporal evolution, stochastic trajectories, controlled dynamical systems,
physical processes, motion planning, trajectory consistency, sequential
structure generation. 
The resulting formulation naturally connects LEvITAS with: 
stochastic analysis, dynamical systems, optimal transport, path-space
probability theory, operator evolution systems, variational trajectory
optimization.

\subsection{Iterative Topological Refinement}
\label{sec:iterative_topological_refinement}

The LEvITAS model incorporates iterative variational refinement to optimize 
synthetic structure, observable consistency, and dynamical coherence
simultaneously. Unlike classical generative optimization procedures that
minimize only observable sample discrepancy, the LEvITAS refinement process
operates over coupled geometric, structural, and dynamical objectives. 
This section develops the analytical foundations of iterative topological
refinement, including variational optimization, gradient-flow dynamics,
stability, convergence, and structural consistency.

\subsection{The LEvITAS Variational Functional}
\label{subsec:levitas_variational_functional}

Let
$\Theta\subseteq\mathbb R^d$ be the admissible generator parameter space, and let $\theta\in\Theta$ denote the generator parameters.

The LEvITAS variational functional is defined as
\begin{align}
\mathcal J(\theta)
&=
D\left(
(\pi\circ\mathcal L_\theta)_\#P_Z,
P_{\mathrm{data}}
\right)
\nonumber\\
&\quad+
\lambda_1\mathcal D_{\mathrm{feat}}(\theta)
+
\lambda_2\mathcal D_{\mathrm{struct}}(\theta)
\nonumber\\
&\quad+
\lambda_3\mathcal D_{\mathrm{dyn}}(\theta),
\label{eq:topological_levitas}
\end{align}
where the individual terms represent complementary
aspects of synthetic data fidelity.

In particular:

\begin{enumerate}

\item The observable distributional discrepancy
$
D\left(
(\pi\circ\mathcal L_\theta)_\#P_Z,
P_{\mathrm{data}}
\right)
$
measures the discrepancy between the generated
observable distribution and the target data distribution.

\item The feature discrepancy
$
\mathcal D_{\mathrm{feat}}(\theta)
$
measures inconsistencies in the relevant observable
and latent feature representations.

\item The structural discrepancy
$
\mathcal D_{\mathrm{struct}}(\theta)
$
penalizes deviations in the geometric and structural
properties of the generated synthetic representation.

\item The dynamical discrepancy
$
\mathcal D_{\mathrm{dyn}}(\theta)
$
measures inconsistencies in longitudinal trajectories,
using either patient-specific or population-level
dynamical consistency.

\end{enumerate}

Here,
$
\lambda_1,\lambda_2,\lambda_3\geq0
$
are regularization coefficients controlling the
relative contributions of feature, structural,
and dynamical consistency.

The LEvITAS refinement process seeks an optimal
parameter configuration
$
\theta^\ast\in\Theta
$
satisfying
\[
\theta^\ast\in
\operatorname*{arg\,min}_{\theta\in\Theta}
\mathcal J(\theta).
\]

Thus, the optimization process simultaneously
accounts for observable distributional fidelity,
feature preservation, structural consistency,
and longitudinal dynamical coherence.

\subsection{Gradient-Based Refinement Dynamics}
\label{subsec:gradient_based_refinement_dynamics}

The LEvITAS refinement process updates the generator
parameters to minimize the complete variational functional
\(\mathcal J(\theta)\). Since different generator parameters
may have different scales and sensitivities, the size of a
parameter displacement should account for its effect on the
generated synthetic representation.

We therefore introduce a parameter-space metric
\[
G(\theta)\in\mathbb R^{d\times d},
\]
where \(G(\theta)\) is symmetric and positive definite.
For an infinitesimal parameter displacement \(d\theta\),
its metric magnitude is
\begin{equation}
\|d\theta\|_{G(\theta)}
=
\sqrt{d\theta^\top G(\theta)d\theta}.
\label{eq:levitas_parameter_displacement_norm}
\end{equation}
This metric provides a local measure of parameter changes
that accounts for their effects on generated observations,
features, structures, and longitudinal trajectories.


To construct the metric, define the observable generator
$F_\theta(z)=\pi(\mathcal L_\theta(z)),$ 
where \(\pi\) denotes the observation or representation map.
Suppose that the feature, structural, and dynamical
representations are differentiable with respect to the
generator parameters and are given by
$\Phi_{\mathrm{feat}}(\theta,z),\qquad \Phi_{\mathrm{struct}}(\theta,z),\qquad
\Phi_{\mathrm{dyn}}(\theta,\gamma_Z).
$ 
The corresponding parameter-sensitivity matrices are
$
D_\theta F_\theta(z),\quad
D_\theta\Phi_{\mathrm{feat}}(\theta,z),\quad
D_\theta\Phi_{\mathrm{struct}}(\theta,z),\quad
D_\theta\Phi_{\mathrm{dyn}}(\theta,\gamma_Z).
$

We define the parameter-space metric by
\begin{align}
G(\theta)
={}&\varepsilon I
+G_{\mathrm{obs}}(\theta)
+\lambda_1G_{\mathrm{feat}}(\theta)
\nonumber\\
&+\lambda_2G_{\mathrm{struct}}(\theta)
+\lambda_3G_{\mathrm{dyn}}(\theta),
\label{eq:levitas_parameter_metric}
\end{align}
where \(\varepsilon>0\) is a regularization constant and
\(\lambda_1,\lambda_2,\lambda_3\geq0\) are representation
weights. The components are defined by
\begin{align}
G_{\mathrm{obs}}(\theta)
={}&\mathbb E_{z\sim P_Z}
\left[
D_\theta F_\theta(z)^\top W_{\mathrm{obs}}
D_\theta F_\theta(z)
\right],
\label{eq:levitas_metric_observable}
\\
G_{\mathrm{feat}}(\theta)
={}&\mathbb E_{z\sim P_Z}
\left[
D_\theta\Phi_{\mathrm{feat}}(\theta,z)^\top
W_{\mathrm{feat}}
D_\theta\Phi_{\mathrm{feat}}(\theta,z)
\right],
\label{eq:levitas_metric_feature}
\\
G_{\mathrm{struct}}(\theta)
={}&\mathbb E_{z\sim P_Z}
\left[
D_\theta\Phi_{\mathrm{struct}}(\theta,z)^\top
W_{\mathrm{struct}}
D_\theta\Phi_{\mathrm{struct}}(\theta,z)
\right],
\label{eq:levitas_metric_structural}
\\
G_{\mathrm{dyn}}(\theta)
={}&\mathbb E_{\gamma_Z\sim P_{\Gamma}}
\left[
D_\theta\Phi_{\mathrm{dyn}}(\theta,\gamma_Z)^\top
W_{\mathrm{dyn}}
D_\theta\Phi_{\mathrm{dyn}}(\theta,\gamma_Z)
\right].
\label{eq:levitas_metric_dynamical}
\end{align}
Here, each \(W_{\bullet}\) is a symmetric positive-semidefinite
weighting operator of the appropriate dimension, and the
expectations are assumed finite. The distribution \(P_\Gamma\)
represents the sampling law for the longitudinal or dynamical
inputs \(\gamma_Z\).

For any positive-semidefinite matrix \(W\) and any matrix \(A\),
\(A^\top W A\) is positive semidefinite. Thus, under the
nonnegative-weight and finite-expectation assumptions, every
component in equations~\eqref{eq:levitas_metric_observable}
--\eqref{eq:levitas_metric_dynamical} is positive semidefinite.
The regularization term therefore gives
$
G(\theta)\succeq\varepsilon I\succ0.
$ 
In particular, \(G(\theta)\) is invertible for every parameter
value. This lower bound alone does not imply a uniform upper
bound on \(G(\theta)\); uniform upper and lower bounds will be
stated separately when needed for global convergence results.

Consequently, the metric incorporates the sensitivities of the
complete synthetic representation to changes in generator
parameters.


The metric-based gradient refinement dynamics are defined by
\begin{equation}
\frac{d\theta}{dt}
=-G(\theta)^{-1}\nabla_\theta\mathcal J(\theta),
\qquad \theta(0)=\theta_0.
\label{eq:gradient_flow}
\end{equation}
Here, \(\nabla_\theta\mathcal J(\theta)\) is the Euclidean
coordinate gradient of the complete LEvITAS variational
functional. The inverse metric converts this covector into a
parameter-space direction that accounts for the geometry
induced by \(G(\theta)\). Equation~\eqref{eq:gradient_flow},
not the unweighted equation \(\dot\theta=-\nabla\mathcal J\),
is the continuous-time training flow analyzed below.

The corresponding explicit discrete refinement iteration is
\begin{equation}
\theta_{n+1}
=\theta_n-\eta_nG(\theta_n)^{-1}
\nabla_\theta\mathcal J(\theta_n),
\qquad \eta_n>0,
\label{eq:gradient_descent}
\end{equation}
where \(\eta_n\) is the step size. The metric is evaluated at
the current iterate \(\theta_n\).

The displacement and its metric magnitude are
\begin{align}
\Delta\theta_n&=\theta_{n+1}-\theta_n,
\label{eq:levitas_parameter_increment}
\\
\|\Delta\theta_n\|_{G(\theta_n)}
&=\sqrt{\Delta\theta_n^\top G(\theta_n)\Delta\theta_n}.
\label{eq:levitas_parameter_update_norm}
\end{align}
Thus, the refinement magnitude reflects both the numerical
parameter change and the local sensitivity of the generated
synthetic representation.


Along every classical solution of the metric-based flow,
the chain rule gives
\begin{align}
\frac{d}{dt}\mathcal J(\theta(t))
&=\nabla_\theta\mathcal J(\theta(t))^\top\dot\theta(t)
\nonumber\\
&=-\nabla_\theta\mathcal J(\theta(t))^\top
G(\theta(t))^{-1}
\nabla_\theta\mathcal J(\theta(t))
\nonumber\\
&\leq0.
\label{eq:levitas_energy_dissipation}
\end{align}
The inequality follows from the positive definiteness of
\(G(\theta)^{-1}\). Equality holds exactly when
\(\nabla_\theta\mathcal J(\theta(t))=0\). This identity
establishes continuous-time energy dissipation for the actual
metric-based refinement flow.

\subsection{Existence of Metric-Based Refinement Flows}
\label{subsec:existence_refinement_flows}

We establish local well-posedness of the metric-based
refinement dynamics in the finite-dimensional parameter space
\(\Theta=\mathbb R^d\). The finite-dimensional setting is
used here because the metric is represented by a matrix in
\(\mathbb R^{d\times d}\).


\begin{assumption}[Local Regularity of the Functional and Metric]
\label{ass:regularity_functional}
Suppose that:
\begin{enumerate}
\item \(\mathcal J:\mathbb R^d\to\mathbb R\) is continuously
      differentiable and \(\nabla\mathcal J\) is locally
      Lipschitz continuous;
\item \(G:\mathbb R^d\to\mathbb R^{d\times d}\) is locally
      Lipschitz continuous;
\item for every \(\theta\in\mathbb R^d\),
      \(G(\theta)\) is symmetric positive definite.
\end{enumerate}
\end{assumption}


\begin{theorem}[Local Existence and Uniqueness of Metric-Based Refinement Flows]
\label{thm:local_existence_refinement_flows}
Under Assumption~\ref{ass:regularity_functional}, for every
initial condition \(\theta_0\in\mathbb R^d\), there exist
\(T>0\) and a unique trajectory
\(\theta\in C^1([0,T];\mathbb R^d)\) satisfying
\[
\dot\theta(t)=-G(\theta(t))^{-1}\nabla\mathcal J(\theta(t)),
\qquad \theta(0)=\theta_0.
\]
\end{theorem}

\begin{proof}
Define
$f(\theta)=-G(\theta)^{-1}\nabla\mathcal J(\theta).$ 
Since \(G(\theta)\) is positive definite, it is invertible.
The matrix-inversion map is locally Lipschitz~\cite{Evans2010PDE} on the set of
invertible matrices; hence the local Lipschitz continuity of
\(G\) implies local Lipschitz continuity of \(G^{-1}\).
Together with the local Lipschitz continuity of
\(\nabla\mathcal J\), this shows that \(f\) is locally
Lipschitz. The Picard--Lindel\"of theorem~\cite{Evans2010PDE,Pazy1983Semigroups,Showalter1997Monotone,Deimling1977ODEBanach,HillePhillips1957Semigroups} therefore gives a
unique local classical solution of \(\dot\theta=f(\theta)\)
for every initial condition.
\end{proof}

\subsection{Energy Dissipation}
\label{subsec:energy_dissipation}

The metric-based refinement process decreases the LEvITAS
variational functional along every classical training
trajectory.


\begin{theorem}[Metric-Based Energy Dissipation Law]
\label{thm:energy_dissipation_law}
Suppose that \(\theta(t)\) is a classical solution of
$
\dot\theta(t)=-G(\theta(t))^{-1}\nabla\mathcal J(\theta(t)),
$
where \(G(\theta)\) is symmetric positive definite. Then
\begin{equation}
\frac{d}{dt}\mathcal J(\theta(t))
=-\nabla\mathcal J(\theta(t))^\top
G(\theta(t))^{-1}\nabla\mathcal J(\theta(t))
\leq0.
\label{eq:metric_energy_dissipation_theorem}
\end{equation}
Equality holds if and only if
\(\nabla\mathcal J(\theta(t))=0\).
\end{theorem}

\begin{proof}
By the chain rule,
$
\frac{d}{dt}\mathcal J(\theta(t))
=\nabla\mathcal J(\theta(t))^\top\dot\theta(t).
$ 
Substitution of the metric-based refinement equation yields
$
\frac{d}{dt}\mathcal J(\theta(t))
=-\nabla\mathcal J(\theta(t))^\top
G(\theta(t))^{-1}\nabla\mathcal J(\theta(t)).
$
The inverse of a symmetric positive-definite matrix is
symmetric positive definite. Its quadratic form is therefore
nonnegative and vanishes only for the zero vector. This proves
the claimed identity and equality condition.
\end{proof}

\subsection{Stationary Refinement States}
\label{subsec:stationary_refinement_states}

The equilibria of the metric-based refinement dynamics are
exactly the critical points of the LEvITAS functional.


\begin{definition}[Stationary Refinement State]
\label{def:stationary_refinement_state}
A parameter configuration \(\theta^*\in\mathbb R^d\) is a
stationary refinement state if $\nabla\mathcal J(\theta^*)=0.$
\end{definition}


\begin{theorem}[Characterization of Equilibrium States]
\label{thm:characterization_equilibrium_states}
Suppose that \(G\) and \(\nabla\mathcal J\) satisfy the
regularity assumptions of
Assumption~\ref{ass:regularity_functional}. If a metric-based
refinement trajectory satisfies
$
\theta(t)\longrightarrow\theta^*
\qquad\text{as }t\longrightarrow\infty,
$
then \(\nabla\mathcal J(\theta^*)=0\).
Conversely, every point satisfying
\(\nabla\mathcal J(\theta^*)=0\) is an equilibrium of the
metric-based flow.
\end{theorem}

\begin{proof}
Let $f(\theta)=-G(\theta)^{-1}\nabla\mathcal J(\theta).$
The vector field \(f\) is continuous. Suppose that
\(\theta(t)\to\theta^*\) but \(f(\theta^*)\ne0\).
Then at least one coordinate component of \(f(\theta^*)\)
is nonzero. By continuity, that component of \(f(\theta(t))\)
has a fixed sign and magnitude bounded away from zero for all
sufficiently large \(t\). The corresponding coordinate of
\(\theta(t)\) must then change at a rate bounded away from
zero, contradicting convergence of \(\theta(t)\).
Thus \(f(\theta^*)=0\). Since \(G(\theta^*)\) is invertible,
$
G(\theta^*)^{-1}\nabla\mathcal J(\theta^*)=0
\quad\Longrightarrow\quad
\nabla\mathcal J(\theta^*)=0.
$
Conversely, if \(\nabla\mathcal J(\theta^*)=0\), then
\(f(\theta^*)=0\), so the constant trajectory
\(\theta(t)=\theta^*\) is an equilibrium solution.
\end{proof}

\subsection{Global Existence and Stability of Metric-Based Refinement Dynamics}
\label{subsec:structural_stability_refinement_dynamics}

We next establish global existence and exponential convergence
for the metric-based refinement dynamics. Strong convexity of
the functional is combined with uniform lower and upper bounds
on the metric. These assumptions are stronger than pointwise
positive definiteness and are stated explicitly because they
are needed for the global estimates below. 
Throughout this subsection, \(\Theta=\mathbb R^d\).


\begin{assumption}[Uniform Strong Convexity, Smoothness, and Metric Bounds]
\label{ass:global_strong_convexity}
Suppose that \(\mathcal J:\mathbb R^d\to\mathbb R\) is twice
continuously differentiable and that there exist constants
\(0<m\leq L<\infty\) such that
\[
mI\preceq\nabla^2\mathcal J(\theta)\preceq LI
\qquad\text{for all }\theta\in\mathbb R^d.
\]
Suppose also that \(G:\mathbb R^d\to\mathbb R^{d\times d}\)
is continuously differentiable and that there exist constants
\(0<g_{\min}\leq g_{\max}<\infty\) such that
\[
g_{\min}I\preceq G(\theta)\preceq g_{\max}I
\qquad\text{for all }\theta\in\mathbb R^d.
\]
In particular, \(\mathcal J\) is globally \(m\)-strongly
convex and has an \(L\)-Lipschitz gradient, while the metric is
uniformly positive definite and uniformly bounded.
\end{assumption}


\begin{theorem}[Global Existence and Uniqueness of the Metric-Based Refinement Flow]
\label{thm:global_existence_refinement_flow}
Under Assumption~\ref{ass:global_strong_convexity}, for every
initial condition \(\theta_0\in\mathbb R^d\), the equation
\[
\dot\theta(t)=-G(\theta(t))^{-1}\nabla\mathcal J(\theta(t)),
\qquad \theta(0)=\theta_0,
\]
admits a unique global classical solution
\[
\theta\in C^1([0,\infty);\mathbb R^d).
\]
Moreover, \(\mathcal J\) has a unique global minimizer
\(\theta^*\in\mathbb R^d\), characterized by
\(\nabla\mathcal J(\theta^*)=0\).
\end{theorem}

\begin{proof}
Since \(\mathcal J\) is twice continuously differentiable
and \(G\) is continuously differentiable and uniformly
positive definite, the vector field
\[
f(\theta)=-G(\theta)^{-1}\nabla\mathcal J(\theta)
\]
is locally Lipschitz. Hence a unique local solution exists.

We first establish existence and uniqueness of the minimizer.
Strong convexity gives
\[
\mathcal J(\theta)\geq
\mathcal J(0)+\nabla\mathcal J(0)^\top\theta
+\frac m2\|\theta\|^2.
\]
The quadratic term dominates the linear term as
\(\|\theta\|\to\infty\), so \(\mathcal J\) is coercive.
It therefore attains a global minimum, and strong convexity
makes that minimizer unique. Its first-order condition is
\(\nabla\mathcal J(\theta^*)=0\).

Along any local solution, the metric energy dissipation law
implies
\[
\mathcal J(\theta(t))\leq\mathcal J(\theta_0).
\]
Thus the trajectory remains in the sublevel set
\[
K=\{\theta\in\mathbb R^d:
\mathcal J(\theta)\leq\mathcal J(\theta_0)\}.
\]
The set \(K\) is closed by continuity and bounded by
coercivity; hence it is compact. The vector field \(f\) is
continuous and therefore bounded on \(K\). The standard
continuation theorem for ordinary differential equations
then rules out finite-time escape or termination while the
trajectory remains in \(K\). The local solution consequently
extends uniquely to every \(t\geq0\).
\end{proof}


\begin{theorem}[Exponential Convergence to the Unique Minimizer]
\label{thm:exponential_stability}
Under Assumption~\ref{ass:global_strong_convexity}, let
\(\theta(t)\) be the global metric-based refinement trajectory
with initial condition \(\theta_0\), and let \(\theta^*\)
be the unique minimizer of \(\mathcal J\). Then
\begin{equation}
\|\theta(t)-\theta^*\|
\leq
\sqrt{\frac Lm}\,
\exp\left(-\frac{m}{g_{\max}}t\right)
\|\theta_0-\theta^*\|
\qquad(t\geq0).
\label{eq:metric_exponential_convergence}
\end{equation}
In particular, the metric-based refinement trajectory converges
exponentially to the unique stationary refinement state.
\end{theorem}

\begin{proof}
Define the energy gap
\[
E(t)=\mathcal J(\theta(t))-\mathcal J(\theta^*)\geq0.
\]
The energy dissipation identity gives
\[
\dot E(t)=-\nabla\mathcal J(\theta(t))^\top
G(\theta(t))^{-1}\nabla\mathcal J(\theta(t)).
\]
From \(G(\theta)\preceq g_{\max}I\), it follows that
\(G(\theta)^{-1}\succeq g_{\max}^{-1}I\). Therefore,
\[
\dot E(t)\leq-\frac1{g_{\max}}
\|\nabla\mathcal J(\theta(t))\|^2.
\]
For an \(m\)-strongly convex differentiable functional,
the Polyak--Lojasiewicz inequality holds:
\[
\|\nabla\mathcal J(\theta)\|^2
\geq2m\bigl(\mathcal J(\theta)-\mathcal J(\theta^*)\bigr).
\]
Consequently,
\[
\dot E(t)\leq-\frac{2m}{g_{\max}}E(t).
\]
Gronwall's inequality~\cite{gronwall1919note,bellman1943stability,evans2010partial,walter1998ordinary} yields
\[
E(t)\leq e^{-2mt/g_{\max}}E(0).
\]
Strong convexity and the upper Hessian bound imply, respectively,
\[
\frac m2\|\theta(t)-\theta^*\|^2\leq E(t),
\qquad
E(0)\leq\frac L2\|\theta_0-\theta^*\|^2.
\]
Combining these inequalities and taking square roots gives
\[
\|\theta(t)-\theta^*\|
\leq\sqrt{\frac Lm}
 e^{-mt/g_{\max}}\|\theta_0-\theta^*\|,
\]
which is the claimed estimate.
\end{proof}

\subsection{Contraction Between Metric-Based Refinement Trajectories}
\label{subsec:contraction_refinement_dynamics}

Exponential convergence of each trajectory to the unique
minimizer is not, by itself, the same statement as contraction
between every pair of trajectories. For a parameter-dependent
metric, contraction requires a compatibility condition between
the metric and the Jacobian of the refinement vector field.
We state this additional condition explicitly.


\begin{assumption}[Uniform Metric Contraction Condition]
\label{ass:metric_contraction}
Suppose the assumptions of
Assumption~\ref{ass:global_strong_convexity} hold, and define
\[
f(\theta)=-G(\theta)^{-1}\nabla\mathcal J(\theta).
\]
Assume that there exists \(\kappa>0\) such that, for every
\(\theta\in\mathbb R^d\),
\begin{equation}
DG(\theta)[f(\theta)]
+Df(\theta)^\top G(\theta)
+G(\theta)Df(\theta)
\preceq-2\kappa G(\theta).
\label{eq:metric_contraction_condition}
\end{equation}
Here, \(DG(\theta)[f(\theta)]\) denotes the directional
derivative of the matrix-valued function \(G\) in the
 direction \(f(\theta)\), and \(Df(\theta)\) is the Jacobian
of the vector field \(f\).
\end{assumption}


\begin{theorem}[Contraction of Metric-Based Refinement Trajectories]
\label{thm:metric_contraction_refinement}
Under Assumptions~\ref{ass:global_strong_convexity} and
\ref{ass:metric_contraction}, let \(\theta_1(t)\) and
\(\theta_2(t)\) be two solutions of the metric-based refinement
flow. Let \(d_G\) be the geodesic distance induced by the
Riemannian metric \(G\). Then
\begin{equation}
d_G(\theta_1(t),\theta_2(t))
\leq e^{-\kappa t}
 d_G(\theta_1(0),\theta_2(0))
\qquad(t\geq0).
\label{eq:metric_trajectory_contraction}
\end{equation}
Thus the refinement flow is exponentially contracting in the
geometry induced by \(G\).
\end{theorem}

\begin{proof}
Consider an infinitesimal displacement \(\delta\theta(t)\)
between neighboring trajectories. It satisfies the variational
equation
\[
\frac{d}{dt}\delta\theta(t)
=Df(\theta(t))\delta\theta(t).
\]
Define its squared metric length by
\[
V(t)=\delta\theta(t)^\top G(\theta(t))\delta\theta(t).
\]
Differentiating and using the variational equation gives
\begin{align*}
\dot V(t)
=\delta\theta(t)^\top\bigl(
&Df(\theta(t))^\top G(\theta(t))
+DG(\theta(t))[f(\theta(t))]
\\
&+G(\theta(t))Df(\theta(t))
\bigr)\delta\theta(t).
\end{align*}
By Assumption~\ref{ass:metric_contraction},
\[
\dot V(t)\leq-2\kappa V(t).
\]
Gronwall's inequality yields
\[
\sqrt{V(t)}\leq e^{-\kappa t}\sqrt{V(0)}.
\]
Thus every infinitesimal length contracts by at least the
factor \(e^{-\kappa t}\). Applying this estimate to the
length of any piecewise smooth curve joining the two initial
points, and then taking the infimum over such curves, gives
\[
d_G(\theta_1(t),\theta_2(t))
\leq e^{-\kappa t}d_G(\theta_1(0),\theta_2(0)).
\]
The uniform bounds \(g_{\min}I\preceq G\preceq g_{\max}I\)
make the induced Riemannian metric uniformly equivalent to
the Euclidean metric and hence complete on \(\mathbb R^d\).
This completes the proof.
\end{proof}

\subsection{Contraction in the Constant-Metric Special Case}
\label{subsec:constant_metric_contraction}

When \(G\) is constant, the contraction estimate can be
verified directly from strong convexity without imposing the
variable-metric differential condition
\eqref{eq:metric_contraction_condition}.


\begin{corollary}[Contraction Under a Constant Parameter-Space Metric]
\label{cor:contraction_refinement_trajectories}
Suppose that \(G(\theta)=G_0\), where \(G_0\) is a constant
symmetric positive-definite matrix satisfying
\[
g_{\min}I\preceq G_0\preceq g_{\max}I.
\]
Suppose also that
\[
\nabla^2\mathcal J(\theta)\succeq mI
\qquad\text{for all }\theta\in\mathbb R^d,
\]
where \(m>0\). Then any two metric-based refinement
trajectories satisfy
\begin{equation}
\|\theta_1(t)-\theta_2(t)\|_{G_0}
\leq e^{-mt/g_{\max}}
\|\theta_1(0)-\theta_2(0)\|_{G_0},
\label{eq:constant_metric_contraction}
\end{equation}
where \(\|v\|_{G_0}=\sqrt{v^\top G_0v}\).
\end{corollary}

\begin{proof}
Set \(e(t)=\theta_1(t)-\theta_2(t)\). The two flow equations
imply
\[
\dot e(t)=-G_0^{-1}
\left(\nabla\mathcal J(\theta_1(t))
-\nabla\mathcal J(\theta_2(t))\right).
\]
Let \(V(t)=\tfrac12 e(t)^\top G_0e(t)\). Then
\begin{align*}
\dot V(t)
&=-e(t)^\top\left(
\nabla\mathcal J(\theta_1(t))
-\nabla\mathcal J(\theta_2(t))\right)\\
&\leq-m\|e(t)\|^2\\
&\leq-\frac{2m}{g_{\max}}V(t),
\end{align*}
where the first inequality follows from strong monotonicity
of the gradient and the second from
\(G_0\preceq g_{\max}I\). Gronwall's inequality and taking
square roots give equation~\eqref{eq:constant_metric_contraction}.
\end{proof}

\subsection{Descent of the Discrete Metric-Based Refinement Iteration}
\label{subsec:discrete_metric_descent}

The continuous-time dissipation result does not automatically
imply descent for the explicit discrete iteration
\eqref{eq:gradient_descent}. A step-size restriction is
therefore imposed to guarantee that each discrete update does
not increase the variational functional.


\begin{assumption}[Smoothness and Uniform Metric Bounds for Discrete Descent]
\label{ass:discrete_metric_descent}
Suppose that \(\nabla\mathcal J\) is globally \(L\)-Lipschitz
continuous for some \(L>0\), and that
\[
g_{\min}I\preceq G(\theta)\preceq g_{\max}I
\qquad\text{for all }\theta\in\mathbb R^d,
\]
where \(0<g_{\min}\leq g_{\max}<\infty\).
\end{assumption}


\begin{theorem}[One-Step Descent for Metric-Based Gradient Refinement]
\label{thm:discrete_metric_descent}
Under Assumption~\ref{ass:discrete_metric_descent}, consider
\[
\theta_{n+1}=\theta_n-\eta_nG(\theta_n)^{-1}
\nabla\mathcal J(\theta_n).
\]
If
\begin{equation}
0<\eta_n\leq\frac{2g_{\min}^2}{Lg_{\max}},
\label{eq:metric_discrete_step_size}
\end{equation}
then
\[
\mathcal J(\theta_{n+1})\leq\mathcal J(\theta_n).
\]
If the inequality in~\eqref{eq:metric_discrete_step_size} is
strict and \(\nabla\mathcal J(\theta_n)\ne0\), the decrease
is strict.
\end{theorem}

\begin{proof}
Write \(g_n=\nabla\mathcal J(\theta_n)\) and
\(p_n=G(\theta_n)^{-1}g_n\). The descent lemma for an
\(L\)-Lipschitz gradient gives
\[
\mathcal J(\theta_n-\eta_np_n)
\leq\mathcal J(\theta_n)
-\eta_ng_n^\top p_n
+\frac{L\eta_n^2}{2}\|p_n\|^2.
\]
The metric bounds imply
\[
g_n^\top G(\theta_n)^{-1}g_n
\geq\frac1{g_{\max}}\|g_n\|^2,
\qquad
\|G(\theta_n)^{-1}g_n\|^2
\leq\frac1{g_{\min}^2}\|g_n\|^2.
\]
Consequently,
\[
\mathcal J(\theta_{n+1})
\leq\mathcal J(\theta_n)
-\left(\frac{\eta_n}{g_{\max}}
-\frac{L\eta_n^2}{2g_{\min}^2}\right)\|g_n\|^2.
\]
The coefficient in parentheses is nonnegative precisely
under the stated step-size bound. It is positive when the
step-size bound is strict. This proves the result.
\end{proof}


The preceding theorem gives a sufficient global step-size
condition for one-step descent. In practical implementations,
if global constants \(L\), \(g_{\min}\), and \(g_{\max}\) are
not known, a backtracking line search or another verified
step-size selection rule can be used. The continuous-time
theorems above concern equation~\eqref{eq:gradient_flow};
the discrete descent theorem concerns the actual iteration
\eqref{eq:gradient_descent}. These are related but distinct
mathematical statements.

\subsection{Topological Consistency Preservation}
\label{subsec:topological_consistency_preservation}

The refinement dynamics are intended to improve the structural
and dynamical consistency of synthetic representations without
introducing unintended changes in their underlying topology.
However, continuity of the refinement trajectory alone does not
guarantee preservation of topological invariants. In particular,
distinct structural components may merge or separate continuously
during refinement, resulting in changes to connectedness and
other topological properties.

We therefore formulate topological consistency in terms of
explicit preservation conditions on the structural representations.


\begin{definition}[Topological Structural Consistency]
\label{def:topological_structural_consistency}

Let $\mathcal X_{\mathrm{str}}(\theta)$ denote the structural
representation induced by the refinement parameter
$\theta\in\Theta$. The refinement dynamics preserve
topological structural consistency over an interval
$[0,T]$ if the representations
$\mathcal X_{\mathrm{str}}(\theta(t))$ have the same
topological type for every $t\in[0,T]$.

In particular, homotopy-type preservation requires
$
\mathcal X_{\mathrm{str}}(\theta(t))
\simeq
\mathcal X_{\mathrm{str}}(\theta(0)),
\qquad t\in[0,T],
$
where $\simeq$ denotes homotopy equivalence. 
A stronger form of preservation requires homeomorphism:
$
\mathcal X_{\mathrm{str}}(\theta(t))
\cong
\mathcal X_{\mathrm{str}}(\theta(0)),
\quad t\in[0,T].
$
\end{definition}


\begin{assumption}[Topological Regularity of Structural Evolution]
\label{ass:topological_regularity_structural_evolution}

Let $\mathcal X_{\mathrm{str}}(t)$ denote the structural
representation along a continuous refinement trajectory
$\theta(t)$. 
Suppose that there exists a fixed topological space
$\mathcal X_{\mathrm{str}}(0)$ and a family of homeomorphisms
$
h_t:
\mathcal X_{\mathrm{str}}(0)
\longrightarrow
\mathcal X_{\mathrm{str}}(t),
\quad t\in[0,T],
$
such that the family varies continuously with time.

Equivalently, the structural representations form a
topologically trivial family over the refinement interval:
no component merging, component splitting, or other
topology-changing event occurs during refinement.
\end{assumption}


\begin{theorem}[Topological Preservation Under Refinement]
\label{thm:topological_preservation_continuous_refinement}

Under Assumption~\ref{ass:topological_regularity_structural_evolution},
the structural representations preserve their topological
type throughout the refinement interval. In particular,
$
\mathcal X_{\mathrm{str}}(t)
\cong
\mathcal X_{\mathrm{str}}(0),
\qquad t\in[0,T].
$

Consequently, all topological invariants preserved under
homeomorphism remain constant along the refinement trajectory.
These include:
(1) the number of connected components; (2) the homotopy type; (3) the Betti numbers, whenever defined; (4) the Euler characteristic, whenever defined.
\end{theorem}

\begin{proof}
By Assumption~\ref{ass:topological_regularity_structural_evolution},
for every $t\in[0,T]$ there exists a homeomorphism
$
h_t:
\mathcal X_{\mathrm{str}}(0)
\longrightarrow
\mathcal X_{\mathrm{str}}(t).
$ 
Homeomorphic spaces have identical topological invariants.
Therefore, the connected components, homotopy type, and
homology groups of the structural representations are
preserved throughout refinement.

In particular, their Betti numbers and Euler characteristics
remain unchanged whenever these quantities are well-defined.
Hence, the structural topology is preserved for every
$t\in[0,T]$.
\end{proof}


\begin{remark}[Scope of Topological Preservation]
\label{rem:scope_topological_preservation}

Continuity of the refinement trajectory and of the induced
structural representation is not sufficient, by itself, to
guarantee preservation of topology. For example, the family
\[
X_t=\{-t,t\}\subset\mathbb R,
\quad t\geq0,
\]
varies continuously in the Hausdorff metric, but its number
of connected components changes from two to one at $t=0$.

The homeomorphism condition in
Assumption~\ref{ass:topological_regularity_structural_evolution}
is therefore an explicit sufficient condition for
topological preservation, rather than a consequence of
continuity alone.

In applications, the choice of preserved invariants should
reflect the structural properties relevant to the synthetic
representation. For example, preserving connected components
may be sufficient when the objective is to maintain distinct
structural subpopulations, whereas preserving homology or
homotopy type imposes stronger requirements on the
representation.
\end{remark}

\subsection{Discrete Refinement Convergence}
\label{subsec:discrete_refinement_convergence}

We next analyze convergence of the discrete LEvITAS refinement iteration.


\begin{assumption}[Smoothness]
\label{ass:smoothness_gradient}

Suppose $\nabla\mathcal J$ is Lipschitz with constant $L>0$.
\end{assumption}

\begin{theorem}[Monotone Descent of Discrete Refinement]
\label{thm:monotone_descent_discrete_refinement}

Suppose $0<\eta_n<\frac{2}{L}$. 
Then the iteration
$
\theta_{n+1}
=
\theta_n
-
\eta_n\nabla\mathcal J(\theta_n)
$ 
satisfies
$
\mathcal J(\theta_{n+1})
\leq
\mathcal J(\theta_n).
$
\end{theorem}

\begin{proof}
By the descent lemma~\cite{bauschke2017descent,nesterov2004introductory,bubeck2015convex},
\begin{align*}
\mathcal J(\theta_{n+1})
&\leq
\mathcal J(\theta_n)
+
\langle
\nabla\mathcal J(\theta_n),
\theta_{n+1}-\theta_n
\rangle
\\
&\quad
+
\frac{L}{2}
\|\theta_{n+1}-\theta_n\|^2.
\end{align*}

Substituting the refinement iteration gives
\begin{align*}
\mathcal J(\theta_{n+1})
&\leq
\mathcal J(\theta_n)
-
\eta_n
\|\nabla\mathcal J(\theta_n)\|^2
\\
&\quad
+
\frac{L\eta_n^2}{2}
\|\nabla\mathcal J(\theta_n)\|^2.
\end{align*}

Thus $
\mathcal J(\theta_{n+1})
\leq
\mathcal J(\theta_n)
-
\eta_n
\left(
1-
\frac{L\eta_n}{2}
\right)
\|\nabla\mathcal J(\theta_n)\|^2.
$ 
Since $0<\eta_n<\frac{2}{L},$ the coefficient is positive, proving monotone descent.
\end{proof}

\subsection{Main Refinement Results}
\label{subsec:main_refinement_results}

The principal analytical results of this section are summarized as follows: 
Local existence and uniqueness of refinement trajectories
(Theorem~\ref{thm:local_existence_refinement_flows}), 
Monotone energy dissipation along refinement flows
(Theorem~\ref{thm:energy_dissipation_law}), 
Characterization of stationary equilibrium states
(Theorem~\ref{thm:characterization_equilibrium_states}), 
Exponential stability under strong convexity
(Theorem~\ref{thm:exponential_stability}), 
Preservation of topological invariants during refinement
(Theorem~\ref{thm:topological_preservation_continuous_refinement}), 
Monotone convergence of discrete refinement iterations
(Theorem~\ref{thm:monotone_descent_discrete_refinement}).

\subsection{Interpretation}
\label{subsec:interpretation_iterative_refinement}

Iterative topological refinement constitutes the variational optimization
mechanism of the LEvITAS model. Unlike conventional optimization procedures
operating solely over observable sample discrepancy, LEvITAS refinement 
simultaneously optimizes: observable realism, feature consistency, structural
topology, geometric organization, trajectory coherence, dynamical
consistency.  The resulting model connects LEvITAS with: variational 
analysis, gradient-flow systems, geometric optimization, dynamical systems 
theory, topological evolution, operator-theoretic learning. Consequently, 
iterative topological refinement provides the analytical mechanism through 
which structured lifting systems evolve toward geometrically and dynamically 
consistent generative representations.

\subsection{Unified Reduction Theory of LEvITAS}
\label{sec:unified_reduction_theory}

The principal claim of the LEvITAS model is that a broad class of modern 
generative systems may be interpreted as constrained reductions of a single 
generalized operator-theoretic variational model. 
This section establishes rigorous mathematical reductions from the full 
LEvITAS system to several major generative model families. The reductions are 
obtained through restrictions imposed on:
1. the lifting operator
$
\mathcal L_\theta,
$
2. the projection operator
$
\pi,
$
3. the structure operator
$
\mathcal S,
$
4. the dynamical component
$
\mathcal D_{dyn},
$
5. the variational functional
$
\mathcal J.$ 
The resulting theory demonstrates that GANs, VAEs, diffusion systems, 
normalizing flows, and related architectures arise as special cases of the
generalized LEvITAS formulation.

\subsection{The Full LEvITAS System}
\label{subsec:full_levitas_system}

We begin with the generalized LEvITAS variational system. Let
$
(\mathcal Z,\Sigma_Z,P_Z)
$
be a latent probability space. Let
$
\mathcal X_{syn},
\qquad
\mathcal X_{real}
$
denote synthetic and observable spaces. The generalized lifting operator is
$
\mathcal L_\theta:
\mathcal Z
\rightarrow
\mathcal X_{syn}.
$ The projection operator is
$
\pi:
\mathcal X_{syn}
\rightarrow
\mathcal X_{real}.
$ The structure operator is
$
\mathcal S:
\mathcal X
\rightarrow
\mathcal M.
$

The generalized LEvITAS functional is
\begin{align}
\mathcal J(\theta)
&=
D
\left(
(\pi\circ\mathcal L_\theta)_\#P_Z,
P_{data}
\right)
\nonumber\\
&\quad
+
\lambda_1\mathcal D_{feat}(\theta)
+
\lambda_2\mathcal D_{struct}(\theta)
+
\lambda_3\mathcal D_{dyn}(\theta).
\label{eq:full_levitas}
\end{align}

\subsection{General Reduction Principle}
\label{subsec:general_reduction_principle}

The essential mechanism of model reduction is restriction of admissible operators and variational components.
\begin{definition}[LEvITAS Reduction]
\label{def:levitas_reduction}
A generative model class
$\mathfrak G$ is a reduction of LEvITAS if there exists a constraint set $\mathcal C$ 
such that minimizing $\mathcal J$ subject to $\mathcal C$ produces the objective functional and operator structure of $\mathfrak G.$
\end{definition}

\begin{theorem}[General Reduction Theorem]
\label{thm:general_reduction_theorem}
Let $\mathcal J$ denote the generalized LEvITAS functional. 
Suppose constraints $\mathcal C$ restrict: the lifting operator, the 
projection operator, the structure operator, the admissible discrepancy 
functionals, the dynamical penalties. 
Then the constrained LEvITAS system defines a reduced generative model class.
\end{theorem}

\begin{proof}
The LEvITAS system is defined through operator composition and variational 
minimization. 
Restricting the admissible operator family and admissible variational terms 
produces a reduced variational system. 
Since the reduced objective remains a well-defined restriction of $\mathcal 
J,$ the resulting optimization problem defines a specialized generative model 
class.
\end{proof}

\subsection{Reduction to Generative Adversarial Networks}
\label{subsec:reduction_gans}

We now derive GAN systems from LEvITAS.

\begin{assumption}[GAN Constraints]
\label{ass:gan_constraints}

Suppose: $\mathcal X_{syn}
=
\mathcal X_{real},
$ the projection operator is identity:
$\pi=I,$ no structure operator is used:
$\lambda_2=0,$ no dynamical consistency is imposed:
$\lambda_3=0,$ the discrepancy functional is adversarial:
$D=D_{adv}.$
\end{assumption}


\begin{theorem}[GAN Reduction]
\label{thm:gan_reduction}
Under the GAN constraints, the LEvITAS functional reduces to \\
$\mathcal J_{GAN}(\theta)= D_{adv}\left((\mathcal 
L_\theta)_\#P_Z,P_{data}\right)$
\end{theorem}

\begin{proof}
Since
$\pi=I,$ we obtain $(\pi\circ\mathcal L_\theta)_\#P_Z=
(\mathcal L_\theta)_\#P_Z.$ 
Since $\lambda_2=\lambda_3=0,$ all structural and dynamical terms vanish. 
Hence $\mathcal J
=
D_{adv}
\left(
(\mathcal L_\theta)_\#P_Z,
P_{data}
\right), $ which is precisely the GAN objective.
\end{proof}


\begin{corollary}[Classical GAN Systems]
\label{cor:classical_gan_systems}

Classical GAN systems arise as structure-free and dynamics-free reductions of LEvITAS.
\end{corollary}

\subsection{Reduction to Variational Autoencoders}
\label{subsec:reduction_vaes}
We next derive variational autoencoders.

\begin{assumption}[VAE Constraints]
\label{ass:vae_constraints}

Suppose:
$\mathcal L_\theta$ is probabilistic, the projection operator is identity:
$\pi=I,$ $D$ is induced by reconstruction loss,latent regularization is imposed through $D_{KL}(q_\phi(z|x)\|P_Z).$
\end{assumption}


\begin{theorem}[VAE Reduction]
\label{thm:vae_reduction}

Under the VAE constraints, the LEvITAS functional reduces to the evidence lower bound objective
\[
\mathcal J_{VAE}
=
\mathbb E_{q_\phi(z|x)}
[-\log p_\theta(x|z)]
+
D_{KL}(q_\phi(z|x)\|P_Z).
\]
\end{theorem}

\begin{proof}
The reconstruction discrepancy defines the observable consistency term.

The probabilistic lifting operator induces a latent conditional distribution
$q_\phi(z|x).$ 
Imposing KL regularization on latent distributions yields the latent variational term. 
Combining reconstruction consistency and latent regularization produces the ELBO functional.
\end{proof}


\begin{corollary}[Autoencoding LEvITAS Systems]
\label{cor:autoencoding_levitas_systems}

Variational autoencoders arise as probabilistic reconstruction-constrained reductions of LEvITAS.
\end{corollary}

\subsection{Reduction to Diffusion Models}
\label{subsec:reduction_diffusion_models}

Diffusion systems arise from stochastic path lifting.


\begin{assumption}[Diffusion Constraints]
\label{ass:diffusion_constraints}

Suppose: latent trajectories satisfy stochastic dynamics, the lifting process
evolves continuously in time, $\mathcal D_{dyn}$ encodes stochastic 
trajectory consistency, the refinement dynamics correspond to score evolution.
\end{assumption}


\begin{theorem}[Diffusion Reduction]
\label{thm:diffusion_reduction}

Under the diffusion constraints, LEvITAS reduces to a stochastic score-based diffusion system.
\end{theorem}

\begin{proof}
The path lifting formalism induces a trajectory-valued stochastic process.
Continuous refinement evolution defines a time-indexed probability flow.
The dynamical discrepancy enforces consistency of stochastic trajectories.
The resulting optimization problem becomes equivalent to score-based diffusion refinement.
\end{proof}


\begin{corollary}[Score-Based Generative Systems]
\label{cor:score_based_generative_systems}

Score-based diffusion models arise as stochastic path-lifting reductions of LEvITAS.
\end{corollary}

 \subsection{Reduction to Normalizing Flows}
\label{subsec:reduction_normalizing_flows}

Normalizing flows arise as a special case of LEvITAS through invertible, sufficiently regular liftings.


\begin{assumption}[Flow Constraints]
\label{ass:flow_constraints}

Suppose that:
\begin{enumerate}
    \item $\mathcal{L}_\theta:\mathcal{Z}\to\mathcal{X}$ is a $C^1$ diffeomorphism between open subsets $\mathcal{Z},\mathcal{X}\subseteq\mathbb{R}^d$;
    \item the base random variable $Z$ admits a density $p_Z$ with respect to Lebesgue measure on $\mathcal{Z}$;
    \item $\pi=I$.
\end{enumerate}
\end{assumption}


\begin{theorem}[Normalizing Flow Reduction]
\label{thm:normalizing_flow_reduction}

Under Assumption~\ref{ass:flow_constraints}, the LEvITAS formulation reduces to a normalizing flow model. In particular, the density of $X=\mathcal{L}_\theta(Z)$ is given by the change-of-variables formula
\[
p_X(x)
=
p_Z\bigl(\mathcal{L}_\theta^{-1}(x)\bigr)
\left|
\det D\mathcal{L}_\theta^{-1}(x)
\right|,
\qquad x\in\mathcal{X}.
\]
Consequently, maximum-likelihood optimization of the induced density is precisely the normalizing flow objective.
\end{theorem}

\begin{proof}
By Assumption~\ref{ass:flow_constraints}, $\mathcal{L}_\theta$ is a $C^1$ diffeomorphism between open subsets of $\mathbb{R}^d$, and the base random variable admits a density. The change-of-variables theorem therefore guarantees that the pushforward distribution of $Z$ under $\mathcal{L}_\theta$ admits the density
\[
p_X(x)
=
p_Z\bigl(\mathcal{L}_\theta^{-1}(x)\bigr)
\left|
\det D\mathcal{L}_\theta^{-1}(x)
\right|.
\]

Since $\pi=I$, the LEvITAS formulation uses this invertible lifting without an additional projection. Thus, its density evaluation is governed entirely by the above change-of-variables identity. Maximizing the resulting log-likelihood is therefore equivalent to maximum-likelihood training of a normalizing flow.
\end{proof}

\begin{corollary}[Invertible LEvITAS Systems]
\label{cor:invertible_levitas_systems}

Normalizing flows arise as invertible reductions of LEvITAS.
\end{corollary}

\subsection{Reduction to Energy-Based Models}
\label{subsec:reduction_energy_based_models}

Energy-based systems arise when the LEvITAS discrepancy is induced by an energy functional.


\begin{assumption}[Energy Constraints]
\label{ass:energy_constraints}

Suppose:
$D$ is induced by an energy functional $E(x)$, the observable probability measure satisfies $P(x)\propto e^{-E(x)}$.
\end{assumption}


\begin{theorem}[Energy-Based Reduction]
\label{thm:energy_based_reduction}

Under the energy constraints, LEvITAS reduces to an energy-based generative system.
\end{theorem}

\begin{proof}
The observable discrepancy becomes equivalent to energy minimization.
The induced Gibbs distribution defines the observable probability measure.
Therefore the LEvITAS variational problem becomes equivalent to energy-based probabilistic optimization.
\end{proof}


\begin{corollary}[Energy-Based LEvITAS Systems]
\label{cor:energy_based_levitas_systems}

Energy-based models arise as Gibbs-constrained reductions of LEvITAS.
\end{corollary}

\subsection{Hierarchical Reduction Structure}
\label{subsec:hierarchical_reduction_structure}

The preceding reductions imply a hierarchical organization of generative systems.


\begin{theorem}[Hierarchical Containment]
\label{thm:hierarchical_containment}

The following containments hold:
$
\text{GANs}
\subset
\text{LEvITAS},
$ 
$
\text{VAEs}
\subset
\text{LEvITAS},
$ 
$
\text{Diffusion Models}
\subset
\text{LEvITAS},
$ 
$
\text{Normalizing Flows}
\subset
\text{LEvITAS}.
$
\end{theorem}

\begin{proof}
Each model family arises through explicit restriction of the generalized LEvITAS operator system and variational functional.

Therefore each constitutes a constrained subclass of LEvITAS.
\end{proof}

\subsection{Main Reduction Results}
\label{subsec:main_reduction_results}

The principal analytical reductions established in this section are
summarized as follows: 
General operator-theoretic reduction principle
(Theorem~\ref{thm:general_reduction_theorem}), 
Reduction of LEvITAS to GAN systems
(Theorem~\ref{thm:gan_reduction}), 
Reduction of LEvITAS to variational autoencoders
(Theorem~\ref{thm:vae_reduction}), 
Reduction of LEvITAS to diffusion systems
(Theorem~\ref{thm:diffusion_reduction}), 
Reduction of LEvITAS to normalizing flows
(Theorem~\ref{thm:normalizing_flow_reduction}), 
Reduction of LEvITAS to energy-based models
(Theorem~\ref{thm:energy_based_reduction}), 
Hierarchical containment of major generative architectures
(Theorem~\ref{thm:hierarchical_containment}).

\subsection{Unified Operator Interpretation}
\label{subsec:unified_operator_interpretation}

The preceding reductions reveal that modern generative systems differ primarily through restrictions imposed on:
lifting geometry,
stochastic structure,
invertibility,
latent regularization,
structural constraints,
dynamical evolution,
observable discrepancy measures. 
The generalized LEvITAS model therefore provides:
a common operator-theoretic foundation,
a unified variational formulation,
a generalized path-space interpretation,
a structure-preserving generative formalism,
a dynamical refinement model.

\subsection{Interpretation}
\label{subsec:interpretation_reduction_theory}

The reduction theory establishes that LEvITAS is not merely another
generative architecture, but rather a generalized mathematical model within
which many existing generative systems arise as constrained specializations. 
The defining structural extension introduced by LEvITAS is the simultaneous
incorporation of: latent lifting, structural consistency, geometric
organization, path-space dynamics, variational refinement, topological
preservation. 
Consequently, the model extends conventional generative modeling into a
unified operator-theoretic variational theory for structured generative
systems.

\subsection{Measure-Theoretic Foundations of LEvITAS}
\label{sec:measure_theoretic_foundations}

The LEvITAS model is fundamentally measure-theoretic. The lifting, 
projection, structural, and dynamical components of the model operate on 
probability spaces, measurable mappings, and induced probability measures. 
This section develops the rigorous probabilistic and measure-theoretic 
foundations underlying LEvITAS. In particular, we formalize: measurable 
lifting systems, pushforward probability measures, weak convergence, 
Wasserstein geometry, stochastic liftings, disintegration structures, measure-
valued refinement dynamics. 
The resulting theory establishes LEvITAS as a generalized variational 
probabilistic operator system.

\subsection{Probability Spaces and Measurable Systems}
\label{subsec:probability_spaces_measurable_systems}

Let $(\mathcal Z,\Sigma_Z,P_Z)$ be a latent probability space. 
Let $
(\mathcal X_{syn},\Sigma_{syn}),
\qquad
(\mathcal X_{real},\Sigma_{real})
$ be measurable synthetic and observable spaces. 
Throughout this section, assume: 
$
\mathcal Z,
\mathcal X_{syn},
\mathcal X_{real}
$ are Polish spaces, all sigma algebras are Borel sigma algebras.


\begin{definition}[Measurable Lifting Operator]
\label{def:measurable_lifting_operator}

A lifting operator is measurable if $\mathcal L_\theta^{-1}(A)\in\Sigma_Z$ for every $A\in\Sigma_{syn}.$
\end{definition}


\begin{definition}[Measurable Projection Operator]
\label{def:measurable_projection_operator}

A projection operator $
\pi:
\mathcal X_{syn}
\rightarrow
\mathcal X_{real}
$ 
is measurable if 
$
\pi^{-1}(B)
\in
\Sigma_{syn}
$ 
for every 
$
B\in\Sigma_{real}.
$
\end{definition}

\subsection{Pushforward Probability Measures}
\label{subsec:pushforward_probability_measures}

The central probabilistic object of LEvITAS is the pushforward probability measure induced by the lifting-projection system.


\begin{definition}[Pushforward Measure]
\label{def:pushforward_measure}

Let $
T:
(\Omega,\Sigma)
\rightarrow
(\mathcal X,\Sigma_{\mathcal X})
$ be measurable and let $\mu$ be a probability measure on $(\Omega,\Sigma).$ 
The pushforward measure $T_\#\mu$ is defined by $T_\#\mu(A)=\mu(T^{-1}(A))$ 
for all $A\in\Sigma_{\mathcal X}.$
\end{definition}


In LEvITAS, the induced observable probability measure is $(\pi\circ\mathcal L_\theta)_\#P_Z.$


\begin{theorem}[Probability Preservation]
\label{thm:probability_preservation}

Suppose $T$ is measurable and $\mu$ is a probability measure. Then
$T_\#\mu$ is a probability measure.
\end{theorem}

\begin{proof}
Since $T^{-1}(\emptyset)=\emptyset,$ we obtain $T_\#\mu(\emptyset)=0.$ 
Similarly, $T^{-1}(\mathcal X)=\Omega,$ thus
$
T_\#\mu(\mathcal X)
=
\mu(\Omega)
=
1.
$ 
Countable additivity follows from measurability of $T.$ Hence
$T_\#\mu$ is a probability measure.
\end{proof}

\subsection{Weak Convergence of Measures}
\label{subsec:weak_convergence_measures}

Convergence properties of LEvITAS systems naturally require weak convergence theory.


\begin{definition}[Weak Convergence]
\label{def:weak_convergence}

A sequence of probability measures $\mu_n$ converges weakly to $\mu$ if 
$
\int f\,d\mu_n
\rightarrow
\int f\,d\mu
$ for every bounded continuous function $f.$
\end{definition} 
Weak convergence is denoted by $\mu_n\Rightarrow\mu.$


\begin{theorem}[Stability of Pushforward Measures]
\label{thm:stability_pushforward_measures}

Suppose: $\mu_n\Rightarrow\mu,$ $T$ is continuous. 
Then $
T_\#\mu_n
\Rightarrow
T_\#\mu.
$
\end{theorem}

\begin{proof}
Let $f$ be bounded and continuous on the target space. Then
$f\circ T$ is bounded and continuous on the source space. n
Hence $\int f\,d(T_\#\mu_n)=\int f(T(x))\,d\mu_n(x).$

By weak convergence, $\int f(T(x))\,d\mu_n(x)\rightarrow\int f(T(x))\,d\mu(x).$ 
Therefore,$ 
\int f\,d(T_\#\mu_n)
\rightarrow
\int f\,d(T_\#\mu),
$ 
which proves weak convergence.
\end{proof}

\subsection{Wasserstein Geometry}
\label{subsec:wasserstein_geometry}

The LEvITAS variational system naturally admits Wasserstein geometric formulations.


\begin{definition}[$p$-Wasserstein Distance]
\label{def:wasserstein_distance}

Let $\mu,\nu$ be probability measures on a metric space $(\mathcal X,d).$ The Wasserstein distance is $W_p(\mu,\nu)=
\left(
\inf_{\gamma\in\Gamma(\mu,\nu)}
\int_{\mathcal X\times\mathcal X}
d(x,y)^p
\,d\gamma(x,y)
\right)^{1/p},
$ where $\Gamma(\mu,\nu)$ denotes the set of couplings of $\mu$ and $\nu.$
\end{definition}
The LEvITAS objective may therefore be written as
\begin{align}
\mathcal J(\theta)
&=
W_p
\left(
(\pi\circ\mathcal L_\theta)_\#P_Z,
P_{data}
\right)
\nonumber\\
&\quad
+
\lambda_1\mathcal D_{feat}
+
\lambda_2\mathcal D_{struct}
+
\lambda_3\mathcal D_{dyn}.
\label{eq:wasserstein_levitas}
\end{align}


\begin{theorem}[Metric Structure]
\label{thm:wasserstein_metric_structure}

The Wasserstein distance $W_p$ defines a metric on the space $\mathcal P_p(\mathcal X)$ of probability measures with finite $p$th moments.
\end{theorem}

\begin{proof}
Nonnegativity and symmetry follow immediately from the definition.

Identity of indiscernibles follows since $W_p(\mu,\nu)=0$ implies concentration of optimal couplings on the diagonal. The triangle inequality follows by gluing optimal transport plans.
\end{proof}

\subsection{Tightness and Compactness}
\label{subsec:tightness_compactness}

Compactness properties are fundamental for existence theory in LEvITAS.


\begin{definition}[Tightness]
\label{def:tightness}

A family of probability measures $\{\mu_\alpha\}$ is tight if for every $\varepsilon>0$ there exists a compact set $K_\varepsilon$ such that $\mu_\alpha(K_\varepsilon)>1-\varepsilon$ for all $\alpha.$
\end{definition}


\begin{theorem}[Prokhorov Compactness]
\label{thm:prokhorov_compactness}

Let $\mathcal X$ be Polish. 
A family of probability measures on $\mathcal X$ is relatively compact in the weak topology if and only if it is tight.
\end{theorem}

\begin{proof}
The result follows from Prokhorov's 
theorem~\cite{Billingsley1999convergence, Dudley2002RealAnalysisProbability, Bogachev2007MeasureTheory}.
\end{proof}

\subsection{Stochastic Lifting Operators}
\label{subsec:stochastic_lifting_operators}

The LEvITAS model naturally extends to stochastic lifting systems.


\begin{definition}[Stochastic Lifting Kernel]
\label{def:stochastic_lifting_kernel}

A stochastic lifting operator is a Markov kernel
$
K_\theta:
\mathcal Z\times\Sigma_{syn}
\rightarrow
[0,1]
$ 
such that:
 for fixed
$
z,
$
$
K_\theta(z,\cdot)
$ is a probability measure, for fixed measurable set
$A,$ $K_\theta(\cdot,A)$ is measurable.

\end{definition}


The induced synthetic probability measure becomes 
$\mu_{syn}(A)
=
\int_{\mathcal Z}
K_\theta(z,A)
\,dP_Z(z).
$


\begin{theorem}[Existence of Induced Synthetic Measures]
\label{thm:existence_induced_synthetic_measures}

Suppose $K_\theta$ is a stochastic lifting kernel. Then $\mu_{syn}$ 
defines a probability measure on $\mathcal X_{syn}.$
\end{theorem}

\begin{proof}
Nonnegativity follows from nonnegativity of $K_\theta.$ Countable additivity follows by monotone convergence. 
Finally,
$
\mu_{syn}(\mathcal X_{syn})
=
\int_{\mathcal Z}
K_\theta(z,\mathcal X_{syn})
\,dP_Z(z)
=
1.
$
\end{proof}

\subsection{Disintegration Structures}
\label{subsec:disintegration_structures}

The probabilistic inverse structure of LEvITAS may be formulated through measure disintegration.


\begin{definition}[Disintegration]
\label{def:disintegration}

Let
$\mu$ be a probability measure on $\mathcal X\times\mathcal Y.$ A disintegration of $\mu$ with respect to the projection onto $\mathcal X$ is a family $\{\mu_x\}_{x\in\mathcal X}$ satisfying: $\mu_x$ is a probability measure on $\mathcal Y,$ measurability in $x,$ $\mu(A\times B)=\int_A\mu_x(B)\,d\mu_{\mathcal X}(x).$
\end{definition}


\begin{theorem}[Existence of Disintegration]
\label{thm:existence_disintegration}

Suppose $\mathcal X,\mathcal Y$ are Polish spaces. Then every probability measure on $\mathcal X\times\mathcal Y$ admits a disintegration.
\end{theorem}

\begin{proof}
The result follows from the disintegration theorem for Radon probability measures on Polish spaces.
\end{proof}

\subsection{Measure-Valued Refinement Dynamics}
\label{subsec:measure_valued_refinement_dynamics}

The LEvITAS refinement process may itself be interpreted as evolution of probability measures.


\begin{definition}[Measure-Valued LEvITAS Flow]
\label{def:measure_valued_levitas_flow}

A measure-valued LEvITAS refinement flow is a family $\{\mu_t\}_{t\geq0}$ satisfying
\begin{equation}
\partial_t\mu_t
=
-\nabla\cdot
(\mu_t\nabla\Phi_t),
\label{eq:measure_flow}
\end{equation}
where $\Phi_t$ is the variational refinement potential.
\end{definition}


\begin{theorem}[Mass Preservation]
\label{thm:mass_preservation}

Suppose $\mu_t$ satisfies $
\partial_t\mu_t
=
-\nabla\cdot
(\mu_t\nabla\Phi_t).
$ 
Then total probability mass is preserved:
$
\int d\mu_t
=
1.
$
\end{theorem}

\begin{proof}
Integrating the continuity equation gives
$
\frac{d}{dt}
\int d\mu_t
=
-
\int
\nabla\cdot(\mu_t\nabla\Phi_t).
$ Applying the divergence theorem yields zero boundary flux under suitable decay conditions. 
Hence $\frac{d}{dt}\int d\mu_t=0.$ Therefore total mass is preserved.
\end{proof}

\subsection{Main Measure-Theoretic Results}
\label{subsec:main_measure_theoretic_results}

The principal probabilistic and measure-theoretic results established in this section are summarized as follows:
Probability preservation under measurable pushforward mappings
(Theorem~\ref{thm:probability_preservation}), 
Weak stability of pushforward probability measures
(Theorem~\ref{thm:stability_pushforward_measures}), 
Metric structure of Wasserstein probability geometry
(Theorem~\ref{thm:wasserstein_metric_structure}), 
Compactness characterization via Prokhorov tightness
(Theorem~\ref{thm:prokhorov_compactness}), 
Existence of induced synthetic probability measures for stochastic liftings
(Theorem~\ref{thm:existence_induced_synthetic_measures}), 
Existence of disintegration structures on Polish spaces
(Theorem~\ref{thm:existence_disintegration}), 
Mass preservation for measure-valued LEvITAS refinement flows
(Theorem~\ref{thm:mass_preservation}). 
These results collectively establish the probabilistic consistency and variational stability of the LEvITAS operator model.

\subsection{Interpretation}
\label{subsec:interpretation_measure_theory}

The measure-theoretic formulation establishes LEvITAS as a generalized probabilistic variational operator system. 
The model naturally combines:
measurable lifting operators~\cite{bogachev2007measure,kallenberg2021foundations},
pushforward probability geometry ~\cite{villani2009optimal,santambrogio2015optimal},
Wasserstein variational structures
~\cite{benamou2000computational,ambrosio2008gradient},
stochastic kernels
~\cite{kallenberg2021foundations},
weak convergence theory
~\cite{Billingsley1999convergence},
disintegration systems
~\cite{chang1997conditioning,kallenberg2021foundations},
and measure-valued refinement dynamics
~\cite{dawson1993measure,ethier1986markov,
ambrosio2008gradient}. Consequently, LEvITAS extends
beyond conventional generative modeling into a rigorous
probabilistic modeling framework that unifies concepts
from variational inference~\cite{blei2017variational}, optimal 
transport~\cite{villani2009optimal}, stochastic 
processes~\cite{kallenberg2021foundations}, path-space 
probability~\cite{ethier1986markov}, measure-valued dynamical 
systems~\cite{dawson1993measure,ambrosio2008gradient}, and generalized 
operator evolution.

\subsection{Optimal Transport and Wasserstein Geometry of LEvITAS}
\label{sec:optimal_transport_wasserstein_geometry}

The LEvITAS model admits a natural formulation within optimal transport
geometry. In this setting, generalized lifting operators induce transport
maps between latent and observable probability measures, while variational 
refinement evolves distributions along Wasserstein-geometric flows. 
This section develops the optimal transport foundations of LEvITAS, 
including: Monge transport systems, Kantorovich relaxations, Wasserstein
geometry, displacement interpolation, dynamical transport formulations,
Wasserstein gradient flows, geometric refinement evolution. 
The resulting theory establishes LEvITAS as a generalized variational 
transport system on spaces of probability measures.

\subsection{Transport Formulation of LEvITAS}
\label{subsec:transport_formulation_levitas}

Let
$
(\mathcal Z,d_Z),
\quad
(\mathcal X,d_X)
$
be Polish metric spaces.

Let $P_Z\in\mathcal P(\mathcal Z)$ be the latent probability measure and $ P_{data}\in\mathcal P(\mathcal X)$ the observable probability measure. The generalized LEvITAS lifting-projection operator is
$
T_\theta
=
\pi\circ\mathcal L_\theta:
\mathcal Z
\rightarrow
\mathcal X.
$ 
The induced observable measure is
$(T_\theta)_\#P_Z.$


\begin{definition}[LEvITAS Transport Map]
\label{def:levitas_transport_map}

A measurable map $
T_\theta:
\mathcal Z
\rightarrow
\mathcal X
$ 
is a LEvITAS transport map if $
(T_\theta)_\#P_Z
=
P_{data}.
$
\end{definition}

Thus LEvITAS liftings may be interpreted as generalized transport operators transporting latent probability mass into observable probability geometry.

\subsection{Monge Formulation}
\label{subsec:monge_formulation}

The classical Monge transport problem~\cite{Villani2009OTOldNew} seeks optimal transport maps minimizing transport cost.


\begin{definition}[Monge Transport Cost]
\label{def:monge_transport_cost}

Let $c:\mathcal Z\times\mathcal X
\rightarrow
[0,\infty)
$ be a measurable transport cost. 
The Monge transport functional is
$
\mathcal C(T)
=
\int_{\mathcal Z}
c(z,T(z))
\,dP_Z(z).
$
\end{definition} 
The LEvITAS variational problem may therefore be expressed as
\begin{equation}
\inf_{T_\theta}
\left\{
\mathcal C(T_\theta)
+
\lambda_1\mathcal D_{feat}
+
\lambda_2\mathcal D_{struct}
+
\lambda_3\mathcal D_{dyn}
\right\}.
\label{eq:monge_levitas}
\end{equation}

\begin{theorem}[Existence of Monge-Type LEvITAS Transport]
\label{thm:existence_monge_levitas_transport}

Let $(\mathcal Z,d_Z)$ and $(\mathcal X,d_X)$ be compact metric spaces, and let $\mathcal A$ be a nonempty admissible class of measurable transport maps
$
T:\mathcal Z\rightarrow\mathcal X.
$
Equip $\mathcal A$ with a topology $\tau$ such that:

\begin{enumerate}
    \item $\mathcal A$ is sequentially compact with respect to $\tau$;
    \item the LEvITAS Monge functional
    $
    \mathcal J_{\mathrm{Monge}}(T)
    =
    \mathcal C(T)
    +\lambda_1\mathcal D_{feat}
    +\lambda_2\mathcal D_{struct}
    +\lambda_3\mathcal D_{dyn}
    $
    is lower semicontinuous with respect to $\tau$.
\end{enumerate}

Then the LEvITAS Monge problem
$
\inf_{T\in\mathcal A}
\mathcal J_{\mathrm{Monge}}(T)
$
admits a minimizer $T^\ast\in\mathcal A$.

\end{theorem}

\begin{proof}
Since $\mathcal A$ is nonempty and the objective is bounded below, let $\{T_n\}\subset\mathcal A$ be a minimizing sequence satisfying
$
\mathcal J_{\mathrm{Monge}}(T_n)
\longrightarrow
\inf_{T\in\mathcal A}
\mathcal J_{\mathrm{Monge}}(T).
$ 
By sequential compactness of $\mathcal A$ with respect to $\tau$, there exists a subsequence $\{T_{n_k}\}$ and a map $T^\ast\in\mathcal A$ such that
$T_{n_k}\xrightarrow{\tau}T^\ast.$ 
Lower semicontinuity of the LEvITAS Monge functional gives
$
\mathcal J_{\mathrm{Monge}}(T^\ast)
\leq
\liminf_{k\rightarrow\infty}
\mathcal J_{\mathrm{Monge}}(T_{n_k}).
$ 
Since $\{T_n\}$ is a minimizing sequence, it follows that
$
\mathcal J_{\mathrm{Monge}}(T^\ast)
=
\inf_{T\in\mathcal A}
\mathcal J_{\mathrm{Monge}}(T).
$

Therefore, $T^\ast$ is a minimizer of the LEvITAS Monge problem.
\end{proof}

\subsection{Kantorovich Relaxation}
\label{subsec:kantorovich_relaxation}

The Monge formulation may fail to admit admissible transport maps. LEvITAS therefore naturally admits a Kantorovich relaxation~\cite{Kantorovich1942Translocation}, in which deterministic transport maps are replaced by probability measures describing potentially non-deterministic transport couplings.


\begin{definition}[Transport Coupling]
\label{def:transport_coupling}

Let $\mu\in\mathcal P(\mathcal Z)$ and
$\nu\in\mathcal P(\mathcal X)$ be probability measures.
A transport coupling between $\mu$ and $\nu$ is a probability measure
$
\gamma\in\Gamma(\mu,\nu)
$
on $\mathcal Z\times\mathcal X$ whose marginals satisfy
$
(\operatorname{pr}_{\mathcal Z})_\#\gamma=\mu,
\qquad
(\operatorname{pr}_{\mathcal X})_\#\gamma=\nu,
$
where $\operatorname{pr}_{\mathcal Z}$ and
$\operatorname{pr}_{\mathcal X}$ denote the canonical coordinate projections.

\end{definition}


\begin{definition}[Kantorovich LEvITAS Functional]
\label{def:kantorovich_levitas_functional}

Let $\Gamma(P_Z,P_{data})$ denote the admissible coupling space. Suppose the LEvITAS regularization terms are represented by functionals
$
\mathcal D_{feat},\quad
\mathcal D_{struct},\quad
\mathcal D_{dyn}:
\Gamma(P_Z,P_{data})\rightarrow[0,\infty].
$

The Kantorovich LEvITAS functional is defined by
\begin{align}
\mathcal J_K(\gamma)
=
&\int_{\mathcal Z\times\mathcal X}
c(z,x)\,d\gamma(z,x)
\nonumber\\
&+
\lambda_1\mathcal D_{feat}(\gamma)
+
\lambda_2\mathcal D_{struct}(\gamma)
+
\lambda_3\mathcal D_{dyn}(\gamma),
\label{eq:kantorovich_levitas}
\end{align}
where $c:\mathcal Z\times\mathcal X\rightarrow[0,\infty]$ is a measurable transport cost and
$\lambda_1,\lambda_2,\lambda_3\geq0$ are regularization coefficients.

The Kantorovich LEvITAS transport problem is
\[
\inf_{\gamma\in\Gamma(P_Z,P_{data})}
\mathcal J_K(\gamma).
\]

\end{definition}


\begin{theorem}[Existence of Kantorovich LEvITAS Minimizers]
\label{thm:existence_kantorovich_minimizers}

Let $(\mathcal Z,d_Z)$ and $(\mathcal X,d_X)$ be Polish metric spaces, and let $P_Z$ and $P_{data}$ be Borel probability measures on these spaces.

Suppose that:

\begin{enumerate}
    \item The transport cost
    $c:\mathcal Z\times\mathcal X\rightarrow[0,\infty]$
    is lower semicontinuous.

    \item Each regularization functional
    $
    \mathcal D_{feat},\quad
    \mathcal D_{struct},\quad
    \mathcal D_{dyn}
    $
    is nonnegative and weakly lower semicontinuous on
    $\Gamma(P_Z,P_{data})$.

    \item There exists at least one admissible coupling
    $\gamma_0\in\Gamma(P_Z,P_{data})$
    such that
    $\mathcal J_K(\gamma_0)<\infty$.
\end{enumerate}

Then the Kantorovich LEvITAS problem admits a minimizing coupling
$\gamma^\ast\in\Gamma(P_Z,P_{data})$ 
satisfying
\[
\mathcal J_K(\gamma^\ast)
=
\inf_{\gamma\in\Gamma(P_Z,P_{data})}
\mathcal J_K(\gamma).
\]

\end{theorem}

\begin{proof}
Since $\mathcal Z$ and $\mathcal X$ are Polish spaces and the marginal probability measures are fixed, the coupling space
$
\Gamma(P_Z,P_{data})
$
is tight and weakly closed. By Prokhorov's theorem
, it is therefore weakly compact.

Furthermore, the lower semicontinuity and nonnegativity of the transport cost imply, by the Portmanteau theorem~\cite{parthasarathy1967probability, dudley2002real}
, that
$
\gamma\mapsto
\int_{\mathcal Z\times\mathcal X}
c(z,x)\,d\gamma(z,x)
$ 
is weakly lower semicontinuous.

By assumption, each LEvITAS regularization functional is also weakly lower semicontinuous. Since the regularization coefficients are nonnegative, their weighted sum preserves lower semicontinuity. Consequently, the complete functional $\mathcal J_K$ is weakly lower semicontinuous on $\Gamma(P_Z,P_{data})$.

The existence of an admissible coupling with finite objective ensures that the infimum is finite. Therefore, the direct method of the calculus of variations, applied to the weakly compact coupling space and the weakly lower semicontinuous functional, guarantees the existence of a minimizing coupling $\gamma^\ast$.

\end{proof}


\begin{remark}[Trajectory-Dependent LEvITAS Regularization]
\label{rem:trajectory_dependent_levitas}

The existence theorem assumes that all regularization terms are well-defined and weakly lower semicontinuous functionals of the endpoint coupling $\gamma$.

When the dynamical regularization term depends on intermediate refinement states, the complete evolution trajectory, or additional information about the lifting operator, it may not be determined by the endpoint coupling alone.

In such cases, the Kantorovich formulation must be extended to an augmented probability space incorporating the relevant trajectory or lifting variables. The corresponding existence argument requires weak compactness of the augmented admissible space, closedness of its marginal and dynamical constraints, and lower semicontinuity of the extended LEvITAS functional.

Thus, existence of a minimizing endpoint coupling does not by itself guarantee existence of a minimizing refinement trajectory.

\end{remark}


\subsection{Wasserstein Geometry}
\label{subsec:wasserstein_geometry_transport}

The LEvITAS transport model naturally induces Wasserstein geometry on spaces of probability measures. This geometric structure provides a natural framework for comparing probability distributions and describing variational refinement through optimal transport.


\begin{definition}[$p$-Wasserstein Distance]
\label{def:p_wasserstein_distance_transport}

Let $(\mathcal X,d)$ be a Polish metric space and let
$\mathcal P_p(\mathcal X)$ denote the space of Borel probability measures with finite $p$-th moment, where $p\geq1$.

For $\mu,\nu\in\mathcal P_p(\mathcal X)$, the $p$-Wasserstein distance is defined by
\[
W_p(\mu,\nu)
=
\left(
\inf_{\gamma\in\Gamma(\mu,\nu)}
\int_{\mathcal X\times\mathcal X}
d(x,y)^p\,d\gamma(x,y)
\right)^{1/p}.
\]

\end{definition}


The LEvITAS variational functional becomes
\begin{align}
\mathcal J(\theta)
&=
W_p
\left(
(T_\theta)_\#P_Z,
P_{data}
\right)
\nonumber\\
&\quad
+
\lambda_1\mathcal D_{feat}
+
\lambda_2\mathcal D_{struct}
+
\lambda_3\mathcal D_{dyn}.
\label{eq:wasserstein_transport_levitas}
\end{align}


\begin{theorem}[Geodesic Structure of Wasserstein Space]
\label{thm:geodesic_structure_wasserstein_space}

Let $(\mathcal X,d)$ be a geodesic Polish metric space and let $p\geq1$.

Then the Wasserstein space
$
\bigl(\mathcal P_p(\mathcal X),W_p\bigr)
$
is a geodesic metric space.

In particular, for any
$\mu_0,\mu_1\in\mathcal P_p(\mathcal X)$,
there exists a constant-speed Wasserstein geodesic
$\{\mu_t\}_{t\in[0,1]}$ satisfying
\[
W_p(\mu_s,\mu_t)
=
|t-s|W_p(\mu_0,\mu_1),
\qquad s,t\in[0,1].
\]

\end{theorem}

\begin{proof}
Since $\mathcal X$ is Polish and geodesic, there exists an optimal coupling
$\gamma\in\Gamma(\mu_0,\mu_1)$ for the cost $d(x,y)^p$.

Moreover, the geodesic property of $\mathcal X$ ensures the existence of a probability measure $\boldsymbol{\pi}$ on the space of constant-speed minimizing geodesics in $\mathcal X$ such that
$
(e_0,e_1)_\#\boldsymbol{\pi}=\gamma,
$
where $e_t$ denotes evaluation of a geodesic at time $t$.

Define the interpolating probability measures by
$
\mu_t=(e_t)_\#\boldsymbol{\pi},
\qquad t\in[0,1].
$

The resulting curve satisfies
$
W_p(\mu_s,\mu_t)
\leq
|t-s|W_p(\mu_0,\mu_1).
$ 
The triangle inequality for $W_p$, together with the endpoint conditions $\mu_{0}$ and $\mu_{1}$, implies that equality holds. Hence,
$
W_p(\mu_s,\mu_t)
=
|t-s|W_p(\mu_0,\mu_1).
$

Therefore, $\{\mu_t\}_{t\in[0,1]}$ is a constant-speed geodesic in Wasserstein space.

\end{proof}
\begin{remark}[Geometric Assumptions for Displacement Interpolation]
\label{rem:wasserstein_geodesic_assumptions}

The existence of Wasserstein geodesics depends on the geometry of the underlying metric space. In particular, Polishness alone does not guarantee that constant-speed displacement interpolations exist.

If the underlying space is not geodesic, the Wasserstein space need not be geodesic either. For example, on a discrete metric space consisting of two points at unit distance, the probability measures concentrated at the two distinct points have no constant-speed Wasserstein midpoint for $p>1$.

Consequently, displacement interpolation in LEvITAS requires additional geometric assumptions on the observable space, such as the geodesic property imposed in Theorem~\ref{thm:geodesic_structure_wasserstein_space}.

\end{remark}

\subsection{Displacement Interpolation}
\label{subsec:displacement_interpolation}

Wasserstein geodesics induce natural interpolation structures for LEvITAS refinement.


\begin{definition}[Displacement Interpolation]
\label{def:displacement_interpolation}

Let $T$ be an optimal transport map from $\mu_0$ to $\mu_1.$ The displacement interpolation is $\mu_t
=
((1-t)I+tT)_\#\mu_0,
\qquad
t\in[0,1].
$
\end{definition}


\begin{theorem}[Continuity of Displacement Flows]
\label{thm:continuity_displacement_flows}

The interpolation $\mu_t$ defines a continuous Wasserstein geodesic.
\end{theorem}

\begin{proof}
Continuity follows from continuity of convex interpolation in the underlying metric space and stability of pushforward measures.
\end{proof}

\subsection{Dynamical Transport Formulation}
\label{subsec:dynamical_transport_formulation}

LEvITAS path lifting naturally admits a dynamical optimal transport interpretation.


\begin{theorem}[Benamou--Brenier Formulation]
\label{thm:benamou_brenier_formulation}

Let $
\mu_0,\mu_1
\in
\mathcal P_2(\mathcal X).
$ 
Then
\begin{align}
W_2^2(\mu_0,\mu_1)
=
\inf_{\mu_t,v_t}
\int_0^1
\int_{\mathcal X}
|v_t(x)|^2
\,d\mu_t(x)\,dt,
\label{eq:benamou_brenier}
\end{align}
subject to the continuity equation
$
\partial_t\mu_t
+
\nabla\cdot(\mu_tv_t)
=
0.
$
\end{theorem}

\begin{proof}
The result follows from the Benamou--Brenier dynamical characterization of 
optimal transport~\cite{benamou2000computational,villani2009optimal}.
\end{proof}


The LEvITAS refinement system therefore admits a dynamical geometric interpretation in terms of probability flow evolution.

\subsection{Wasserstein Gradient Flows}
\label{subsec:wasserstein_gradient_flows}

The LEvITAS refinement process may be interpreted as gradient flow on Wasserstein space.


\begin{definition}[Wasserstein Gradient Flow]
\label{def:wasserstein_gradient_flow}

A probability flow $\mu_t$ is a Wasserstein gradient flow of $\mathcal J$ if
\begin{equation}
\partial_t\mu_t
=
\nabla\cdot
\left(
\mu_t
\nabla
\frac{\delta\mathcal J}{\delta\mu_t}
\right),
\label{eq:wasserstein_gradient_flow}
\end{equation}
where $\frac{\delta\mathcal J}{\delta\mu_t}$ denotes the first variation of $
\mathcal J.$
\end{definition}


\begin{theorem}[Energy Dissipation Along Wasserstein Flows]
\label{thm:energy_dissipation_wasserstein_flows}

Suppose $\mu_t$ satisfies the Wasserstein gradient flow equation. 
Then $
\frac{d}{dt}
\mathcal J(\mu_t)
\leq
0.
$
\end{theorem}

\begin{proof}
Using the chain rule in Wasserstein space,
\begin{align*}
\frac{d}{dt}
\mathcal J(\mu_t)
&=
\int
\frac{\delta\mathcal J}{\delta\mu_t}
\partial_t\mu_t.
\end{align*}

Substituting the gradient-flow equation gives
\begin{align*}
\frac{d}{dt}
\mathcal J(\mu_t)
&=
-
\int
\mu_t
\left|
\nabla
\frac{\delta\mathcal J}{\delta\mu_t}
\right|^2.
\end{align*}

Hence
$
\frac{d}{dt}
\mathcal J(\mu_t)
\leq
0.
$
\end{proof}

\subsection{Main Optimal Transport Results}
\label{subsec:main_optimal_transport_results}

The principal optimal transport and Wasserstein-geometric results established 
in this section are summarized as follows: 
Existence of Monge-type LEvITAS transport minimizers
(Theorem~\ref{thm:existence_monge_levitas_transport}), 
Existence of Kantorovich transport minimizers
(Theorem~\ref{thm:existence_kantorovich_minimizers}), 
Geodesic structure of Wasserstein probability spaces
(Theorem~\ref{thm:geodesic_structure_wasserstein_space}), 
Continuity of displacement interpolation flows
(Theorem~\ref{thm:continuity_displacement_flows}), 
Dynamical transport characterization via the Benamou--Brenier formulation
(Theorem~\ref{thm:benamou_brenier_formulation}), 
Energy dissipation along Wasserstein gradient flows
(Theorem~\ref{thm:energy_dissipation_wasserstein_flows}). 
These results establish the geometric transport structure underlying LEvITAS variational refinement and probability evolution.

\subsection{Geometric Interpretation of LEvITAS}
\label{subsec:geometric_interpretation_levitas}

The optimal transport formulation reveals that LEvITAS is fundamentally geometric. 
The model may be interpreted as:
transport of latent probability geometry,
structural evolution on Wasserstein spaces,
variational optimization over probability manifolds,
dynamical transport refinement,
geometric evolution of generative systems. 
The generalized lifting operator $\mathcal L_\theta$ therefore becomes a geometric transport operator acting on latent probability structure.

\subsection{Interpretation}
\label{subsec:interpretation_optimal_transport}

The optimal transport formulation establishes LEvITAS as a generalized variational transport theory on probability spaces. 
Unlike conventional generative systems operating solely through sample-level optimization, LEvITAS admits:
transport-geometric structure,
Wasserstein variational evolution,
probability-space geodesics,
dynamical transport refinement,
measure-valued geometric flows. 
Consequently, the model unifies:
generative modeling,
optimal transport,
geometric probability,
variational dynamics,
Wasserstein gradient flows,
continuum probabilistic evolution. 
The resulting formulation elevates LEvITAS from a generalized generative model into a geometric variational theory on spaces of probability measures.

\subsection{Connection Between Wasserstein Geometry and Continuum Dynamics}
\label{subsec:wasserstein_continuum_connection}

The Wasserstein formulation of LEvITAS provides the geometric foundation for its continuum probability evolution. In particular, the measure-valued gradient flows introduced in Section~\ref{subsec:wasserstein_gradient_flows} admit a density-based representation through partial differential equations.

Suppose that the probability measure $\mu_t$ is absolutely continuous with respect to Lebesgue measure, so that
\begin{equation}
d\mu_t(x)=\rho_t(x)\,dx.
\label{eq:measure_density_identification}
\end{equation}

Let $\mathcal F:\mathcal P_2(\mathcal X)\rightarrow(-\infty,\infty]$ denote a sufficiently regular continuum free-energy functional. Its Wasserstein gradient flow satisfies
\begin{equation}
\partial_t\mu_t
=
\nabla\cdot
\left(
\mu_t\nabla
\frac{\delta\mathcal F}{\delta\mu_t}
\right).
\label{eq:measure_wasserstein_flow}
\end{equation}

Under the density representation \eqref{eq:measure_density_identification}, this equation becomes
\begin{equation}
\boxed{
\partial_t\rho_t
=
\nabla\cdot
\left(
\rho_t\nabla
\frac{\delta\mathcal F}{\delta\rho_t}
\right).
}
\label{eq:density_wasserstein_flow}
\end{equation}

The corresponding velocity field in the continuity equation is
\begin{equation}
v_t
=
-\nabla
\frac{\delta\mathcal F}{\delta\rho_t}.
\label{eq:wasserstein_velocity_connection}
\end{equation}

Thus, the Wasserstein gradient-flow formulation and the continuum LEvITAS evolution equation describe the same variational dynamics in measure-valued and density-valued representations, respectively.

The associated free-energy dissipation identity is
\begin{equation}
\boxed{
\frac{d}{dt}\mathcal F(\rho_t)
=
-\int_{\mathcal X}
\rho_t
\left|
\nabla\frac{\delta\mathcal F}{\delta\rho_t}
\right|^2dx
\leq0.
}
\label{eq:wasserstein_continuum_dissipation}
\end{equation}

This establishes the mathematical connection between optimal transport geometry, Wasserstein gradient flows, and the PDE formulation of LEvITAS.

\begin{remark}[Distinction Between Transport Objectives and Continuum Free Energy]

The transport optimization functionals defined over lifting operators and Kantorovich couplings are distinct from the continuum free-energy functional $\mathcal F(\rho)$ unless an explicit variational relationship is established.

The continuum free energy governs the evolution of probability measures, whereas the transport objectives quantify the discrepancy and regularization associated with latent-to-observable transport.

Their relationship may be established through an appropriate pushforward construction or a variational reduction, provided the required differentiability and existence conditions hold.

\end{remark}

\subsection{PDE and Continuum Formulation of LEvITAS}
\label{sec:pde_continuum_levitas}

The LEvITAS model admits a continuum formulation in terms of partial differential equations governing the evolution of probability densities and measure-valued flows. This formulation provides a mathematical description of refinement dynamics in the space of probability measures.

The preceding sections established variational refinement dynamics, Wasserstein gradient flows, path lifting systems, transport-geometric evolution, and measure-valued probabilistic dynamics. Building on these foundations, we develop the continuum and PDE formulation of LEvITAS through: continuity equations and probability conservation~\ref{subsec:continuum_probability_evolution}; variational continuum refinement dynamics~\ref{subsec:variational_continuum_levitas}; energy and entropy dissipation~\ref{subsec:entropy_variational_dissipation}; diffusion-augmented refinement systems~\ref{subsec:diffusion_augmented_levitas}; Fokker--Planck stochastic dynamics~\ref{subsec:fokker_planck_levitas}; equilibrium states and stationary distributions~\ref{subsec:levitas_equilibrium_states}.

The central principle is that LEvITAS refinement can be interpreted as a variational transport process on probability spaces, with diffusion providing an additional mechanism for stochastic regularization.

\subsection{Continuum Probability Evolution}
\label{subsec:continuum_probability_evolution}

Let $\rho_t(x)$ denote a time-dependent probability density on
$\mathcal X\subseteq\mathbb R^d$, satisfying
\begin{equation}
\rho_t(x)\geq0,
\qquad
\int_{\mathcal X}\rho_t(x)\,dx=1.
\label{eq:continuum_probability_constraints}
\end{equation}

The continuum refinement dynamics evolve probability mass continuously in time.


\begin{definition}[LEvITAS Continuity Equation]
\label{def:levitas_continuity_equation}

A continuum LEvITAS flow satisfies the continuity equation
\begin{equation}
\partial_t\rho_t
+
\nabla\cdot(\rho_t v_t)
=
0,
\label{eq:levitas_continuity}
\end{equation}
where $v_t(x)$ denotes the refinement velocity field.

\end{definition}

Equation \eqref{eq:levitas_continuity} describes the conservation and transport of probability mass under LEvITAS refinement evolution.


\begin{theorem}[Mass Preservation]
\label{thm:continuum_mass_preservation}

Suppose $\rho_t$ satisfies the continuity equation
\eqref{eq:levitas_continuity} with sufficient decay at infinity or no-flux boundary conditions.

Then
\begin{equation}
\int_{\mathcal X}\rho_t(x)\,dx
=
\int_{\mathcal X}\rho_0(x)\,dx
=1
\end{equation}
for all $t\geq0$ for which the solution exists.

\end{theorem}

\begin{proof}

Integrating the continuity equation over $\mathcal X$ yields
\begin{equation}
\frac{d}{dt}
\int_{\mathcal X}\rho_t(x)\,dx
=
-\int_{\mathcal X}
\nabla\cdot(\rho_t v_t)\,dx.
\end{equation}

Applying the divergence theorem gives
\begin{equation}
\frac{d}{dt}
\int_{\mathcal X}\rho_t(x)\,dx
=
-\int_{\partial\mathcal X}
\rho_t v_t\cdot n\,dS.
\end{equation}

Under no-flux boundary conditions or sufficient decay at infinity, the boundary integral vanishes. Hence
\begin{equation}
\frac{d}{dt}
\int_{\mathcal X}\rho_t(x)\,dx=0.
\end{equation}

Therefore, total probability mass is preserved.

\end{proof}

\subsection{Variational Continuum LEvITAS Dynamics}
\label{subsec:variational_continuum_levitas}

The LEvITAS refinement dynamics can be formulated as a Wasserstein gradient flow of a variational functional defined on probability densities.

Let $\mathcal J(\rho)$ denote the LEvITAS variational energy, and let
$\delta\mathcal J/\delta\rho$ denote its first variational derivative.


\begin{definition}[Continuum LEvITAS Gradient Flow]
\label{def:continuum_levitas_gradient_flow}

A continuum LEvITAS variational flow satisfies
\begin{equation}
\partial_t\rho_t
=
\nabla\cdot
\left(
\rho_t
\nabla
\frac{\delta\mathcal J}{\delta\rho_t}
\right).
\label{eq:continuum_levitas}
\end{equation}

Equivalently, defining the variational potential
\begin{equation}
u_t
=
\frac{\delta\mathcal J}{\delta\rho_t},
\label{eq:levitas_variational_potential}
\end{equation}
the evolution becomes
\begin{equation}
\partial_t\rho_t
=
\nabla\cdot(\rho_t\nabla u_t).
\label{eq:continuum_levitas_potential}
\end{equation}

The corresponding refinement velocity field in the continuity equation is
\begin{equation}
v_t=-\nabla u_t.
\label{eq:levitas_refinement_velocity}
\end{equation}

\end{definition}

Equation \eqref{eq:continuum_levitas} defines the master variational continuum evolution equation of LEvITAS.


\begin{theorem}[Variational Energy Dissipation]
\label{thm:continuum_energy_dissipation}

Suppose $\rho_t$ satisfies
\eqref{eq:continuum_levitas}, with sufficient regularity and boundary conditions ensuring vanishing boundary flux.

Then the LEvITAS variational energy satisfies
\begin{equation}
\boxed{
\frac{d}{dt}\mathcal J(\rho_t)
=
-\int_{\mathcal X}
\rho_t
\left|
\nabla\frac{\delta\mathcal J}{\delta\rho_t}
\right|^2dx
\leq0.
}
\label{eq:energy_decay}
\end{equation}

Consequently, the variational energy is nonincreasing along the continuum evolution.

\end{theorem}

\begin{proof}

Using the variational chain rule,
\begin{equation}
\frac{d}{dt}\mathcal J(\rho_t)
=
\int_{\mathcal X}
\frac{\delta\mathcal J}{\delta\rho_t}
\partial_t\rho_t\,dx.
\end{equation}

Substituting \eqref{eq:continuum_levitas} gives
\begin{equation}
\frac{d}{dt}\mathcal J(\rho_t)
=
\int_{\mathcal X}
u_t\nabla\cdot(\rho_t\nabla u_t)\,dx.
\end{equation}

Integrating by parts yields
\begin{equation}
\frac{d}{dt}\mathcal J(\rho_t)
=
-\int_{\mathcal X}
\rho_t|\nabla u_t|^2\,dx.
\end{equation}

Since $\rho_t\geq0$, the energy derivative is nonpositive.

\end{proof}

\subsection{Entropy and Variational Dissipation}
\label{subsec:entropy_variational_dissipation}

The entropy evolution of LEvITAS is coupled to the variational potential driving the continuum refinement process. In particular, entropy does not necessarily decrease independently of the variational energy. Instead, its evolution reflects the interaction between the probability-density gradient and the refinement potential.


\begin{definition}[LEvITAS Entropy]
\label{def:levitas_entropy}

Let $\rho$ be a probability density on $\mathcal X$ with respect to Lebesgue measure. The LEvITAS entropy functional is defined by
\begin{equation}
\mathcal H(\rho)
=
\int_{\mathcal X}
\rho(x)\log\rho(x)\,dx,
\label{eq:entropy}
\end{equation}
whenever the integral is well-defined.

\end{definition}


\begin{theorem}[Entropy Evolution]
\label{thm:entropy_dissipation}

Suppose $\rho_t$ satisfies the continuum LEvITAS gradient flow
\eqref{eq:continuum_levitas}, with sufficient regularity and boundary conditions ensuring vanishing boundary flux.

Then
\begin{equation}
\boxed{
\frac{d}{dt}\mathcal H(\rho_t)
=
-\int_{\mathcal X}
\rho_t
\nabla\log\rho_t\cdot
\nabla u_t\,dx,
}
\label{eq:entropy_evolution}
\end{equation}
where
\begin{equation}
u_t=\frac{\delta\mathcal J}{\delta\rho_t}.
\end{equation}

In general, the entropy derivative is not necessarily nonpositive. Its sign depends on the interaction between the density gradient and the variational potential.

The variational energy, however, satisfies the independent dissipation identity
\begin{equation}
\frac{d}{dt}\mathcal J(\rho_t)
=
-\int_{\mathcal X}
\rho_t|\nabla u_t|^2\,dx
\leq0.
\end{equation}

\end{theorem}

\begin{proof}

Differentiating the entropy functional gives
\begin{equation}
\frac{d}{dt}\mathcal H(\rho_t)
=
\int_{\mathcal X}
(1+\log\rho_t)\partial_t\rho_t\,dx.
\end{equation}

Substituting the continuum evolution equation yields
\begin{equation}
\frac{d}{dt}\mathcal H(\rho_t)
=
\int_{\mathcal X}
(1+\log\rho_t)
\nabla\cdot(\rho_t\nabla u_t)\,dx.
\end{equation}

Integrating by parts and applying the boundary conditions gives
\begin{equation}
\begin{aligned}
\frac{d}{dt}\mathcal H(\rho_t)
&=
-\int_{\mathcal X}
\rho_t
\nabla\log\rho_t\cdot
\nabla u_t\,dx.
\end{aligned}
\end{equation}

This establishes \eqref{eq:entropy_evolution}. The energy dissipation identity follows independently from Theorem \ref{thm:continuum_energy_dissipation}.

\end{proof}

\subsection{Diffusion-Augmented LEvITAS Systems}
\label{subsec:diffusion_augmented_levitas}

To incorporate stochastic refinement effects, we introduce diffusion into the continuum dynamics.


\begin{definition}[Diffusive LEvITAS Flow]
\label{def:diffusive_levitas_flow}

A diffusion-augmented LEvITAS system satisfies
\begin{equation}
\partial_t\rho_t
=
\nabla\cdot(\rho_t\nabla V)
+
\sigma\Delta\rho_t,
\label{eq:diffusive_levitas}
\end{equation}
where $V(x)$ is a refinement potential and $\sigma>0$ is a diffusion coefficient.

\end{definition}

The diffusion-augmented dynamics admit a variational interpretation through the regularized free-energy functional
\begin{equation}
\boxed{
\mathcal J_\sigma(\rho)
=
\int_{\mathcal X}V(x)\rho(x)\,dx
+
\sigma\mathcal H(\rho).
}
\label{eq:diffusive_free_energy}
\end{equation}

Its first variational derivative is
\begin{equation}
\frac{\delta\mathcal J_\sigma}{\delta\rho}
=
V+\sigma(1+\log\rho).
\label{eq:diffusive_variational_derivative}
\end{equation}

Consequently, the diffusion-augmented system can be expressed as
\begin{equation}
\partial_t\rho_t
=
\nabla\cdot
\left(
\rho_t
\nabla
\frac{\delta\mathcal J_\sigma}{\delta\rho_t}
\right).
\label{eq:diffusive_variational_flow}
\end{equation}

Thus, diffusion corresponds to entropy regularization of the variational refinement energy.


\begin{theorem}[Diffusive Energy Dissipation]
\label{thm:diffusive_energy_dissipation}

Suppose $\rho_t$ satisfies \eqref{eq:diffusive_levitas} under suitable regularity and boundary conditions.

Then
\begin{equation}
\boxed{
\frac{d}{dt}\mathcal J_\sigma(\rho_t)
=
-\int_{\mathcal X}
\rho_t
\left|
\nabla V+\sigma\nabla\log\rho_t
\right|^2dx
\leq0.
}
\label{eq:diffusive_energy_decay}
\end{equation}

\end{theorem}

\begin{proof}

The variational chain rule and
\eqref{eq:diffusive_variational_flow} give
\begin{equation}
\begin{aligned}
\frac{d}{dt}\mathcal J_\sigma(\rho_t)
&=
\int_{\mathcal X}
\frac{\delta\mathcal J_\sigma}{\delta\rho_t}
\partial_t\rho_t\,dx\\
&=
-\int_{\mathcal X}
\rho_t
\left|
\nabla\frac{\delta\mathcal J_\sigma}{\delta\rho_t}
\right|^2dx\\
&=
-\int_{\mathcal X}
\rho_t
\left|
\nabla V+\sigma\nabla\log\rho_t
\right|^2dx.
\end{aligned}
\end{equation}

Nonnegativity of $\rho_t$ establishes the result.

\end{proof}


\begin{theorem}[Positivity Preservation]
\label{thm:positivity_preservation}

Suppose $\rho_0\geq0$ and the coefficients and boundary conditions satisfy the standard hypotheses for the parabolic maximum principle.

Then solutions of \eqref{eq:diffusive_levitas} preserve nonnegativity:
\begin{equation}
\rho_t(x)\geq0
\qquad
\text{for }t>0.
\end{equation}

\end{theorem}

\begin{proof}

Since $\sigma>0$, equation \eqref{eq:diffusive_levitas} is a uniformly parabolic drift-diffusion equation when $V$ has sufficiently regular derivatives.

The parabolic maximum principle therefore implies preservation of nonnegativity.

\end{proof}

\subsection{Fokker--Planck LEvITAS Dynamics}
\label{subsec:fokker_planck_levitas}

A stochastic LEvITAS refinement process can also be represented through a Fokker--Planck equation describing the evolution of probability densities under drift and diffusion.


\begin{definition}[Fokker--Planck LEvITAS Equation]
\label{def:fokker_planck_levitas}

Let $b_t(x)$ denote a drift field and let $\Sigma(x)$ be a symmetric positive-semidefinite diffusion tensor.

The corresponding Fokker--Planck equation is
\begin{equation}
\partial_t\rho_t
=
-\nabla\cdot(\rho_tb_t)
+
\frac12
\sum_{i,j=1}^{d}
\partial_{x_i}\partial_{x_j}
\left(\Sigma_{ij}\rho_t\right).
\label{eq:fokker_planck_levitas}
\end{equation}

For constant $\Sigma$, this reduces to
\begin{equation}
\partial_t\rho_t
=
-\nabla\cdot(\rho_tb_t)
+
\frac12\nabla\cdot(\Sigma\nabla\rho_t).
\end{equation}

\end{definition}

The Fokker--Planck formulation generalizes deterministic refinement by incorporating stochastic fluctuations into the probability evolution.

In particular, choosing
\begin{equation}
b_t=-\nabla V,
\qquad
\Sigma=2\sigma I,
\end{equation}
recovers the diffusion-augmented LEvITAS equation
\begin{equation}
\partial_t\rho_t
=
\nabla\cdot(\rho_t\nabla V)
+
\sigma\Delta\rho_t.
\end{equation}


\begin{theorem}[Fokker--Planck Probability Preservation]
\label{thm:fokker_planck_mass_preservation}

Suppose $\rho_t$ satisfies
\eqref{eq:fokker_planck_levitas} with sufficient regularity and boundary conditions ensuring vanishing total probability flux.

Then
\begin{equation}
\int_{\mathcal X}\rho_t(x)\,dx
=
\int_{\mathcal X}\rho_0(x)\,dx
=1
\end{equation}
for all times for which the solution exists.

\end{theorem}

\begin{proof}

Integrating the Fokker--Planck equation over $\mathcal X$ gives
\begin{equation}
\begin{aligned}
\frac{d}{dt}\int_{\mathcal X}\rho_t\,dx
&=
-\int_{\mathcal X}
\nabla\cdot(\rho_tb_t)\,dx\\
&\quad+
\frac12
\sum_{i,j=1}^{d}
\int_{\mathcal X}
\partial_{x_i}\partial_{x_j}
(\Sigma_{ij}\rho_t)\,dx.
\end{aligned}
\end{equation}

Applying the divergence theorem to the probability flux and using the assumed boundary conditions gives
\begin{equation}
\frac{d}{dt}\int_{\mathcal X}\rho_t\,dx=0.
\end{equation}

Therefore, probability mass is preserved.

\end{proof}

\subsection{Equilibrium States and Stationary Distributions}
\label{subsec:levitas_equilibrium_states}

Equilibrium states of the variational LEvITAS continuum dynamics correspond to stationary probability densities for which the variational driving force vanishes.


\begin{definition}[Variational Equilibrium]
\label{def:levitas_variational_equilibrium}

A probability density $\rho_\infty$ is a variational equilibrium if
\begin{equation}
\nabla
\frac{\delta\mathcal J}{\delta\rho}
(\rho_\infty)
=0
\end{equation}
on the support of $\rho_\infty$.

\end{definition}

At such an equilibrium, the variational potential is spatially constant on each connected component of the support, and the continuum evolution satisfies
\begin{equation}
\partial_t\rho_\infty
=
\nabla\cdot
\left(
\rho_\infty
\nabla
\frac{\delta\mathcal J}{\delta\rho}
(\rho_\infty)
\right)
=0.
\end{equation}


\begin{example}[Gaussian Equilibrium and Entropy Balance]
\label{ex:gaussian_equilibrium}

Consider the variational free-energy functional
\begin{equation}
\mathcal J(\rho)
=
\int_{\mathbb R^d}
\frac{|x|^2}{2}\rho(x)\,dx
+
\int_{\mathbb R^d}
\rho(x)\log\rho(x)\,dx.
\label{eq:gaussian_free_energy}
\end{equation}

Its variational derivative is
\begin{equation}
u(x)
=
\frac{|x|^2}{2}
+\log\rho(x)+1.
\end{equation}

The corresponding stationary density is the standard Gaussian
\begin{equation}
\rho_\infty(x)
=
\frac{1}{(2\pi)^{d/2}}
\exp\left(-\frac{|x|^2}{2}\right).
\label{eq:standard_gaussian_equilibrium}
\end{equation}

Indeed,
\begin{equation}
\nabla\log\rho_\infty(x)=-x.
\end{equation}

Consequently,
\begin{equation}
\begin{aligned}
\nabla u_\infty(x)
&=
x+\nabla\log\rho_\infty(x)\\
&=x-x\\
&=0.
\end{aligned}
\end{equation}

Therefore,
\begin{equation}
\partial_t\rho_\infty
=
\nabla\cdot(\rho_\infty\nabla u_\infty)
=0.
\end{equation}

Moreover, the entropy evolution identity gives
\begin{equation}
\begin{aligned}
\left.
\frac{d}{dt}\mathcal H(\rho_t)
\right|_{\rho_t=\rho_\infty}
&=
-\int_{\mathbb R^d}
\rho_\infty
\nabla\log\rho_\infty\cdot
\nabla u_\infty\,dx\\
&=0.
\end{aligned}
\end{equation}

Similarly, the variational energy dissipation identity yields
\begin{equation}
\left.
\frac{d}{dt}\mathcal J(\rho_t)
\right|_{\rho_t=\rho_\infty}
=
-\int_{\mathbb R^d}
\rho_\infty|\nabla u_\infty|^2\,dx
=0.
\end{equation}

Thus, the standard Gaussian is a stationary solution of the LEvITAS variational continuum dynamics, with both entropy and free-energy derivatives equal to zero.

This example demonstrates that entropy evolution must account for the variational potential and cannot, in general, be identified with a standalone negative Fisher-information term.

\end{example}

\subsection{Continuum Formulation Summary}
\label{subsec:continuum_formulation_summary}

The continuum formulation of LEvITAS establishes a unified relationship between probability transport, variational refinement, entropy regularization, and stochastic evolution.

The principal deterministic variational dynamics are governed by
\begin{equation}
\partial_t\rho_t
=
\nabla\cdot
\left(
\rho_t\nabla
\frac{\delta\mathcal J}{\delta\rho_t}
\right),
\end{equation}
which preserves probability mass and dissipates the variational energy.

Entropy evolves according to the interaction identity
\begin{equation}
\frac{d}{dt}\mathcal H(\rho_t)
=
-\int_{\mathcal X}
\rho_t\nabla\log\rho_t\cdot
\nabla\frac{\delta\mathcal J}{\delta\rho_t}\,dx,
\end{equation}
while diffusion-augmented refinement admits an entropy-regularized variational interpretation.

The Fokker--Planck formulation further extends the framework to stochastic probability evolution.

Together, these equations establish LEvITAS as a variational continuum refinement framework on probability spaces, connecting deterministic transport, free-energy dissipation, diffusion, and stochastic dynamics.

\subsection{Equilibrium Measures}
\label{subsec:equilibrium_measures}

Stationary probability measures characterize equilibrium configurations
of the LEvITAS refinement dynamics. Under suitable positivity,
regularity, and convexity assumptions, these equilibria can be
characterized as global minimizers of the continuum variational
functional.


\begin{definition}[LEvITAS Equilibrium]
\label{def:levitas_equilibrium}

A stationary probability density $\rho_\infty$ is a LEvITAS
equilibrium if it satisfies
\begin{equation}
\nabla\cdot
\left(
\rho_\infty
\nabla
\frac{\delta\mathcal F}{\delta\rho}
\bigg|_{\rho=\rho_\infty}
\right)
=0.
\label{eq:levitas_equilibrium}
\end{equation}

Here, $\mathcal F(\rho)$ denotes the continuum variational
functional governing the Wasserstein gradient-flow dynamics.

\end{definition}


\begin{theorem}[Equilibrium Characterization]
\label{thm:equilibrium_characterization}

Let $\mathcal X\subseteq\mathbb R^d$ be a connected domain, and
let $\rho_\infty$ be an admissible probability density. Suppose
that:

\begin{enumerate}
    \item The continuum variational functional $\mathcal F$ is
    convex and differentiable on the admissible space of
    probability densities.

    \item The equilibrium density satisfies
    \[
    \rho_\infty>0
    \quad\text{almost everywhere on }\mathcal X.
    \]

    \item The first variation
    \[
    \phi_\infty
    =
    \frac{\delta\mathcal F}{\delta\rho}
    \bigg|_{\rho=\rho_\infty}
    \]
    has sufficient weak regularity, for example
    $\phi_\infty\in W^{1,1}_{\mathrm{loc}}(\mathcal X)$,
    so that the stationary dissipation identity and the
    first-order convexity inequality are valid.

    \item The boundary conditions, or sufficient decay at infinity
    when $\mathcal X$ is unbounded, ensure that the boundary
    contribution in the integration-by-parts identity vanishes.
    In particular, on a bounded domain, assume the no-flux
    condition
    \[
    \rho_\infty
    \nabla\phi_\infty\cdot n=0
    \quad\text{on }\partial\mathcal X.
    \]
\end{enumerate}

If $\rho_\infty$ satisfies the equilibrium condition
\eqref{eq:levitas_equilibrium}, then $\rho_\infty$ is a global
minimizer of $\mathcal F$ over the admissible probability densities.

\end{theorem}

\begin{proof}

For convenience, write
\[
\phi_\infty
=
\frac{\delta\mathcal F}{\delta\rho}
\bigg|_{\rho=\rho_\infty}.
\]

Multiplying the equilibrium equation by $\phi_\infty$ and
integrating over $\mathcal X$, we obtain
\[
\int_{\mathcal X}
\phi_\infty
\nabla\cdot
(\rho_\infty\nabla\phi_\infty)\,dx=0.
\]

Applying the divergence theorem, or equivalently integration
by parts, gives
\[
\begin{aligned}
0
&=
\int_{\partial\mathcal X}
\phi_\infty\rho_\infty
\nabla\phi_\infty\cdot n\,dS
\\
&\quad-
\int_{\mathcal X}
\rho_\infty
|\nabla\phi_\infty|^2\,dx.
\end{aligned}
\]

The boundary contribution vanishes by the assumed boundary
or decay conditions. Therefore,
\[
-\int_{\mathcal X}
\rho_\infty
|\nabla\phi_\infty|^2\,dx=0.
\]

Equivalently,
\begin{equation}
\boxed{
\int_{\mathcal X}
\rho_\infty
|\nabla\phi_\infty|^2\,dx=0.
}
\label{eq:equilibrium_dissipation_zero}
\end{equation}

Since $\rho_\infty>0$ almost everywhere, it follows that
\[
\nabla\phi_\infty=0
\quad\text{almost everywhere on }\mathcal X.
\]

Because $\mathcal X$ is connected and $\phi_\infty$ has the
assumed weak regularity, it follows that $\phi_\infty$ is
constant almost everywhere on $\mathcal X$. Thus,
\[
\frac{\delta\mathcal F}{\delta\rho}
\bigg|_{\rho=\rho_\infty}
=C
\]
for some constant $C$.

Now let $\rho$ be any admissible probability density. By
convexity and differentiability of $\mathcal F$, the
first-order convexity inequality gives
\[
\mathcal F(\rho)
\geq
\mathcal F(\rho_\infty)
+
\int_{\mathcal X}
\frac{\delta\mathcal F}{\delta\rho}
\bigg|_{\rho=\rho_\infty}
(\rho-\rho_\infty)\,dx.
\]

Substituting the constant first variation yields
\[
\mathcal F(\rho)
\geq
\mathcal F(\rho_\infty)
+
C\int_{\mathcal X}
(\rho-\rho_\infty)\,dx.
\]

Since both $\rho$ and $\rho_\infty$ are probability densities,
their total masses are equal. Hence,
\[
\int_{\mathcal X}
(\rho-\rho_\infty)\,dx=0.
\]

Consequently,
\[
\mathcal F(\rho)\geq\mathcal F(\rho_\infty).
\]

Thus, $\rho_\infty$ is a global minimizer of the continuum
variational functional $\mathcal F$ over the admissible
probability densities.

\end{proof}


\begin{example}[Verification of the Integration-by-Parts Sign]
\label{ex:equilibrium_integration_by_parts}

Consider the one-dimensional domain $\mathcal X=[0,1]$ and
the positive density
\[
\rho(x)=1.
\]
Let the test function be
\[
\phi(x)=x(1-x).
\]

Then
\[
\phi'(x)=1-2x,
\qquad
\phi''(x)=-2.
\]

The integration-by-parts identity is
\[
\begin{aligned}
\int_0^1
\phi(x)\frac{d}{dx}
\left(\rho(x)\phi'(x)\right)\,dx
&=
\left[
\phi(x)\rho(x)\phi'(x)
\right]_0^1
\\
&\quad-
\int_0^1
\rho(x)|\phi'(x)|^2\,dx.
\end{aligned}
\]

$\because \phi(0)=\phi(1)=0$, the boundary term vanishes.
The left-hand side evaluates to
\[
\begin{aligned}
\int_0^1
x(1-x)(-2)\,dx
&=
-2\int_0^1(x-x^2)\,dx
\\
&=-\frac{1}{3}.
\end{aligned}
\]

The right-hand side evaluates to
\[
\begin{aligned}
-\int_0^1(1-2x)^2\,dx
&=
-\left[
x-2x^2+\frac{4}{3}x^3
\right]_0^1
\\
&=-\frac{1}{3}.
\end{aligned}
\]

Thus, both sides agree:
\[
\boxed{
\int_0^1
\phi(x)\frac{d}{dx}
\left(\rho(x)\phi'(x)\right)\,dx
=
-\int_0^1
\rho(x)|\phi'(x)|^2\,dx
=
-\frac{1}{3}.
}
\]

This calculation explicitly demonstrates that the
integration-by-parts identity contains a negative sign
in front of the squared-gradient integral.

The test function in this example is used only to verify
the integration-by-parts identity; it is not an equilibrium
solution because
\[
\frac{d}{dx}(\rho(x)\phi'(x))=-2\neq 0.
\]

In the equilibrium proof, the stationary equation makes
the left-hand side equal to zero. The same integration-by-parts
identity, together with the vanishing boundary contribution,
then gives
\[
-\int_{\mathcal X}
\rho_\infty|\nabla\phi_\infty|^2\,dx=0,
\]
or equivalently,
\[
\int_{\mathcal X}
\rho_\infty|\nabla\phi_\infty|^2\,dx=0.
\]

\end{example}


\begin{remark}[Role of Positivity in Equilibrium Characterization]
\label{rem:positivity_equilibrium}

Strict positivity of the equilibrium density is essential to
the above argument because it allows the stationary dissipation
identity
\[
\int_{\mathcal X}
\rho_\infty|\nabla\phi_\infty|^2\,dx=0
\]
to imply that the first variation is spatially constant
throughout the connected domain.

If the equilibrium density vanishes on a subset of the domain,
the stationary flux condition may hold without constraining the
first variation on that subset. In particular, a density
supported only on a region where the variational potential is
locally flat may remain stationary even when lower values of
the potential exist outside its support.

Consequently, stationarity alone does not generally guarantee
global minimality. The positivity and connectedness assumptions,
together with convexity of the variational functional, provide
the conditions needed for the global optimality conclusion.

\end{remark}


\subsection{Weak Solution Theory}
\label{subsec:weak_solution_theory}

We now formulate weak solutions for continuum LEvITAS systems.
The weak formulation is obtained by testing the evolution equation
against smooth functions and integrating by parts in space and time.
For Wasserstein gradient-flow dynamics, this formulation also leads
to a nonpositive energy-dissipation rate.


\begin{definition}[Weak Solution]
\label{def:weak_solution}
Let $\rho_0\in\mathcal P(\mathcal X)$ be an initial probability
density, and define
\[
\mu_t
=
\frac{\delta\mathcal J}{\delta\rho_t}.
\]
A curve of probability densities
$\rho\in L^1((0,T)\times\mathcal X)$ is a weak solution of
\eqref{eq:continuum_levitas} if, for every
$\varphi\in C_c^\infty([0,T)\times\mathcal X)$,
\begin{align}
\int_0^T
\int_{\mathcal X}
\rho_t\,\partial_t\varphi
\,dx\,dt
&+
\int_0^T
\int_{\mathcal X}
\rho_t\nabla\mu_t\cdot\nabla\varphi
\,dx\,dt
\nonumber\\
&+
\int_{\mathcal X}
\rho_0(x)\varphi(x,0)\,dx
=
0.
\label{eq:weak_solution}
\end{align}
The flux is understood in a weakly integrable sense, so that
the second integral is well-defined.
\end{definition}

The positive sign in front of the flux term in
\eqref{eq:weak_solution} corresponds to the Wasserstein
gradient-flow equation
\[
\partial_t\rho_t
=
\nabla\cdot
\left(
\rho_t\nabla
\frac{\delta\mathcal J}{\delta\rho_t}
\right).
\]
It results from transferring the spatial divergence onto the
test function. The negative sign instead appears in the
corresponding energy-dissipation identity.


\begin{theorem}[Energy Dissipation for Weak Solutions]
\label{thm:weak_solution_energy_dissipation}
Suppose that $\rho_t$ is a sufficiently regular solution of
\eqref{eq:continuum_levitas}, that the relevant integrals are
finite, and that the boundary conditions eliminate the boundary
flux. Then
\begin{equation}
\frac{d}{dt}\mathcal J(\rho_t)
=
-
\int_{\mathcal X}
\rho_t
\left|
\nabla\frac{\delta\mathcal J}{\delta\rho_t}
\right|^2
\,dx
\leq 0.
\label{eq:weak_solution_energy_dissipation}
\end{equation}
Consequently,
\begin{equation}
\mathcal J(\rho_T)
+
\int_0^T
\int_{\mathcal X}
\rho_t
\left|
\nabla\frac{\delta\mathcal J}{\delta\rho_t}
\right|^2
\,dx\,dt
=
\mathcal J(\rho_0).
\label{eq:weak_solution_energy_identity}
\end{equation}
For weak solutions for which only lower semicontinuity and an
energy-dissipation inequality are available, equality in
\eqref{eq:weak_solution_energy_identity} is replaced by
\begin{equation}
\mathcal J(\rho_T)
+
\int_0^T
\int_{\mathcal X}
\rho_t
\left|
\nabla\frac{\delta\mathcal J}{\delta\rho_t}
\right|^2
\,dx\,dt
\leq
\mathcal J(\rho_0).
\label{eq:weak_solution_energy_inequality}
\end{equation}
\end{theorem}

\begin{proof}
For a sufficiently regular solution, the chain rule gives
$
\frac{d}{dt}\mathcal J(\rho_t)
=
\int_{\mathcal X}
\frac{\delta\mathcal J}{\delta\rho_t}
\,\partial_t\rho_t\,dx.
$
Substituting the continuum evolution equation yields
$
\frac{d}{dt}\mathcal J(\rho_t)
=
\int_{\mathcal X}
\mu_t
\nabla\cdot(\rho_t\nabla\mu_t)\,dx.
$
Integration by parts, together with the assumed boundary
conditions, gives
$
\frac{d}{dt}\mathcal J(\rho_t)
=
-
\int_{\mathcal X}
\rho_t|\nabla\mu_t|^2\,dx
\leq 0.
$
Integrating this identity over $[0,T]$ gives
\eqref{eq:weak_solution_energy_identity}. For weak solutions,
the corresponding inequality follows under the usual
lower-semicontinuity and approximation assumptions.
\end{proof}


\begin{theorem}[Existence of Weak Solutions]
\label{thm:existence_weak_solutions}
Suppose that $\mathcal X$ is a bounded domain and that
$\rho_0\in\mathcal P(\mathcal X)$ has finite initial energy
$\mathcal J(\rho_0)$. Assume that the continuum LEvITAS equation
admits a sequence of regularized approximate solutions
$\{\rho_t^\varepsilon\}_{\varepsilon>0}$ satisfying the following
properties:

\begin{enumerate}
\item The approximate solutions preserve probability mass and
satisfy the prescribed initial condition in the limit.

\item The energies are uniformly bounded from below and the
approximate solutions satisfy a uniform energy-dissipation
estimate:
\[
\sup_{\varepsilon>0}
\sup_{t\in[0,T]}
\mathcal J(\rho_t^\varepsilon)<\infty,
\]
and
\[
\sup_{\varepsilon>0}
\int_0^T\int_{\mathcal X}
\rho_t^\varepsilon
\left|
\nabla
\frac{\delta\mathcal J}{\delta\rho_t^\varepsilon}
\right|^2
\,dx\,dt
<\infty.
\]

\item These estimates provide sufficient compactness to extract
a subsequence converging to a probability-density curve
$\rho_t$, and the associated fluxes admit a limit that can be
identified with
$\rho_t\nabla(\delta\mathcal J/\delta\rho_t)$ in the weak
formulation.

\item The energy functional $\mathcal J$ is lower
semicontinuous under the convergence obtained above.
\end{enumerate}

Then, after passing to a subsequence, the approximate solutions
converge to a weak solution of
\eqref{eq:continuum_levitas} in the sense of
Definition~\ref{def:weak_solution}. The limiting solution
satisfies the energy-dissipation inequality
\eqref{eq:weak_solution_energy_inequality}.
\end{theorem}

\begin{proof}
The proof follows the compactness method for regularized
gradient-flow equations. The uniform energy and dissipation
bounds yield the required compactness of the approximate
solutions and their fluxes. Thus, along a subsequence,
the approximate densities converge to a limiting probability
density curve $\rho_t$.

By the assumed identification of the limiting flux, the
regularized weak formulations pass to the limit and give
\eqref{eq:weak_solution}. Lower semicontinuity of
$\mathcal J$, together with the uniform energy-dissipation
estimate, yields
\eqref{eq:weak_solution_energy_inequality}.
Therefore, the limiting curve is a weak solution satisfying
the asserted energy-dissipation bound.
\end{proof}

\subsection{Interpretation}
\label{subsec:pde_interpretation}

The continuum formulation establishes LEvITAS as a variational PDE system governing probability evolution on infinite-dimensional spaces. 
The model naturally combines:
 transport equations,
 gradient-flow dynamics,
 diffusion systems,
 Fokker--Planck evolution,
 entropy dissipation,
 continuum variational refinement,
 measure-valued geometric flows. 
Consequently, LEvITAS extends beyond conventional generative modeling into a continuum field theory of structured probabilistic evolution. 
The generalized refinement process may therefore be interpreted simultaneously as:
 variational optimization,
 geometric transport,
 stochastic probability evolution,
 continuum dynamical refinement,
 infinite-dimensional gradient flow.

\subsection{Main PDE and Continuum Formulation Results}
\label{subsec:pde_connections_main_results}

The PDE formulation developed in this section provides the continuum counterpart of the operator-theoretic and variational results established throughout the LEvITAS model. 
In particular:
The continuity equation
\eqref{eq:levitas_continuity}
and continuum gradient flow
\eqref{eq:continuum_levitas}
extend the refinement dynamics introduced in
Section~\ref{sec:iterative_topological_refinement}. 
The Wasserstein-gradient structure connects directly with the transport-geometric formulation established in
Section~\ref{sec:optimal_transport_wasserstein_geometry}. 
The Fokker--Planck dynamics
\eqref{eq:fokker_planck_levitas}
provide the continuum stochastic realization of the stochastic lifting systems developed in
Section~\ref{sec:measure_theoretic_foundations}. 
The entropy dissipation law of
Theorem~\ref{thm:entropy_dissipation}
provides a continuum analogue of the variational energy dissipation principles established for refinement trajectories in
Section~\ref{sec:iterative_topological_refinement}. 
The weak solution model of
Definition~\ref{def:weak_solution}
and
Theorem~\ref{thm:existence_weak_solutions}
supplies the analytical existence theory underlying continuum LEvITAS evolution systems.
The equilibrium characterization of
Theorem~\ref{thm:equilibrium_characterization}
links stationary continuum states with variational minimizers of the generalized LEvITAS functional introduced in
Section~\ref{sec:unified_reduction_theory}. 
Consequently, the PDE formulation unifies:
 variational refinement,
 Wasserstein geometry,
 stochastic probability evolution,
 measure-valued dynamics,
 continuum transport systems,
 infinite-dimensional optimization

within a single continuum operator-theoretic formulation of LEvITAS.

\subsection{Unified LEvITAS Synthesis Theory}
\label{sec:unified_synthesis}

The preceding sections established the LEvITAS model through:
 generalized lifting operators (Section~\ref{subsec:unified_operator_interpretation}),
 structure-preserving mappings (Section~\ref{subsec:structural_consistency}),
 path-space lifting systems (Section~\ref{sec:path_lifting_dynamical_systems}),
 variational refinement dynamics (Section~\ref{subsec:variational_dynamical_optimization}),
 measure-theoretic probability evolution (Section~\ref{sec:measure_theoretic_foundations}),
 optimal transport geometry (Section~\ref{sec:optimal_transport_wasserstein_geometry}),
 continuum PDE refinement systems (Section~\ref{sec:pde_continuum_levitas}). 
We now synthesize these components into a unified mathematical formulation. 
The purpose of this section is to establish:
global consistency of the LEvITAS model, compatibility between operator, 
probabilistic, geometric, and continuum formulations, reduction consistency 
from Section~\ref{sec:unified_reduction_theory}, existence of generalized 
variational refinement evolution. 
Within the generalized LEvITAS formulation, operator-theoretic mappings, 
probabilistic transport, variational evolution, and continuum dynamics are 
interpreted as mutually compatible realizations of a single structural 
generative system.

\subsection{Global LEvITAS System}
\label{subsec:global_levitas_system}

\begin{definition}[Global LEvITAS Operator System]
A generalized LEvITAS system is a tuple \\
$
\mathfrak L
=
(
\mathcal Z,
\mathcal X_{syn},
\mathcal X_{real},
P_Z,
\mathcal L_\theta,
\pi,
\mathcal S,
\mathcal J
),
$
where:
$
(\mathcal Z,\Sigma_Z,P_Z)
$ 
is a latent probability space, 
$
\mathcal L_\theta:
\mathcal Z
\rightarrow
\mathcal X_{syn}
$ 
is a measurable lifting operator (cf.~Definition~\ref{def:measurable_lifting_operator}), 
$
\pi:
\mathcal X_{syn}
\rightarrow
\mathcal X_{real}
$ 
is a measurable projection operator (cf.~Definition~\ref{def:measurable_projection_operator}), 
$
\mathcal S:
\mathcal X_{real}
\rightarrow
\mathcal M
$ 
is a measurable structure operator, 
$
\mathcal J:
\mathcal P(\mathcal X_{real})
\rightarrow
\mathbb R\cup\{+\infty\}
$ 
is a proper lower semicontinuous variational functional. 
\end{definition}

The induced observable probability measure is
$\mu_\theta=(\pi\circ\mathcal L_\theta)_\#P_Z,$ 
as introduced in Section~\ref{sec:measure_theoretic_foundations}. 
The generalized LEvITAS objective functional is
\begin{align}
\mathcal J(\theta)
&=
D(\mu_\theta,P_{data})
+
\lambda_1\mathcal D_{feat}
\nonumber\\
&\quad
+
\lambda_2\mathcal D_{struct}
+
\lambda_3\mathcal D_{dyn}.
\label{eq:global_levitas_energy_updated}
\end{align}

Equation
\eqref{eq:global_levitas_energy_updated}
synthesizes the variational formulations developed in
Sections~\ref{subsec:variational_formulation},
\ref{sec:measure_theoretic_foundations},
and~\ref{sec:optimal_transport_wasserstein_geometry}. 
The generalized refinement problem consists of minimizing
$\mathcal J$ over admissible lifting configurations.

\subsection{Global Variational Evolution}
\label{subsec:global_variational_evolution}

The global LEvITAS model induces variational evolution on spaces of probability measures.


\begin{theorem}[Global LEvITAS Variational Evolution]
\label{thm:global_variational_evolution}
Suppose:
$
\mathcal Z,
\mathcal X_{syn},
\mathcal X_{real}
$ 
are Polish spaces, 
$\mathcal L_\theta$ and $\pi$ are measurable, $\mathcal J$ is proper, coercive, and lower semicontinuous, admissible probability measures are tight. 
Then the LEvITAS model defines a well-posed variational evolution system on
$\mathcal P(\mathcal X_{real}).$
\end{theorem}

\begin{proof}
Measurability of
$
\mathcal L_\theta
$ 
and
$
\pi
$ 
ensures existence of the pushforward measure
$
\mu_\theta
=
(\pi\circ\mathcal L_\theta)_\#P_Z
$
by Theorem~\ref{thm:stability_pushforward_measures}. 
Tightness of admissible measures together with Prokhorov compactness (Theorem~\ref{thm:prokhorov_compactness}) yields weak sequential compactness. 
Lower semicontinuity and coercivity of $\mathcal J$ imply existence of minimizing refinement sequences via the direct method of the calculus of variations. 
Consequently, the generalized LEvITAS model induces a well-defined variational evolution system on probability space.
\end{proof}

\subsection{Transport-Geometric Consistency}
\label{subsec:transport_geometric_consistency}

The probabilistic and continuum formulations of LEvITAS are mutually compatible through Wasserstein transport evolution developed in Section~\ref{sec:optimal_transport_wasserstein_geometry}.


\begin{theorem}[Transport-Geometric Consistency]
\label{thm:transport_geometric_consistency}
Suppose: 
$\mu_t$ satisfies the Wasserstein gradient flow
$
\partial_t\mu_t
=
\nabla\cdot
\left(
\mu_t
\nabla
\frac{\delta\mathcal J}{\delta\mu_t}
\right),
$ 
as introduced in Definition~\ref{def:wasserstein_gradient_flow}, 
$\rho_t$ is the density representation of $\mu_t.$ 
Then $\rho_t$ satisfies the continuum LEvITAS PDE
\begin{equation}
\partial_t\rho_t
=
\nabla\cdot
\left(
\rho_t
\nabla
\frac{\delta\mathcal J}{\delta\rho_t}
\right).
\label{eq:global_continuum_updated}
\end{equation}
\end{theorem}

\begin{proof}
The Wasserstein gradient-flow equation induces a continuity-equation evolution for probability densities (cf.~Section~\ref{sec:pde_continuum_levitas}). 
Defining the refinement velocity field by
$
v_t
=
-
\nabla
\frac{\delta\mathcal J}{\delta\rho_t},
$ 
the continuity equation becomes
$
\partial_t\rho_t
+
\nabla\cdot(\rho_tv_t)
=
0.$ 
Substituting the velocity representation yields
\eqref{eq:global_continuum_updated}.
\end{proof}

\subsection{Energy Dissipation Structure}
\label{subsec:energy_dissipation_structure}

The global LEvITAS evolution possesses a monotone variational dissipation structure compatible with the continuum dynamics of Section~\ref{sec:pde_continuum_levitas}.


\begin{theorem}[Global Energy Dissipation]
\label{thm:global_energy_dissipation}
Suppose $\rho_t$ satisfies
\eqref{eq:global_continuum_updated}. 
Then
\begin{equation}
\frac{d}{dt}\mathcal J(\rho_t)
=
-
\int_{\mathcal X}
\rho_t
\left|
\nabla
\frac{\delta\mathcal J}{\delta\rho_t}
\right|^2dx
\leq0.
\label{eq:global_energy_decay_updated}
\end{equation}
\end{theorem}

\begin{proof}
Applying the variational chain rule gives
$\frac{d}{dt}\mathcal J(\rho_t)
=
\int
\frac{\delta\mathcal J}{\delta\rho_t}
\partial_t\rho_t\,dx.
$

Substituting
\eqref{eq:global_continuum_updated},
$
\frac{d}{dt}\mathcal J(\rho_t)
=
\int
\frac{\delta\mathcal J}{\delta\rho_t}
\nabla\cdot
\left(
\rho_t
\nabla
\frac{\delta\mathcal J}{\delta\rho_t}
\right)dx.
$

Integration by parts yields
$
\frac{d}{dt}\mathcal J(\rho_t)
=
-
\int
\rho_t
\left|
\nabla
\frac{\delta\mathcal J}{\delta\rho_t}
\right|^2dx.
$

Therefore the LEvITAS energy is nonincreasing along refinement evolution.

Equation
\eqref{eq:global_energy_decay_updated}
is the global dissipation law corresponding to the local continuum dissipation identity established in Section~\ref{sec:pde_continuum_levitas}.
\end{proof}

\subsection{Reduction Hierarchy}
\label{subsec:reduction_hierarchy}

The generalized LEvITAS system contains classical generative architectures as constrained reductions, as developed in Section~\ref{sec:unified_reduction_theory}.


\begin{theorem}[LEvITAS Reduction Hierarchy]
\label{thm:levitas_reduction_hierarchy}
Specialized generative systems arise as constrained reductions of LEvITAS under structural restrictions on:
$
\mathcal L_\theta,
\qquad
\pi,
\qquad
\mathcal S,
\qquad
\mathcal J.
$
\end{theorem}

\begin{proof}
The generalized LEvITAS model contains:
 deterministic liftings,
 stochastic liftings,
 path-space transport systems,
 variational latent regularization,
 Wasserstein refinement evolution,
 diffusion refinement dynamics. 
Imposing restrictions on admissible lifting operators, discrepancy functionals, structural penalties, and refinement dynamics yields reduced systems corresponding to specialized generative architectures. 
The reductions constructed in Section~\ref{sec:unified_reduction_theory} include:
 GAN systems,
 variational autoencoders,
 diffusion systems,
 normalizing flows,
 energy-based models. 
Since each reduction preserves the underlying variational structure, every reduced system remains embedded within the generalized LEvITAS model.
\end{proof}

\subsection{Equilibrium and Stability}
\label{subsec:equilibrium_stability}

Equilibrium configurations of LEvITAS correspond to stationary probability measures minimizing the generalized variational energy.


\begin{theorem}[LEvITAS Equilibrium Principle]
\label{thm:levitas_equilibrium_principle}
Suppose: $\mathcal J$ is convex, 
$\rho_\infty$ satisfies
$
\nabla\cdot
\left(
\rho_\infty
\nabla
\frac{\delta\mathcal J}{\delta\rho}
\right)
=
0.
$ 
Then
$\rho_\infty$ is a global minimizer of $\mathcal J.$
\end{theorem}

\begin{proof}
At equilibrium,
$
\nabla
\frac{\delta\mathcal J}{\delta\rho}
=
0.
$ 
Hence the first variation of
$\mathcal J$ vanishes. 
Convexity of $\mathcal J$ therefore implies global minimality of $\rho_\infty.$ 
This equilibrium condition coincides with the stationary continuum formulation derived in Section~\ref{sec:pde_continuum_levitas}.
\end{proof}

\subsection{Unified LEvITAS Meta-Theorem}
\label{subsec:meta_theorem}

We now state the principal synthesis theorem of the generalized LEvITAS model.

\begin{theorem}[Unified LEvITAS Meta-Theorem]
\label{thm:unified_levitas_meta_theorem}
The LEvITAS model defines:
a measurable lifting-projection operator system, a variational probability evolution theory, a Wasserstein transport geometry, a continuum gradient-flow PDE system, a generalized reduction hierarchy for generative models. 
Furthermore, these formulations are mutually compatible under pushforward probability evolution.
\end{theorem}

\begin{proof}
The measurable operator formulation induces pushforward probability measures via Section~\ref{sec:measure_theoretic_foundations}. 
The variational formulation defines refinement evolution on probability spaces through the generalized functional
\eqref{eq:global_levitas_energy_updated}. 
The Wasserstein formulation induces transport-geometric dynamics through the constructions of Section~\ref{sec:optimal_transport_wasserstein_geometry}. 
The continuum PDE formulation follows from the continuity equation associated with Wasserstein gradient flow, developed in Section~\ref{sec:pde_continuum_levitas}. 
Compatibility follows from equivalence between the measure-valued and density-valued representations of probability evolution.
\end{proof}

\newpage

\section{Appendix B: LEvITAS--Biological Validation and Clinical Insights Using ci-fGBD and MTP~\cite{Thakur2025ciFGBD}}

We investigate the biological validity and clinical relevance of the LEvITAS 
operator model across heterogeneous biomedical applications by employing the 
previously developed Cluster-Integrated Fast Generalized Bruhat Decomposition 
(ci-fGBD), a foundational model that establishes a connection with classical 
matrix factorization, together with its extended counterpart, the Longitudinal 
Cluster-Integrated Fast Generalized Bruhat Decomposition (Lci-fGBD). The 
latter provides a unified, modular, and longitudinally extended model that 
explicitly incorporates Bruhat decomposition, extending the methodological 
scope beyond that presented in the preprint. The resulting feature-level 
patterns are examined using our Merkmal Treiber Plot (MTP) visualization 
tool~\cite{Thakur2025ciFGBD}, while the associated low-dimensional clustering 
projections provide a complementary sample-level view of latent organization. 
Our objective is to assess whether LEvITAS-generated synthetic representations 
preserve latent organization, clinically meaningful feature relationships, and 
multimodal dependencies characteristic of real-world biomedical cohorts, and 
to contextualize these findings against the clustering structure obtained 
using ci-fGBD and baseline non-negative matrix factorization (NMF).

The ci-fGBD framework and its associated MTP approach provide complementary 
analytical capabilities for investigating heterogeneous biomedical 
representations, feature hierarchies, and modality-dependent relationships. In 
the present study, MTP is employed as an established visualization and 
interpretability approach to assess the structural correspondence between real 
and LEvITAS-generated cohorts, extending the evaluation beyond conventional 
distributional similarity towards the preservation of coordinated feature 
organization. In addition, the clustering projections compare the ci-fGBD 
representation, denoted by $U$, with the baseline NMF representation, denoted 
by $W$, using two-dimensional principal-component (PC1--PC2) views. MTP 
therefore provides a feature-contribution perspective, whereas the clustering 
plots show the arrangement of observations in the corresponding low-
dimensional representations.

To complement the structural assessment, three quantitative similarity 
measures are considered: mean absolute error (MAE), maximum mean discrepancy 
(MMD), and Kullback--Leibler (KL) divergence. MAE evaluates numerical 
differences between corresponding representation values, MMD quantifies 
distributional discrepancies using kernel-based comparisons, and KL divergence 
measures differences between estimated probability distributions. Together, 
these metrics provide complementary assessments of numerical and 
distributional fidelity, while MTP visualization examines latent feature 
hierarchies and cross-modal relationships. These metrics are associated with 
the real-versus-synthetic comparisons and should be interpreted separately 
from the ci-fGBD--NMF clustering projections, which are descriptive 
visualizations rather than additional similarity scores.

The validation encompasses five biomedical applications: (i) endometriosis 
laparoscopic imaging, (ii) antral follicle count (AFC) ultrasound phenotyping, 
(iii) amyotrophic lateral sclerosis (ALS) multi-omic profiling, (iv) 
Alzheimer's disease (AD) multimodal clinical--cognitive representation 
learning, and (v) antimicrobial resistance (AMR) genomic--morphological 
modelling. The following disease-specific subsections present the 
corresponding real-versus-synthetic comparisons, quantitative similarity 
results, MTP visualizations, ci-fGBD-versus-NMF clustering projections, and 
their biological or clinical relevance.

The clustering panels should be interpreted cautiously. Each method is 
represented in its own latent coordinate system, and the displayed principal 
components may differ in scale, sign, orientation, and variance structure. 
Accordingly, the comparison focuses on the broad within-panel grouping 
geometry, apparent separation, overlap, dispersion, and visible outliers, 
rather than on absolute coordinate magnitudes or direct point-to-point 
distances across methods. The colours identify groups within each panel and 
are not assumed to encode matched cluster identities between ci-fGBD and NMF.

Collectively, these investigations examine the ability of LEvITAS to generate 
structurally consistent synthetic biomedical representations while preserving 
dominant latent organization and clinically relevant multimodal relationships 
across diverse biomedical domains. The accompanying ci-fGBD--NMF comparisons 
provide an additional descriptive assessment of how observations are organized 
in two alternative latent representations, complementing rather than replacing 
the MTP and quantitative real-versus-synthetic evaluations.

\begin{figure}[!ht]
\centering

\begin{subfigure}[b]{0.58\textwidth}
    \centering
    \includegraphics[width=\linewidth]{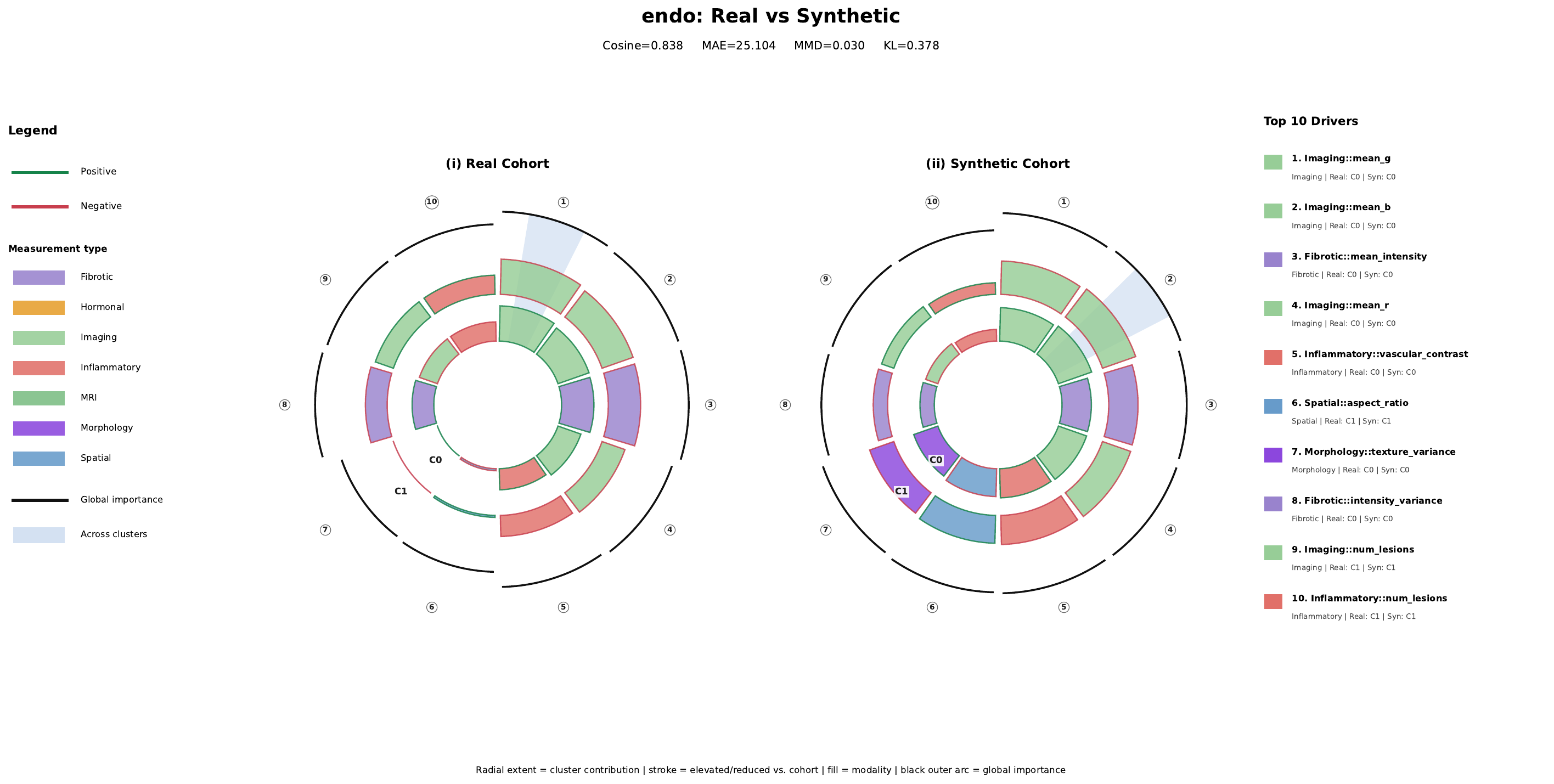}
\caption{MTP feature-contribution summary for the endometriosis cohort, comparing dominant feature patterns in the real GLENDA cohort and LEvITAS-generated synthetic data. The comparison illustrates the preservation of major feature patterns and variation across selected feature groups.}
\label{fig:mtp_endometriosis_cifgbd}
\end{subfigure}
\hfill
\begin{subfigure}[b]{0.38\textwidth}
    \centering
    \includegraphics[width=\linewidth]{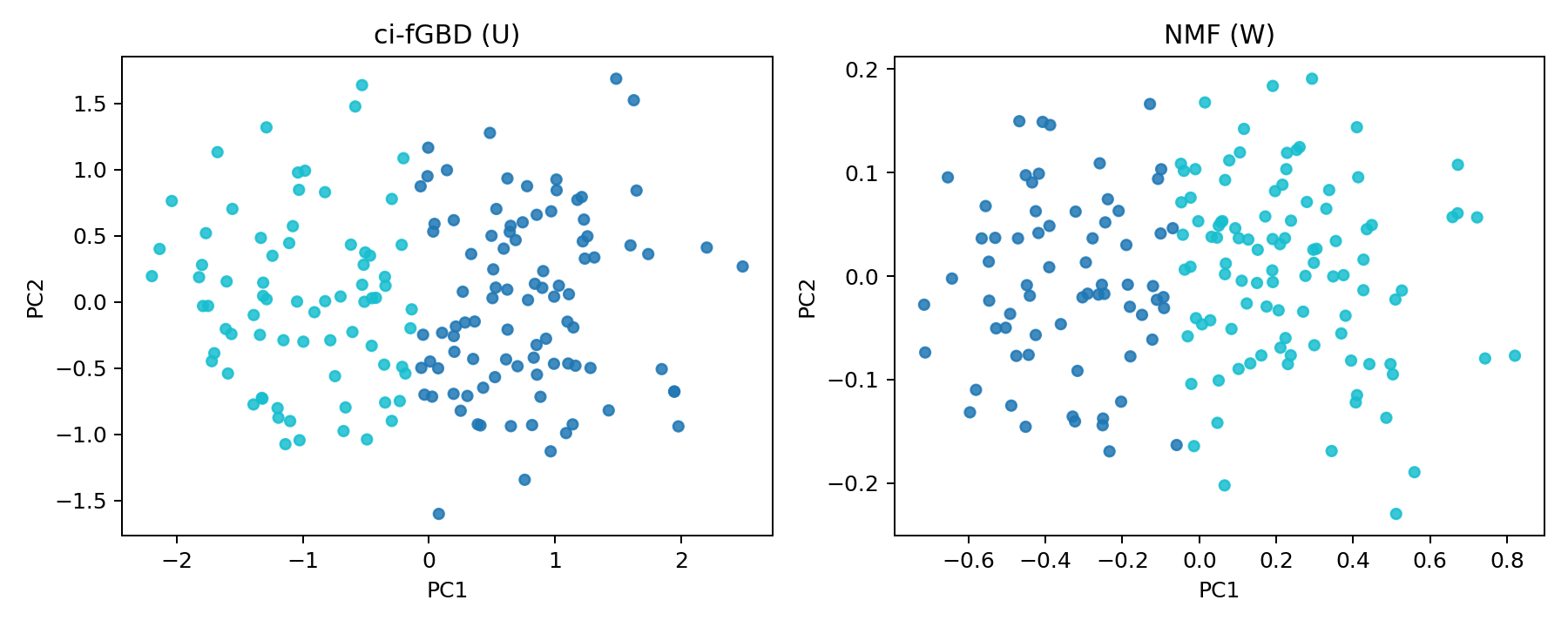}
    \caption{PC1--PC2 clustering projections of the endometriosis cohort using the ci-fGBD $U$ and baseline NMF $W$ representations. Both projections exhibit a broad two-group arrangement, with differences in within-group dispersion and overlap near the central boundary. Colour assignments are method-specific, and the projected coordinate systems are independent; therefore, absolute coordinates and distances should not be compared directly.}
    \label{fig:cluster_endometriosis}
\end{subfigure}
    \caption{Complementary MTP and latent-space clustering analyses of the endometriosis cohort. Panel (a) presents the feature-contribution comparison between real GLENDA and LEvITAS-generated synthetic data, illustrating feature-level fidelity and variation. Panel (b) presents the PC1--PC2 projections of the ci-fGBD $U$ and baseline NMF $W$ representations, highlighting broad latent grouping and differences in sample dispersion and overlap. Together, these analyses provide complementary perspectives on feature-level fidelity and sample-level latent organization in synthetic endometriosis data.}
    \label{fig:endo_exp_analysis}
\end{figure}

\paragraph*{Enodometriosis-Complementary Clustering Comparison:}
The endometriosis projections show two broad colour-coded groups in both the 
ci-fGBD $U$ and baseline NMF $W$ representations Figure~\ref{fig:endo_exp_analysis}(b), with the principal visible 
division occurring along PC1. In the ci-fGBD panel, the groups occupy largely 
opposing sides of the horizontal axis and exhibit appreciable variation along 
PC2, including observations extending beyond the densest portions of each 
group. The NMF panel retains a similar broad left--right organization but 
displays a more compact point distribution, with some proximity between the 
groups near the central boundary. The colour assignments are specific to each 
panel and should not be treated as matched cluster identities across methods. 
Overall, both decompositions expose a broad two-group arrangement, while 
representing within-group spread and the boundary region differently. This 
visual comparison complements the MTP feature-contribution summary by 
providing a sample-level view of latent organization. However, without 
clinical annotations and independent cluster-validity assessment, the groups 
cannot be assigned to specific lesion phenotypes, disease stages, or 
clinically validated endometriosis subtypes.  

Overall, the combined MTP and clustering analyses 
Figure~\ref{fig:endo_exp_analysis} provide complementary views 
of how LEvITAS and ci-fGBD represent the structural, morphological, and 
imaging-associated characteristics of endometriosis. The MTP comparison 
suggests that several dominant feature patterns in the real GLENDA cohort are 
retained in the LEvITAS-generated synthetic data, while reduced variation in 
selected feature groups highlights the need to assess whether clinically 
meaningful patient heterogeneity is also preserved. The ci-fGBD--NMF 
projections additionally show two broad colour-coded groups in both latent 
representations, with differences in their visible dispersion and overlap near 
the central boundary. These projections provide a sample-level view of latent 
grouping, but do not establish that the groups correspond to validated lesion 
phenotypes or clinical disease categories. Accordingly, evaluating synthetic 
data quality requires more than numerical distributional similarity or visual 
cluster separation; it also requires examining structural fidelity, spatial 
organization, relationships among disease-associated features, and biological 
consistency. Together, these complementary analyses provide a more 
comprehensive assessment of the extent to which LEvITAS-generated data retain 
meaningful endometriosis-associated characteristics and support future 
multimodal endometriosis research and precision gynecological analytics.

\begin{figure}[!ht]
\centering

\begin{subfigure}[b]{0.68\textwidth}
    \centering
    \includegraphics[width=\linewidth]{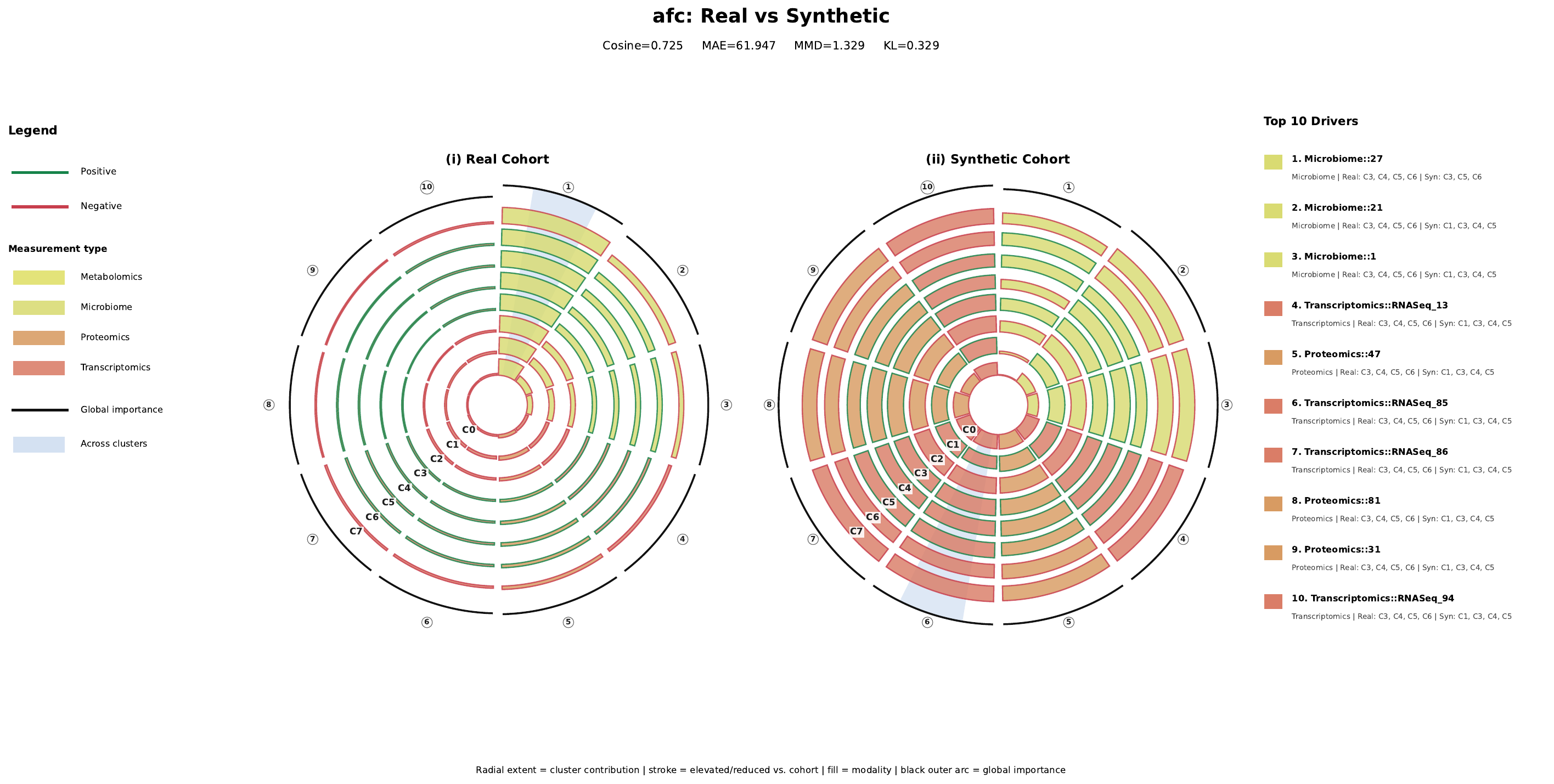}
    \caption{MTP feature-contribution summary for the AFC cohort, comparing dominant feature patterns in the real USOVA3D cohort and LEvITAS-generated synthetic data. The comparison illustrates the preservation of ovarian imaging- and morphology-related feature patterns, while highlighting differences in the distribution and concentration of selected feature groups.}
    \label{fig:mtp_afc}
        
\end{subfigure}
\hfill
\begin{subfigure}[b]{0.28\textwidth}
    \centering
    \includegraphics[width=\linewidth]{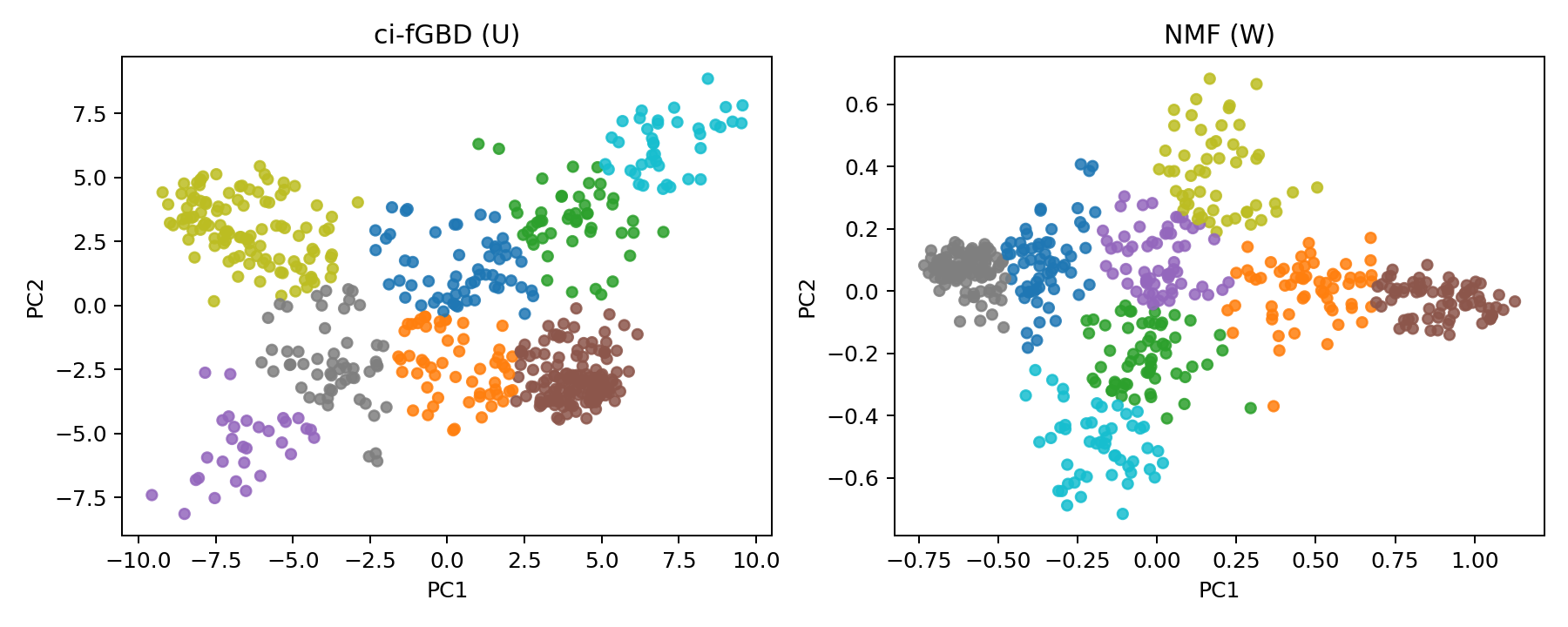}
    \caption{PC1--PC2 clustering projections of the AFC cohort using the ci-fGBD $U$ and baseline NMF $W$ representations. Eight colour-coded groups are visible in each panel, although their relative spatial arrangements, dispersion, and proximity to neighbouring groups differ between methods. The projected coordinate systems are method-specific; therefore, absolute coordinates, distances, and numerical axis ranges should not be compared directly. The colour assignments do not represent clinically established ovarian-reserve categories or follicular developmental stages.}
    \label{fig:cluster_afc_nmf}
\end{subfigure}

\caption{Complementary MTP and latent-space clustering analyses of the AFC cohort. Panel (a) presents the feature-contribution comparison between real USOVA3D and LEvITAS-generated synthetic data, illustrating the preservation of ovarian imaging- and morphology-related feature patterns. Panel (b) presents the PC1--PC2 projections of the ci-fGBD $U$ and baseline NMF $W$ representations, highlighting eight colour-coded groups and differences in their spatial organization, dispersion, and overlap. Together, these analyses provide complementary perspectives on feature-level fidelity and multi-group latent organization in synthetic AFC data. The observed patterns and groupings do not independently establish the preservation of biological follicular populations or clinically validated ovarian-reserve characteristics.}
\label{fig:afc_exp_analysis}

\end{figure}

\paragraph*{AFC-Complementary Clustering Comparison:}
The AFC projections show eight colour-coded groups in both the ci-fGBD $U$ and 
NMF $W$ representations Figure~\ref{fig:afc_exp_analysis}(b), indicating that the displayed latent spaces contain 
multiple groupings rather than a single compact population. In the ci-fGBD 
panel, the groups are distributed across a comparatively broad PC1--PC2 
region, with several peripheral groups and a number of groups approaching one 
another around the central region. The NMF panel also displays multiple 
groups, but their spatial arrangement is different and occupies a 
substantially more compact numerical coordinate range; several groups are 
organized along the horizontal direction, while others are positioned above or 
below the central region. These differences suggest that the two 
decompositions represent the multi-group organization of the AFC observations 
with different latent geometries. The group colours are not accompanied by 
clinical labels and should not be interpreted as direct ovarian-reserve 
categories, follicular developmental stages, or specific patient phenotypes. 
Read alongside the AFC MTP, the projections provide a complementary view of 
latent heterogeneity in ovarian ultrasound representations, but they do not 
establish that the groups correspond to distinct biological follicular 
populations or clinically meaningful ovarian-reserve strata. The difference in 
plotted coordinate ranges is not a direct measure of relative biological 
variability or model quality. 

Taken together, the AFC MTP and clustering analyses
Figure~\ref{fig:afc_exp_analysis} indicate that LEvITAS 
retains important ovarian imaging- and morphology-related feature patterns, 
while the ci-fGBD and NMF projections both reveal multiple colour-coded groups 
with differing spatial arrangements and degrees of overlap. The MTP results 
suggest that aspects of the original ovarian imaging representation are 
maintained without reproducing all of its concentrated latent activation 
patterns, while the clustering projections provide a complementary view of 
multi-group organization in the latent sample space. However, neither the 
preserved MTP patterns nor the displayed clusters establish that actual 
follicular structures, follicle counts, follicular sizes, or clinically 
meaningful ovarian-reserve characteristics have been preserved. Direct 
structural comparisons and validation against clinical ovarian-reserve 
measurements remain important next steps. Overall, these findings support the 
potential utility of LEvITAS for reproductive imaging augmentation and future 
fertility-related imaging research, while highlighting the need for further 
anatomical and clinical validation before establishing biological or clinical 
applicability.

\paragraph*{ALS-Complementary Clustering Comparison:}

\begin{figure}[!ht]
\centering

\begin{subfigure}[b]{0.68\textwidth}
    \centering
    \includegraphics[width=\linewidth]{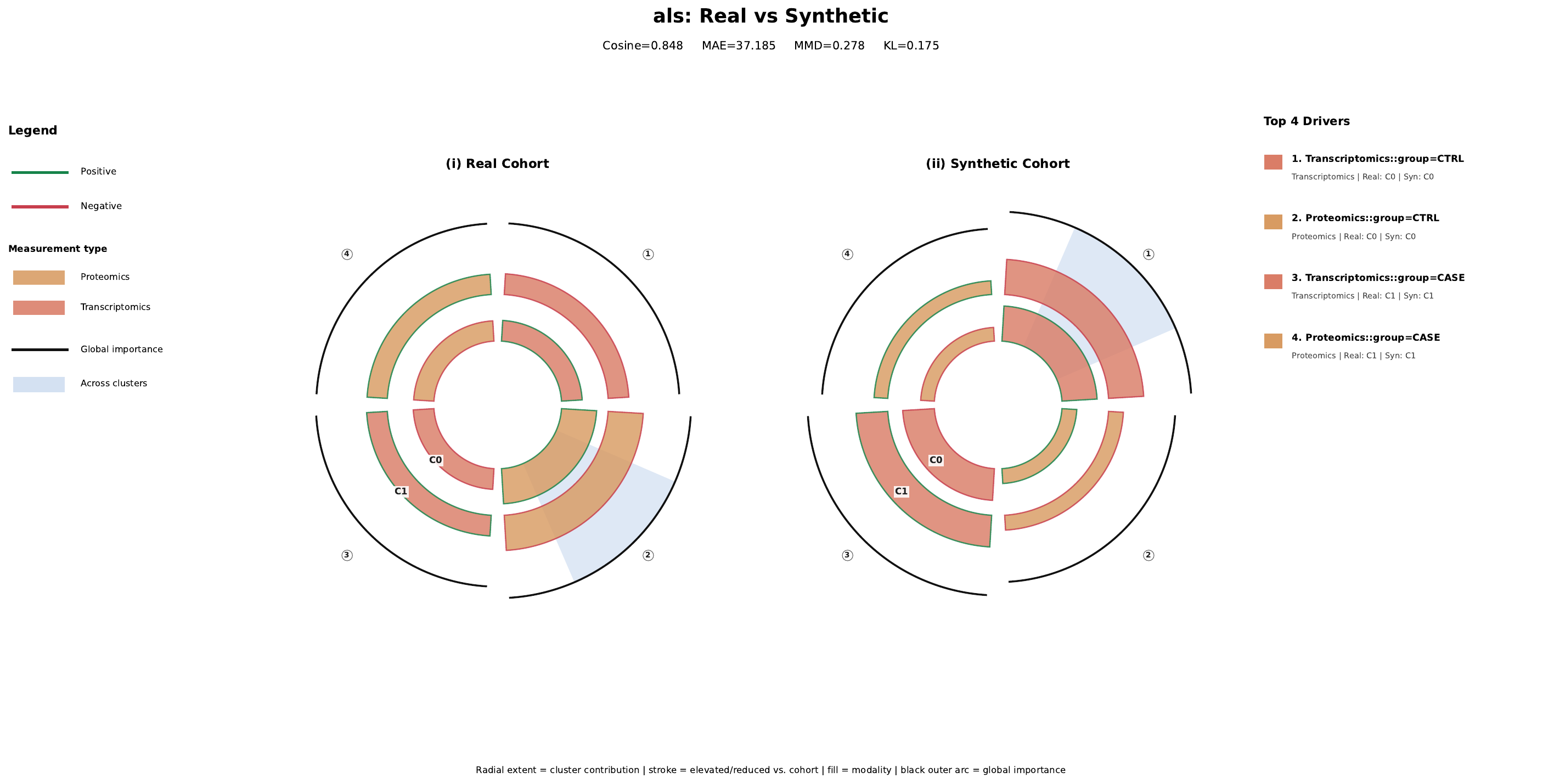}
    \caption{MTP feature-contribution summary for the ALS cohort, comparing dominant genomic and proteomic feature patterns in the real Neuromine cohort and LEvITAS-generated synthetic data. The comparison illustrates the extent to which the synthetic data retain dominant molecular contributors and global feature-importance patterns.}
    \label{fig:cluster_als_cifgbd}
    
\end{subfigure}
\hfill
\begin{subfigure}[b]{0.28\textwidth}
    \centering
    \includegraphics[width=\linewidth]{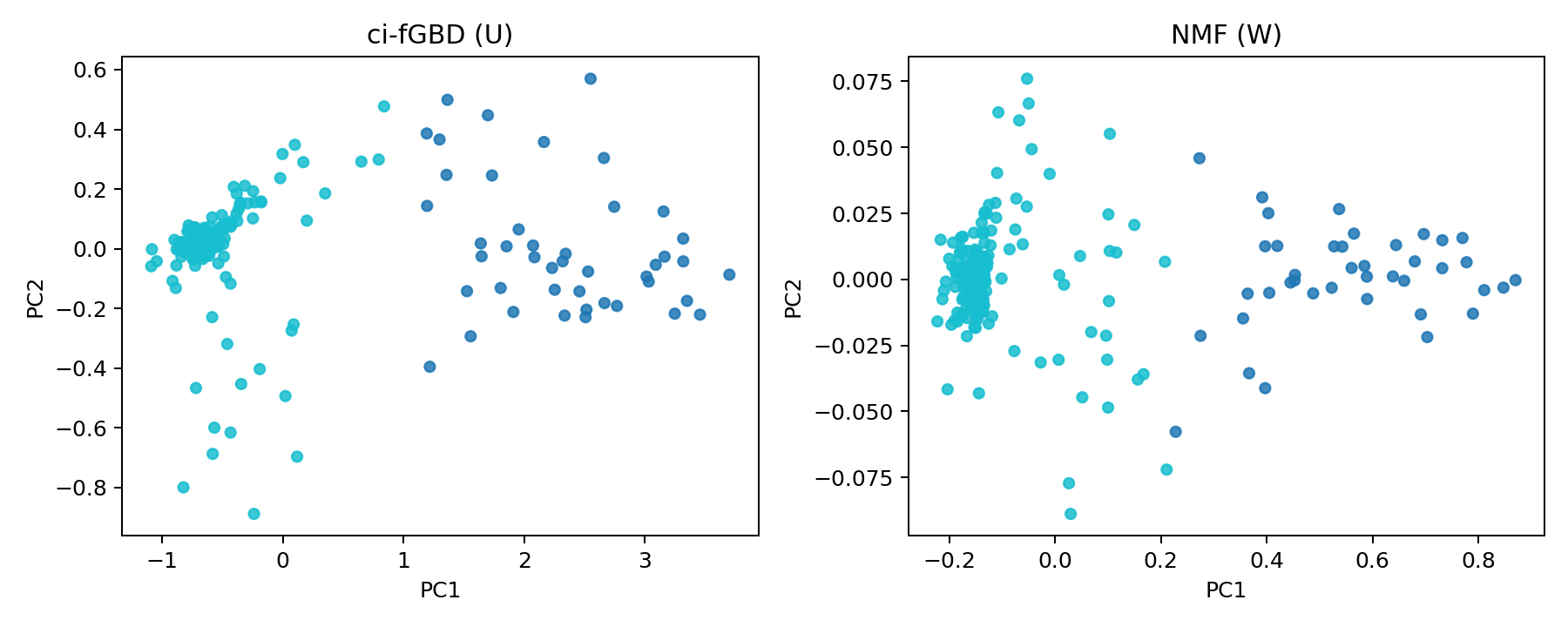}
    \caption{PC1--PC2 clustering projections of the ALS cohort using the ci-fGBD $U$ and baseline NMF $W$ representations. Both projections exhibit two broad colour-coded groups primarily separated along PC1, although differences are visible in within-group dispersion and the distribution of observations away from the densest point concentrations. The projected coordinate systems are method-specific; therefore, absolute coordinates, distances, and numerical axis ranges should not be compared directly. The colour assignments do not represent clinically or molecularly validated ALS subgroups.}
    \label{fig:cluster_als}
    
\end{subfigure}
\caption{Complementary MTP and latent-space clustering analyses of the ALS cohort. Panel (a) presents the feature-contribution comparison between real Neuromine and LEvITAS-generated synthetic data, illustrating the retention of dominant genomic and proteomic feature patterns. Panel (b) presents the PC1--PC2 projections of the ci-fGBD $U$ and baseline NMF $W$ representations, highlighting a broad two-group arrangement and differences in within-group dispersion and outlier structure. Together, these analyses provide complementary perspectives on molecular feature-level fidelity and sample-level latent organization in synthetic ALS data. The observed patterns do not independently establish the preservation of gene--protein interactions, biologically validated molecular relationships, or clinically defined ALS subgroups.}
\label{fig:als_exp_analysis}

\end{figure}

The ALS projections display two dominant colour-coded groups in both 
representations Figure!\ref{fig:als_exp_analysis}(b). In the ci-fGBD $U$
panel, one group forms a relatively dense concentration at negative PC1 
values, with additional observations extending 
vertically and horizontally away from the central concentration; the second 
group is distributed mainly across positive PC1 values and exhibits a broader 
spread. The NMF $W$ panel retains a broad separation along PC1, with a 
concentrated group on the negative side and a second group on the positive 
side, together with several observations displaced from the densest point 
clouds. Thus, both methods reveal a similar high-level two-group arrangement, 
while the apparent within-group dispersion and outlier patterns differ. The 
distinct PC scales mean that the more extended numerical range in the ci-fGBD 
panel cannot be interpreted directly as greater ALS-related biological 
variability. Moreover, the figure does not map the displayed groups to ALS 
subtype, disease stage, progression rate, or molecularly validated patient 
categories. The clustering projections therefore complement the MTP's genomic--
proteomic feature summaries by sample-level latent grouping, 
without independently establishing preservation of gene--protein interactions 
or clinically defined ALS subgroups. 

Overall, the combined ALS MTP and clustering analyses 
Figure~\ref{fig:als_exp_analysis} indicate that LEvITAS 
retains important genomic and proteomic feature representations, while the ci-
fGBD and NMF projections both display two broad colour-coded groups with 
differences in within-group dispersion and outlier structure. The continued 
representation of dominant molecular contributors and the alignment of global 
MTP importance patterns suggest that aspects of the original multi-omic 
organization are maintained; the clustering projections add a sample-level 
view of this latent organization. Nevertheless, these similarities do not 
establish preservation of gene--protein relationships, biologically meaningful 
interactions, clinically defined ALS subgroups, or longitudinal disease 
characteristics. Further validation through quantitative cross-omic 
comparisons, independent molecular measurements, diagnostic evaluation on real 
unseen patients, and appropriately characterized longitudinal and subtype-
specific datasets is therefore necessary. Taken together, these findings 
support the potential utility of LEvITAS for multi-omic data augmentation, 
biomarker investigation, and future systems-level neurodegenerative disease 
research, while highlighting the need for molecular, clinical, and 
longitudinal validation before establishing diagnostic, prognostic, or 
therapeutic applicability.

\begin{figure}[!ht]
\centering

\begin{subfigure}[b]{0.58\textwidth}
    \centering
    \includegraphics[width=\linewidth]{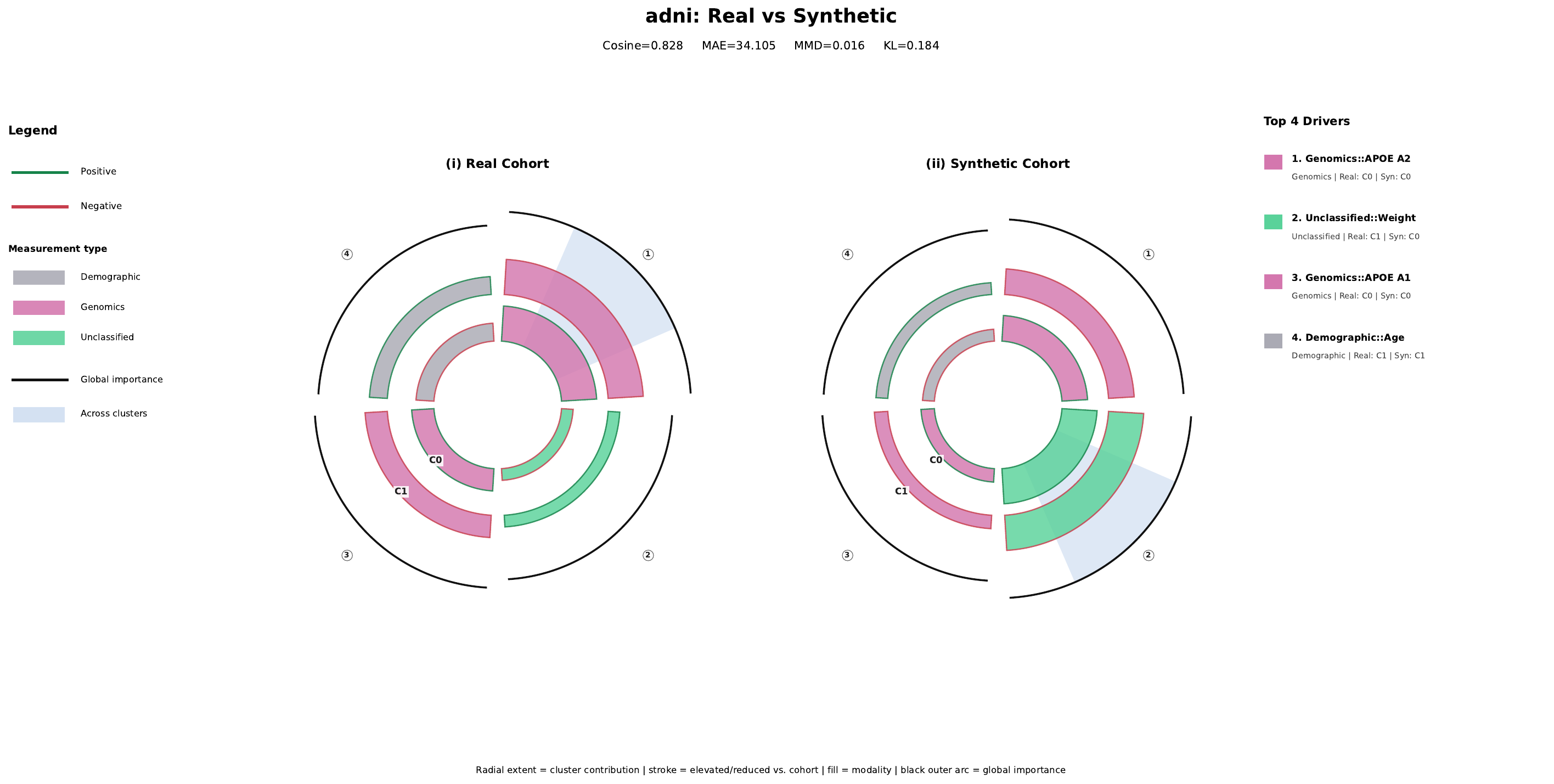}
    \caption{MTP feature-contribution summary for the ADNI cohort, comparing dominant demographic, clinical, and genetic feature patterns in the real ADNI cohort and LEvITAS-generated synthetic data. The comparison illustrates the extent to which the synthetic data retain dominant multimodal feature patterns and their relative contributions.}
    \label{fig:cluster_ad_cifgbd}
    
\end{subfigure}
\hfill
\begin{subfigure}[b]{0.28\textwidth}
    \centering
    \includegraphics[width=\linewidth]{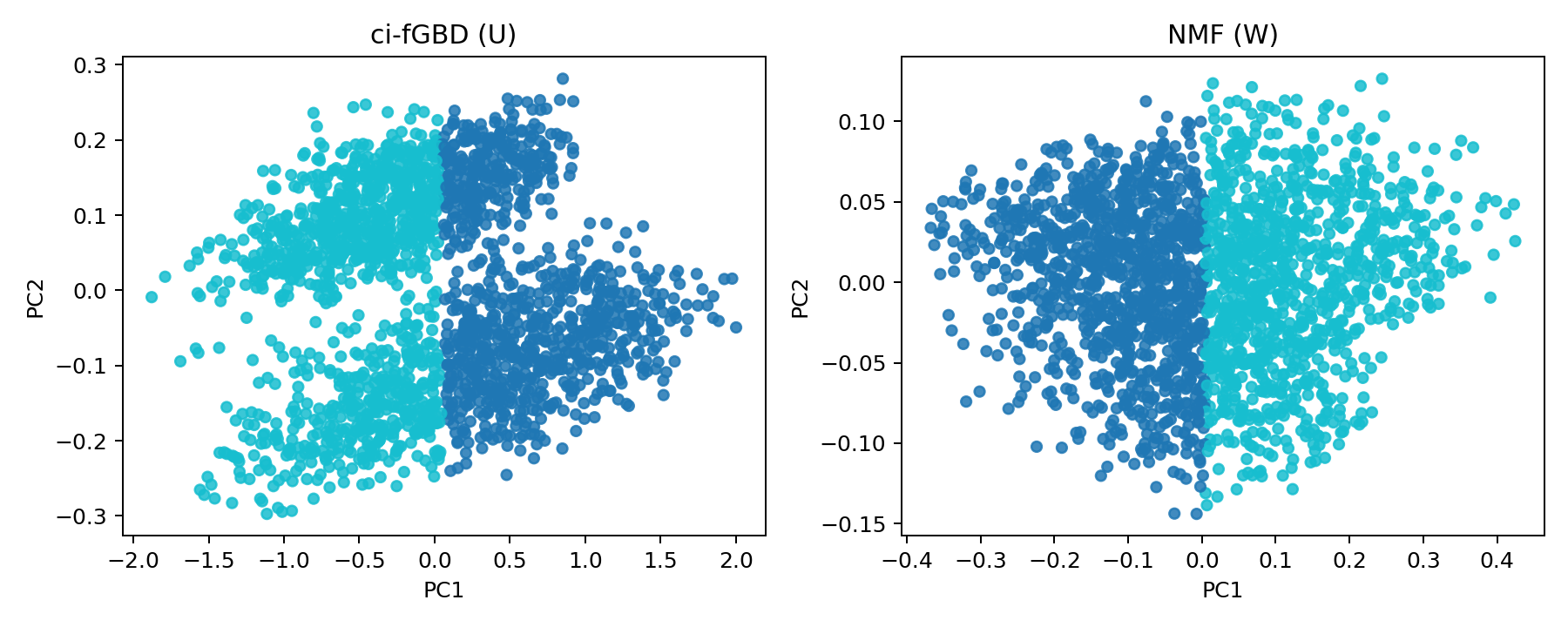}
    \caption{PC1--PC2 clustering projections of the ADNI cohort using the ci-fGBD $U$ and baseline NMF $W$ representations. Both projections exhibit a broad two-group arrangement primarily separated along PC1, although differences are visible in within-group vertical dispersion, internal organization, and proximity near the central boundary. The projected coordinate systems are method-specific; therefore, absolute coordinates, distances, and numerical axis ranges should not be compared directly. The colour assignments do not represent clinically established AD stages, diagnostic classes, or pathological subtypes.}
    \label{fig:cluster_ad}
    
\end{subfigure}
\caption{Complementary MTP and latent-space clustering analyses of the ADNI cohort. Panel (a) presents the feature-contribution comparison between real ADNI and LEvITAS-generated synthetic data, illustrating the retention of dominant demographic, clinical, and genetic feature patterns. Panel (b) presents the PC1--PC2 projections of the ci-fGBD $U$ and baseline NMF $W$ representations, highlighting a broad two-group arrangement and differences in within-group dispersion, internal organization, and boundary overlap. Together, these analyses provide complementary perspectives on multimodal feature-level fidelity and sample-level latent organization in synthetic ADNI data. The observed patterns do not independently establish the preservation of cross-modal biological relationships or clinically validated Alzheimer's disease stages, diagnostic classes, or pathological subtypes.}
\label{fig:ad_exp_analysis}

\end{figure}

\paragraph*{AD-Complementary Clustering Comparison:}
The ADNI projections show two broad colour-coded groups in both ci-fGBD $U$ 
and NMF $W$ Figure~\ref{fig:ad_exp_analysis}(b), with the main visible division occurring along PC1. In the ci-
fGBD panel, the groups occupy opposite sides of the projection and show 
noticeable vertical structure, with observations distributed across both 
positive and negative PC2 values. The NMF panel also presents a broad left--
right separation, although the groups approach one another near the PC1 
boundary and occupy a comparatively compact PC2 range. The broad partition is 
therefore visible in both representations, while the internal geometry, 
vertical spread, and degree of boundary proximity differ. These colour-coded 
groups are not identified in the figure as clinical AD stages, diagnostic 
classes, or longitudinal trajectories. Read alongside the AD MTP, the 
projections add a sample-level view of multimodal latent organization, but do 
not demonstrate that relationships among the shown modalities in the Figure~\ref{fig:ad_exp_analysis}(a). 

Overall, the combined AD MTP and clustering analyses
Figure~\ref{fig:ad_exp_analysis} indicate that LEvITAS 
retains important demographic, clinical, and genetic feature 
representations, while the ci-fGBD and NMF projections both show two broad 
colour-coded groups with differences in their internal spread and boundary 
regions. The continued representation of dominant multimodal contributors and 
the smoother organization of MTP feature patterns suggest that aspects of the 
original multimodal structure are maintained, while the clustering plots add a 
sample-level view of the latent organization. Nevertheless, these visual 
similarities do not establish preservation of relationships among demgraphic 
and clinical assessments, and genetic susceptibility, nor do they identify the 
projected groups as clinical AD stages or pathological subtypes. Further 
validation through quantitative cross-modal comparisons, independent biomarker 
and clinical assessments, diagnostic 
evaluation on real unseen patients, severity-specific comparisons, and 
appropriately characterized longitudinal datasets is therefore necessary. 
Taken together, these findings support the potential utility of LEvITAS for 
multimodal data augmentation, biomarker investigation, and future 
neurodegenerative disease research, while highlighting the importance of 
establishing biological, clinical, and longitudinal validity before extending 
these findings to diagnostic, prognostic, or precision neurology applications.

\begin{figure}[!ht]
\centering

\begin{subfigure}[b]{0.68\textwidth}
    \centering
    \includegraphics[width=\linewidth]{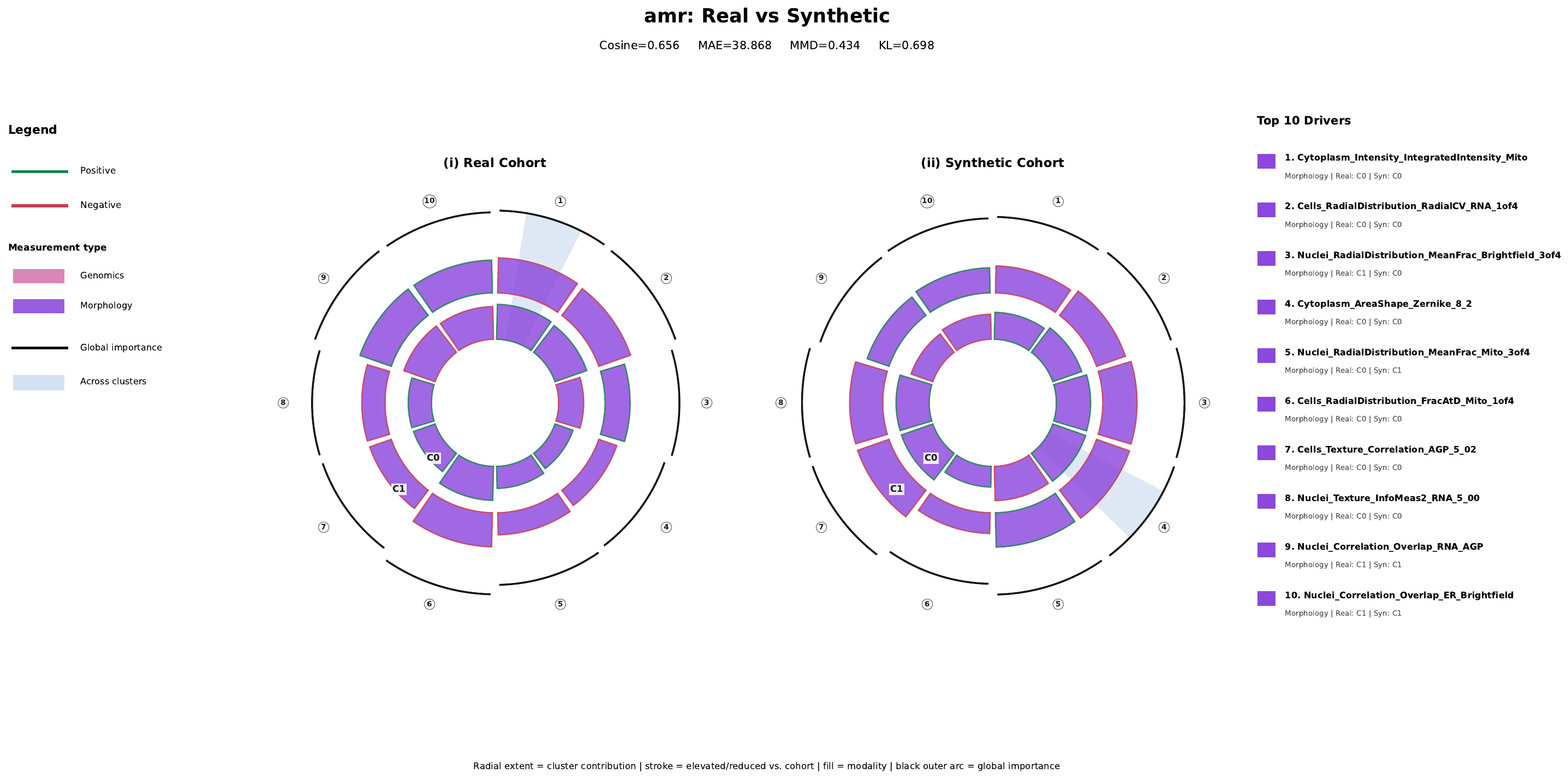}
    \caption{MTP feature-contribution summary for the AMR cohort, comparing dominant genomic k-mer signatures and Cell Painting-derived morphological feature patterns in the real AMR cohort and LEvITAS-generated synthetic data. The comparison illustrates the extent to which the synthetic data retain important genomic and morphology-associated feature contributions and their overall distribution.}
    \label{fig:mtp_amr}
    
\end{subfigure}
\hfill
\begin{subfigure}[b]{0.28\textwidth}
    \centering
    \includegraphics[width=\linewidth]{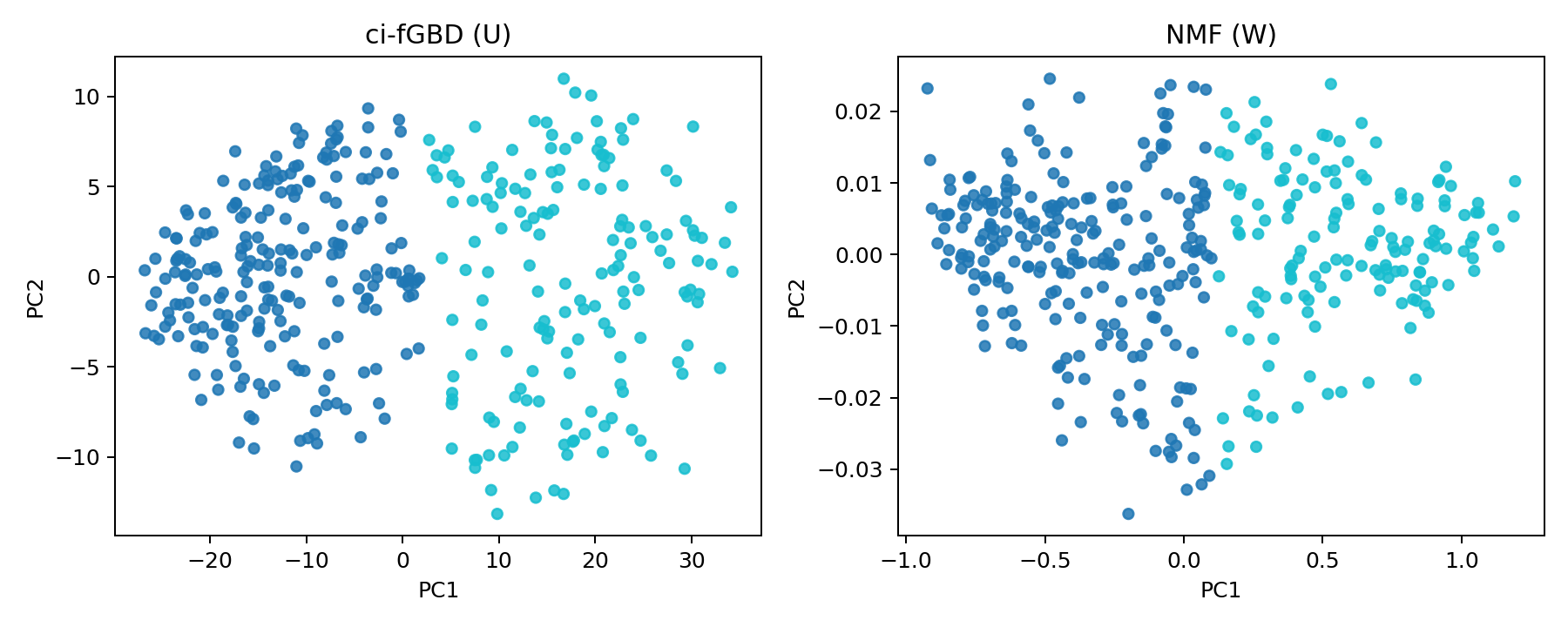}
    \caption{PC1--PC2 clustering projections of the AMR cohort using the ci-fGBD $U$ and baseline NMF $W$ representations. Both projections exhibit two broad colour-coded groups primarily separated along PC1, although differences are visible in within-group dispersion and the distribution of observations away from the densest point concentrations. The projected coordinate systems are method-specific; therefore, absolute coordinates, distances, and numerical axis ranges should not be compared directly. The colour assignments do not represent annotated antimicrobial-resistance phenotypes, species, strain identities, or validated antimicrobial-susceptibility classes.}
    \label{fig:cluster_amr}
    
\end{subfigure}
\caption{Complementary MTP and latent-space clustering analyses of the AMR cohort. Panel (a) presents the feature-contribution comparison between real AMR and LEvITAS-generated synthetic data, illustrating the retention of dominant genomic k-mer signatures and Cell Painting-derived morphological feature patterns. Panel (b) presents the PC1--PC2 projections of the ci-fGBD $U$ and baseline NMF $W$ representations, highlighting a broad two-group arrangement and differences in within-group dispersion and outlier structure. Together, these analyses provide complementary perspectives on genomic and morphology-associated feature fidelity and sample-level latent organization in synthetic AMR data. The observed patterns do not independently establish the preservation of genomic--phenotypic relationships, specific resistance mechanisms, or clinically validated antimicrobial-susceptibility characteristics.}
\label{fig:amr_exp_analysis}

\end{figure}

\paragraph*{AMR-Complementary Clustering Comparison:}
The AMR projections show two broad colour-coded groups in both ci-fGBD $U$ and 
NMF $W$, with the principal visible separation occurring along PC1 Figure~\ref{fig:amr_exp_analysis}(b). In the ci-fGBD panel, the two groups extend across a broad region of the plotted plane, 
with considerable within-group variation along PC2 and observations extending 
beyond the densest portions of each group. The NMF panel retains the broad 
separation on opposite sides of PC1 but displays a substantially more compact 
numerical range, while some observations extend away from the principal 
concentrations. These observations indicate that both decompositions expose a 
similar high-level two-group arrangement, although the visible dispersion and 
outlier structure differ. The large difference in axis ranges reflects the 
independent latent coordinate systems and should not be interpreted as 
evidence that ci-fGBD captures greater biological diversity or that NMF 
removes biologically meaningful variation. The figures do not annotate the 
groups with antimicrobial-resistance phenotypes, species, strain identities, 
or susceptibility labels. Consequently, the clustering comparison complements 
the AMR MTP by illustrating sample-level grouping of the integrated genomic 
and morphology-associated representations, but does not establish that the 
groups correspond to specific resistance mechanisms or validated antimicrobial-
susceptibility classes. 

Overall, the combined AMR MTP and clustering analyses 
Figure~\ref{fig:amr_exp_analysis} indicate that LEvITAS 
retains important genomic and morphology-derived feature representations, 
while the ci-fGBD and NMF projections both display two broad colour-coded 
groups with differences in within-group spread and outlier structure. The 
continued representation of dominant genomic k-mer signatures and Cell 
Painting-derived morphological descriptors suggests that aspects of the 
original multimodal feature organization are maintained, and the clustering 
projections provide a complementary sample-level view of this organization. 
Nevertheless, these similarities do not establish preservation of genomic--
phenotypic relationships, biologically meaningful resistance mechanisms, or 
clinically relevant antimicrobial-susceptibility characteristics. It is also 
important to determine whether reduced extreme latent contributions retain 
rare but meaningful resistance-associated patterns. Further validation through 
quantitative genomic--phenotypic comparisons, independent antibiotic-
susceptibility assessment on previously unseen real samples, and evaluation of 
rare resistance-associated phenotypes is therefore necessary. Taken together, 
these findings support the potential utility of LEvITAS for multimodal AMR 
data augmentation, microbial phenotype analysis, and future antimicrobial-
resistance research, while highlighting the importance of establishing 
genomic, phenotypic, and clinical validity before extending these findings to 
resistance prediction, susceptibility assessment, or translational infectious-
disease applications.

\newpage

\section*{Author Contributions}
- \textbf{Dr. Lokendra S. Thakur} conceived the idea, designed the study, developed 
LEvITAS and its all versions mathematical formulations with theorems and its 
proofs wherever required. Curated data, implemented experiments using ci-fGBD 
and its MTP visualization analysis for demonstrating practical feasibility of 
the LEvITAS, and all sections writing. \\
- \textbf{Dr. James Dale Berry} contributed to the biological and clinical analysis of 
the experimental results.  \\
- All authors reviewed and approved the final version.

\section*{Acknowledgments}
The locally used \textbf{GLENDA} image dataset corresponds to the
\emph{Gynecologic Laparoscopy Endometriosis Dataset (GLENDA)}, originally
described and released by Leibetseder et al. and distributed through the
GLENDA/ITEC resource \cite{Leibetseder2020GLENDA}. The locally retained GLENDA
observations were used for the image baseline mode.

The \textbf{USOVA3D} ultrasound imaging dataset was obtained from the
\textbf{USOVA3D} database, a publicly accessible resource containing annotated
three-dimensional ovarian ultrasound images developed for the validation of
follicle detection algorithms. The locally retained USOVA3D observations
were used for ovarian ultrasound phenotyping and the corresponding image
baseline analyses. The original database and its associated research
contributions are acknowledged through the publication by Poto\v{c}nik et al.
\cite{potocnik2020public}.

The \textbf{Neuromine} transcriptomic resource was obtained from the
\textbf{Answer ALS Data Portal (Neuromine)}, which provides clinical and
multi-omics data generated through the Answer ALS program and associated ALS
research initiatives \cite{Baxi2022AnswerALS}. The locally retained
\texttt{TRANSCRIPTOMICS.txt} representation was derived from this resource and
was used for both the static EHR and multiomics baseline modes.

The \textbf{ADNI} longitudinal clinical data were obtained from the
\textbf{Alzheimer's Disease Neuroimaging Initiative (ADNI)} database and used
for the temporal EHR baseline mode. ADNI data were used in accordance with the
applicable ADNI Data Use Agreement and publication policy. We acknowledge the
ADNI investigators, participating institutions, study participants, and
supporting organizations that contributed to the collection and dissemination
of these data \cite{weiner2013adni}.

The \textbf{AMR cell-painting} resource was used as a locally maintained
multimodal representation of an antimicrobial resistance Cell Painting
dataset described in the associated source study; the corresponding original
publication and source repository are acknowledged and cited separately in
the manuscript.
\\ \\
We are grateful to 
our in-house physician-scientist \textbf{Dr. Nino Dekanoidze} (MD in 
Medicine, David Tvildiani Medical University, GEORGIA) and early career researcher \textbf{Maya Chande} (Ph.D. student at Department of Mathematics, Massachusetts Institute of Technology, Massachusetts, Cambridge, United States) for their insightful suggestions. And

Dr. Russell Rockne Ph.D. (Division of Mathematical Oncology, City of Hope (CoH), Duarte, CA, United States.), Dr. Kamana Porwal Ph.D. (Department of Mathematics, IIT Delhi, India), 
Dr. Jason Flannick Ph.D. (Genetics and Genomics Division, Boston Children's Hospital and Broad Institute, MIT and Harvard), 
Dr. Eyal Y. Kimchi MD-Ph.D. (Department of Neurology, Northwestern University Feinberg School of Medicine, United States), 
Dr. Strajit Ghosh Ph.D. (Massachusetts General Hospital and Harvard Medical School, United Statea), Dr. Fernando Gómez-Baquero Ph.D. (Jacobs Technion–Cornell Institute, CornellTech; NSF Upstate New York Energy Storage Engine, USA) for their valuable guidance, and support. Which greatly contributed to this work.

\section*{Ethics Approval and Consent to Participate}
Not applicable.


\section*{Declaration of interests}
The authors declare no competing interests.


\section*{Funding}
The authors received no specific funding for this work.


\section*{Data sharing}
The locally used resources comprised the \textbf{GLENDA} image dataset,
the \textbf{Neuromine} transcriptomic resource, the \textbf{ADNI}
longitudinal clinical dataset, and the \textbf{AMR} cell-painting
resource. \textbf{GLENDA} (Gynecologic Laparoscopy Endometriosis
Dataset) was obtained from the official \textbf{GLENDA/ITEC} resource
\cite{Leibetseder2020GLENDA}. \textbf{Neuromine} data were obtained
through the \textbf{Answer ALS Data Portal (Neuromine)} and locally
retained in the \texttt{TRANSCRIPTOMICS.txt} representation
\cite{Baxi2022AnswerALS}. \textbf{ADNI} data were obtained through the
official \textbf{Alzheimer's Disease Neuroimaging Initiative (ADNI)}
database and used in accordance with the applicable ADNI Data Use
Agreement and publication requirements \cite{weiner2013adni}. The
\textbf{AMR cell-painting} resource was derived from publicly available
Cell Painting data distributed by the \textbf{JUMP Cell Painting
Consortium} through the \textbf{Cell Painting Gallery} hosted on the
\textbf{Registry of Open Data on AWS}; the corresponding JUMP resource
and source publication are cited in the manuscript
\cite{chandrasekaran2023jump}. The locally maintained multimodal NPZ
representation used for the AMR cell-painting benchmark was derived
from the relevant JUMP Cell Painting resource. These local files
represent locally stored and/or processed representations of the
respective source datasets and are not intended to denote additional
independent datasets.

Access to restricted or controlled datasets, including \textbf{ADNI}
and resources distributed through the \textbf{Answer ALS / Neuromine}
portal and \textbf{USOVA3D} is subject to registration, approval, data-use agreements, and other access requirements specified by the respective data custodians.

Public datasets should be obtained directly from their original or
official repositories described above. Derived feature representations
and benchmark-prepared data used by the LEvITAS pipeline are described
in this work; locally maintained raw files and processed representations
are not necessarily redistributed with the manuscript.
\\
Code supporting this study are available on reasonable request.

\section*{AI use statement}
ChatGPT was used to refine some of the sentences.

\newpage


\section*{References}
\renewcommand\refname{}

\bibliographystyle{unsrt}
\bibliography{references}


\end{document}